\documentclass{gtart}

\def\ifplaintex{\expandafter\ifx\csname documentclass\endcsname\relax}

\ifplaintex 
\else
\fi

\expandafter\ifx\csname beginpicture\endcsname\relax
\expandafter\ifx\csname documentclass\endcsname\relax
\input pictex \else
\input prepictex \input pictex \input postpictex \fi\fi

\def\gt{{\mathsurround=0pt\it $\cal G\mskip-2mu$eometry \&\ 
$\cal T\!\!$opology}}        %  journal title in recommended style

\def\gtp{{\mathsurround=0pt\it $\cal G\mskip-2mu$eometry \&\ 
$\cal T\!\!$opology $\cal P\!$ublications}}  % GT publications

\def\lognumber#1{\def\thelognumber{#1}}
\def\volumenumber#1{\def\thevolumenumber{#1}}
\def\papernumber#1{\def\thepapernumber{#1}}
\def\volumeyear#1{\def\thevolumeyear{#1}}

\def\pagenumbers#1#2{\def\startpage{#1}\def\finishpage{#2}}
\def\published#1{\def\publishdate{#1}}
\def\proposed#1{\def\theproposer{#1}}
\def\seconded#1{\def\theseconders{#1}}
\def\received#1{\def\receiveddate{#1}}
\def\revised#1{\def\reviseddate{#1}}
\def\accepted#1{\def\accepteddate{#1}}

\long\def\asciiabstract#1{\long\def\theasciiabstract{#1}}

\let\\\par\let\thelognumber\relax
\let\thevolumenumber\relax\let\thepapernumber\relax
\let\thevolumeyear\relax\let\thesamplenumber\relax\let\startpage\relax
\let\finishpage\relax\let\publishdate\relax\let\receiveddate\relax
\let\reviseddate\relax\let\accepteddate\relax\let\theasciititle\relax
\let\theasciiauthors\relax
\let\theasciiabstract\relax
\let\theasciiemail\relax\let\theshortauthors\relax\let\theshorttitle\relax

\long\def\maketitlep{   % start of definition of \maketitlep

\count0=\startpage

\gt\hfill      %   Journal title (top left) 
\beginpicture
\setcoordinatesystem units <0.33truein, 0.33truein> point at 2.2 0.9
\setplotsymbol ({$\cal G$})
\plotsymbolspacing=9truept
\circulararc 315 degrees from 0 1 center at 0 0
\setplotsymbol ({$\cal T$})
\circulararc 315 degrees from 1 -1 center at 1 0
\endpicture
\break
{\small\ifx\thesamplenumber\relax % sample?  
Volume \else Sample
\fi\thevolumenumber\ (\thevolumeyear)
\startpage--\finishpage\nl
Published: \publishdate}
\vglue 0.5truein plus 0.4fil minus 0.1truein

{\parskip=0pt\leftskip 0pt plus 1fil\def\\{\par\smallskip}{\ifplaintex\large
\else\Large\fi\bf\thetitle}\par\medskip}   

\vglue 0pt plus 0.1fil 

{\parskip=0pt\leftskip 0pt plus 1fil\def\\{\par}{\sc\theauthors}
\par\medskip}

\vglue 0pt plus 0.1fil 

{\small\parskip=0pt\let\newline\\
{\leftskip 0pt plus 1fil\def\\{\par}{\sl\theaddress}\par}
\expandafter\ifx\theemail\relax    % email address?
\relax\else\vglue 5pt plus 0.02fil minus 2pt\def\\{\stdspace{\rm 
and}\stdspace} 
\cl{Email:\stdspace\tt\theemail}\fi
\ifx\theurl\relax                  % URL given?
\relax\else\vglue 5pt plus 0.02fil minus 2pt\def\\{\stdspace{\rm 
and}\stdspace}
\cl{URL:\stdspace\tt\theurl}\fi\par}

\vglue 7pt plus 0.3fil minus 3pt

{\bf Abstract}
\vglue 5pt plus 0.1fil minus 2pt

\theabstract

\vglue 7pt plus 0.3fil minus 3pt

{\bf AMS Classification numbers}\quad Primary:\quad \theprimaryclass

Secondary:\quad \thesecondaryclass

\vglue 5pt plus 0.3fil minus 2pt

{\bf Keywords}\quad \thekeywords

\vglue 10pt plus 0.5fil minus 5pt

{\small  Proposed: \theproposer\hfill Received: \receiveddate\nl
Seconded: \theseconders\hfill 
\ifx\reviseddate\relax                         % paper revised?
Accepted: \accepteddate                        % no
\else
Revised: \reviseddate                          % yes
\fi}
\eject
}       %  end of definition of \maketitlep

\let\maketitlepage\maketitlep
\let\maketitle\maketitlepage

\font\phead=cmsl9 scaled 950
\font=cmsl9 scaled 1050
\font\pnum=cmbx10 scaled 913
\font\lnum=cmbx10 
\font\pfoot=cmsl9 scaled 950
\font=cmsl9 scaled 1050
\ifplaintex
\headline{\vbox to 0pt{\vskip -4.5mm\line{\small\phead\ifnum
\count0=\startpage ISSN 1364-0380 (on line)
1465-3060 (printed) \hfill {\pnum\folio}\else\ifodd\count0\def\\{ }% 
\ifx\theshorttitle\relax\thetitle\else\theshorttitle\fi\hfill{\pnum\folio}
\else\def\\{ and }{\pnum\folio}\hfill\ifx\theshortauthors\relax\theauthors
\else\theshortauthors\fi\fi\fi}\vss}}
\footline{\vbox to 0pt{\vglue 0mm\line{\small\pfoot\ifnum\count0=\startpage
\copyright\ \gtp\hfill\else
\gt, Volume \thevolumenumber\ (\thevolumeyear)\hfill\fi}\vss
}}
\else
\makeatletter
\def\@oddhead{{\small\lhead\ifnum\count0=\startpage ISSN 1364-0380 (on line)
1465-3060 (printed) \hfill {\lnum\number\count0}\else\ifodd\count0
\def\\{ }\ifx\theshorttitle\relax \thetitle \else\theshorttitle\fi\hfill
{\lnum\number\count0}\else\def\\{ and }{\lnum\number\count0}
\hfill\ifx\theshortauthors\relax 
\theauthors\else\theshortauthors\fi\fi\fi}}\def\@evenhead{\@oddhead}
\def\@oddfoot{\small\lfoot\ifnum\count0=\startpage\copyright\ \gtp\hfill\else
\gt, Volume \thevolumenumber\ (\thevolumeyear)\hfill\fi}
\def\@evenfoot{\@oddfoot}
\makeatother
\fi

\newwrite\gtoutfile
\long\gdef\makeheadfile{  %%% start of definition of \makeheadfile
{\def\\{, }\def\s{ }
\immediate\openout\gtoutfile head.xxx
\immediate\write\gtoutfile{Proxy-for: \ifx\theasciiauthors\relax
\theauthors\else\theasciiauthors\fi\s<\ifx\theasciiemail\relax\theemail\else\theasciiemail\fi>}
\immediate\write\gtoutfile{\noexpand\\}
\immediate\write\gtoutfile{Authors: \ifx\theasciiauthors\relax
\theauthors\else\theasciiauthors\fi}
{\def\\{ }\immediate\write\gtoutfile{Title: \ifx\theasciititle\relax
\thetitle\else\theasciititle\fi}}
\immediate\write\gtoutfile{Subj-class: GT or SG or MG etc}
\immediate\write\gtoutfile{MSC-class: \theprimaryclass\ifx\thesecondaryclass\relax\else, \thesecondaryclass\fi}
\immediate\write\gtoutfile{Journal-ref: Geom. Topol. \thevolumenumber
(\thevolumeyear) \startpage-\finishpage}
\immediate\write\gtoutfile{Comments: Published by Geometry and Topology at}
\immediate\write\gtoutfile{\s\s http://www.maths.warwick.ac.uk/gt/GTVol\thevolumenumber/paper\thepapernumber.abs.html}
\immediate\write\gtoutfile{\noexpand\\}
\immediate\write\gtoutfile{}
\ifx\theasciiabstract\relax
\immediate\write\gtoutfile{\theabstract}\else
\immediate\write\gtoutfile{\theasciiabstract}\fi
\immediate\write\gtoutfile{}
\immediate\write\gtoutfile{\noexpand\\}
\immediate\write\gtoutfile{}
\immediate\closeout\gtoutfile}}  %%% end of definition of \makeheadfile

\def\maketitlepage{\maketitlep\makeheadfile}
\let\maketitle\maketitlepage

\lognumber{301}
\volumenumber{8}\papernumber{5}\volumeyear{2004}
\pagenumbers{205}{275}
\received{22 January 2003}
\revised{12 February 2004}
\published{12 February 2004}
\accepted{17 December 2003}

\proposed{Benson Farb}
\seconded{Walter Neumann, Martin Bridson}

\usepackage{amsmath,amssymb,amsthm}

\newtheorem{thm}{Theorem}[section]
\newtheorem{prop}[thm]{Proposition}
\newtheorem{lem}[thm]{Lemma}
\newtheorem{cor}[thm]{Corollary}

\newtheorem{claim}{Claim}[thm] % Numbering of claims starts over
\theoremstyle{definition}
\newtheorem{defn}[thm]{Definition}
\newtheorem{rem}[thm]{Remark}
\newtheorem{conv}[thm]{Convention}

\newcommand{\abs}[1]{\left\lvert{#1}\right\rvert}
\newcommand{\norm}[1]{\left\lVert{#1}\right\rVert}
\renewcommand{\bar}[1]{\overline{#1}}
\newcommand{\boundary}{\partial}
\newcommand{\set}[2]{\{\,{#1} \mid {#2} \,\}}
\newcommand{\bigset}[2]{ \bigl\{ \, {#1} \bigm| {#2} \, \bigr\} }
\renewcommand{\emptyset}{\varnothing}
\renewcommand{\setminus}{-}

\newcommand{\field}[1]{\mathbb{#1}}
\newcommand{\Z}{\field{Z}}
\newcommand{\R}{\field{R}}
\newcommand{\N}{\field{N}}
\newcommand{\E}{\field{E}}
\newcommand{\Hyp}{\field{H}}

\newcommand{\inclusion}{\hookrightarrow}
\newcommand{\of}{\circ}
\renewcommand{\implies}{\Rightarrow}
\renewcommand{\hat}{\widehat}

\DeclareMathOperator{\CAT}{CAT}
\DeclareMathOperator{\Lk}{Lk}
\DeclareMathOperator{\Image}{Im}
\DeclareMathOperator{\Sh}{Sh}  % Shadow
\newcommand{\ball}[2]{B ( {#1}, {#2} )}   %  Ball{center}{radius}
\newcommand{\bigball}[2]{B \bigl( {#1}, {#2} \bigr)}
\newcommand{\nbd}[2]{\mathcal{N}_{#2}({#1})}  % Neighborhood{center}{radius}
\newcommand{\bignbd}[2]{\mathcal{N}_{#2} \bigl( {#1} \bigr)}
\newcommand{\Set}[1]{\mathcal{#1}}   % The name of a set.
\newcommand{\preflat}{P}
\newcommand{\interior}[1]{\mathring{#1}}

\usepackage{graphicx}

\newcommand{\drawthintriangle}{\begin{center} %
                               \includegraphics{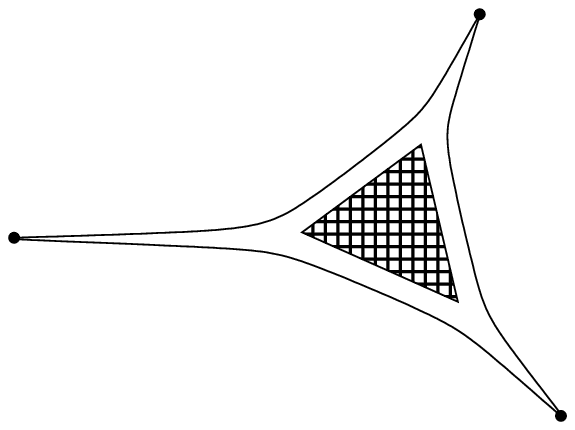}%
                               \end{center}}

\newcommand{\drawrftp}{\begin{center}%
                       \includegraphics{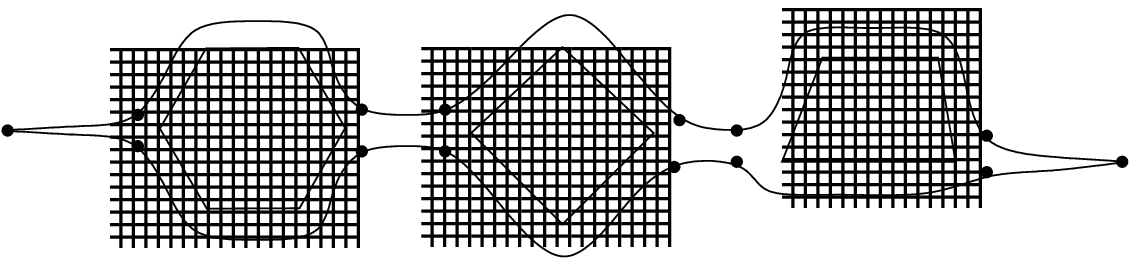}%
                       \end{center}}

\newcommand{\drawCATzero}{\begin{center}
                          \begin{picture}(0,0)%
\includegraphics{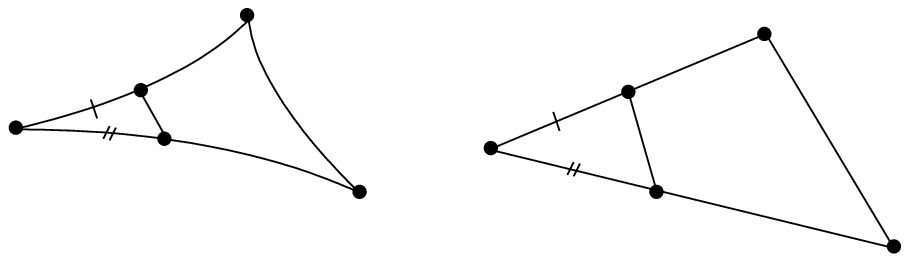}%
\end{picture}%
\setlength{\unitlength}{3947sp}%
\begingroup\makeatletter\ifx\SetFigFont\undefined
\def\x##1##2##3##4##5##6##7\relax{\def\x{##1##2##3##4##5##6}}%
\expandafter\x\fmtname xxxxxx\relax \def\y{splain}%
\ifx\x\y   % LaTeX or SliTeX?
\gdef\SetFigFont##1##2##3{%
  \ifnum ##1<17\tiny\else \ifnum ##1<20\small\else
  \ifnum ##1<24\normalsize\else \ifnum ##1<29\large\else
  \ifnum ##1<34\Large\else \ifnum ##1<41\LARGE\else
     \huge\fi\fi\fi\fi\fi\fi
  \csname ##3\endcsname}%
\else
\gdef\SetFigFont##1##2##3{\begingroup
  \count@##1\relax \ifnum 25<\count@\count@25\fi
  \def\x{\endgroup\@setsize\SetFigFont{##2pt}}%
  \expandafter\x
    \csname \romannumeral\the\count@ pt\expandafter\endcsname
    \csname @\romannumeral\the\count@ pt\endcsname
  \csname ##3\endcsname}%
\fi
\fi\endgroup
\begin{picture}(4463,1639)(98,-1274)
\put(1013,-721){\makebox(0,0)[lb]{\smash{\small{$y$}%
}}}
\put(827,-173){\makebox(0,0)[lb]{\smash{\small{$x$}%
}}}
\put(2011,-923){\makebox(0,0)[lb]{\smash{\small{$r$}%
}}}
\put(3158,-174){\makebox(0,0)[lb]{\smash{\small{$\bar{x}$}%
}}}
\put(3346,-961){\makebox(0,0)[lb]{\smash{\small{$\bar{y}$}%
}}}
\put(3894,119){\makebox(0,0)[lb]{\smash{\small{$\bar{q}$}%
}}}
\put(4530,-1216){\makebox(0,0)[lb]{\smash{\small{$\bar{r}$}%
}}}
\put(1404,209){\makebox(0,0)[lb]{\smash{\small{$q$}%
}}}
\put(2385,-593){\makebox(0,0)[lb]{\smash{\small{$\bar{p}$}%
}}}
\put( 98,-489){\makebox(0,0)[lb]{\smash{\small{$p$}%
}}}
\end{picture}
\end{center}}

\newcommand{\drawdiagram}{\begin{center}%
                          \includegraphics[width=4in]{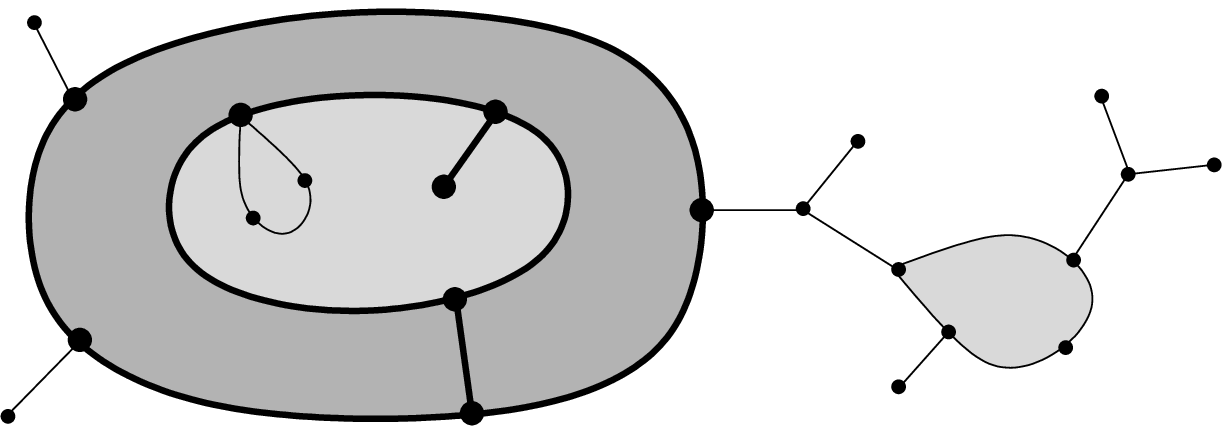}%
                          \end{center} }

\newcommand{\drawcancelable}{\begin{center}%
                             \begin{picture}(0,0)%
\includegraphics{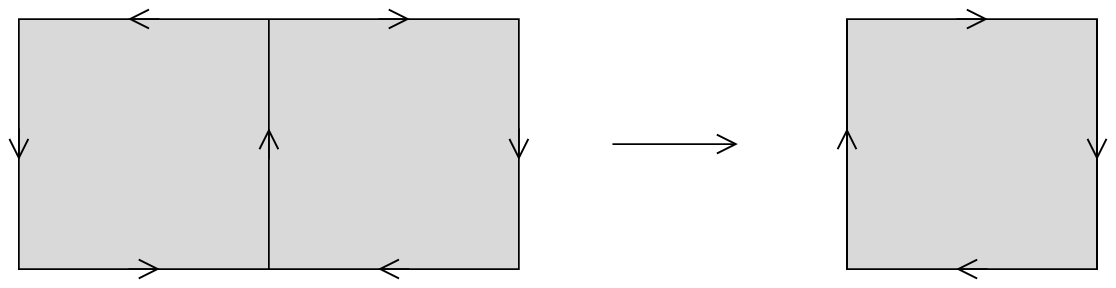}%
\end{picture}%
\setlength{\unitlength}{3947sp}%
\begingroup\makeatletter\ifx\SetFigFont\undefined
\def\x##1##2##3##4##5##6##7\relax{\def\x{##1##2##3##4##5##6}}%
\expandafter\x\fmtname xxxxxx\relax \def\y{splain}%
\ifx\x\y   % LaTeX or SliTeX?
\gdef\SetFigFont##1##2##3{%
  \ifnum ##1<17\tiny\else \ifnum ##1<20\small\else
  \ifnum ##1<24\normalsize\else \ifnum ##1<29\large\else
  \ifnum ##1<34\Large\else \ifnum ##1<41\LARGE\else
     \huge\fi\fi\fi\fi\fi\fi
  \csname ##3\endcsname}%
\else
\gdef\SetFigFont##1##2##3{\begingroup
  \count@##1\relax \ifnum 25<\count@\count@25\fi
  \def\x{\endgroup\@setsize\SetFigFont{##2pt}}%
  \expandafter\x
    \csname \romannumeral\the\count@ pt\expandafter\endcsname
    \csname @\romannumeral\the\count@ pt\endcsname
  \csname ##3\endcsname}%
\fi
\fi\endgroup
\begin{picture}(5550,1756)(376,-1835)
\put(4276,-1036){\makebox(0,0)[lb]{\smash{\small{$a$}%
}}}
\put(1876,-1036){\makebox(0,0)[lb]{\smash{\small{$a$}%
}}}
\put(5101,-211){\makebox(0,0)[lb]{\smash{\small{$b$}%
}}}
\put(2326,-211){\makebox(0,0)[lb]{\smash{\small{$b$}%
}}}
\put(1126,-211){\makebox(0,0)[lb]{\smash{\small{$b$}%
}}}
\put(5926,-1036){\makebox(0,0)[lb]{\smash{\small{$c$}%
}}}
\put(3076,-1036){\makebox(0,0)[lb]{\smash{\small{$c$}%
}}}
\put(376,-1036){\makebox(0,0)[lb]{\smash{\small{$c$}%
}}}
\put(5101,-1786){\makebox(0,0)[lb]{\smash{\small{$d$}%
}}}
\put(2326,-1786){\makebox(0,0)[lb]{\smash{\small{$d$}%
}}}
\put(1126,-1786){\makebox(0,0)[lb]{\smash{\small{$d$}%
}}}
\end{picture}
\end{center} }

\newcommand{\drawBranchedCover}{\begin{center}%
                                \begin{picture}(0,0)%
\includegraphics{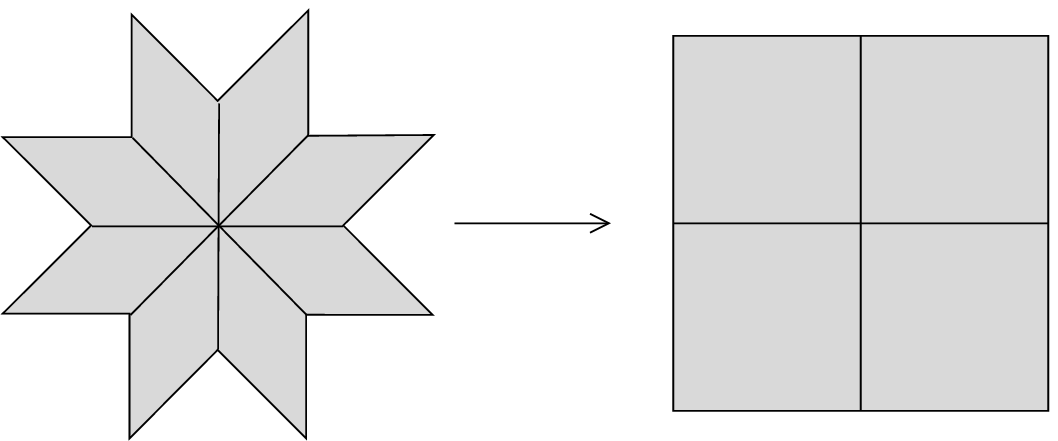}%
\end{picture}%
\setlength{\unitlength}{3947sp}%
\begingroup\makeatletter\ifx\SetFigFont\undefined
\def\x##1##2##3##4##5##6##7\relax{\def\x{##1##2##3##4##5##6}}%
\expandafter\x\fmtname xxxxxx\relax \def\y{splain}%
\ifx\x\y   % LaTeX or SliTeX?
\gdef\SetFigFont##1##2##3{%
  \ifnum ##1<17\tiny\else \ifnum ##1<20\small\else
  \ifnum ##1<24\normalsize\else \ifnum ##1<29\large\else
  \ifnum ##1<34\Large\else \ifnum ##1<41\LARGE\else
     \huge\fi\fi\fi\fi\fi\fi
  \csname ##3\endcsname}%
\else
\gdef\SetFigFont##1##2##3{\begingroup
  \count@##1\relax \ifnum 25<\count@\count@25\fi
  \def\x{\endgroup\@setsize\SetFigFont{##2pt}}%
  \expandafter\x
    \csname \romannumeral\the\count@ pt\expandafter\endcsname
    \csname @\romannumeral\the\count@ pt\endcsname
  \csname ##3\endcsname}%
\fi
\fi\endgroup
\begin{picture}(5043,2079)(370,-1706)
\put(3991,-286){\makebox(0,0)[lb]{\smash{\small{$A$}%
}}}
\put(4869,-286){\makebox(0,0)[lb]{\smash{\small{$B$}%
}}}
\put(4028,-1171){\makebox(0,0)[lb]{\smash{\small{$D$}%
}}}
\put(4868,-1141){\makebox(0,0)[lb]{\smash{\small{$C$}%
}}}
\put(864,-526){\makebox(0,0)[lb]{\smash{\small{$A$}%
}}}
\put(1156,-211){\makebox(0,0)[lb]{\smash{\small{$B$}%
}}}
\put(1561,-211){\makebox(0,0)[lb]{\smash{\small{$C$}%
}}}
\put(1853,-518){\makebox(0,0)[lb]{\smash{\small{$D$}%
}}}
\put(1884,-938){\makebox(0,0)[lb]{\smash{\small{$A$}%
}}}
\put(1576,-1231){\makebox(0,0)[lb]{\smash{\small{$B$}%
}}}
\put(1163,-1216){\makebox(0,0)[lb]{\smash{\small{$C$}%
}}}
\put(864,-916){\makebox(0,0)[lb]{\smash{\small{$D$}%
}}}
\end{picture}
\end{center} }

\newcommand{\drawcurvature}{\begin{center}%
                            \begin{picture}(0,0)%
\includegraphics{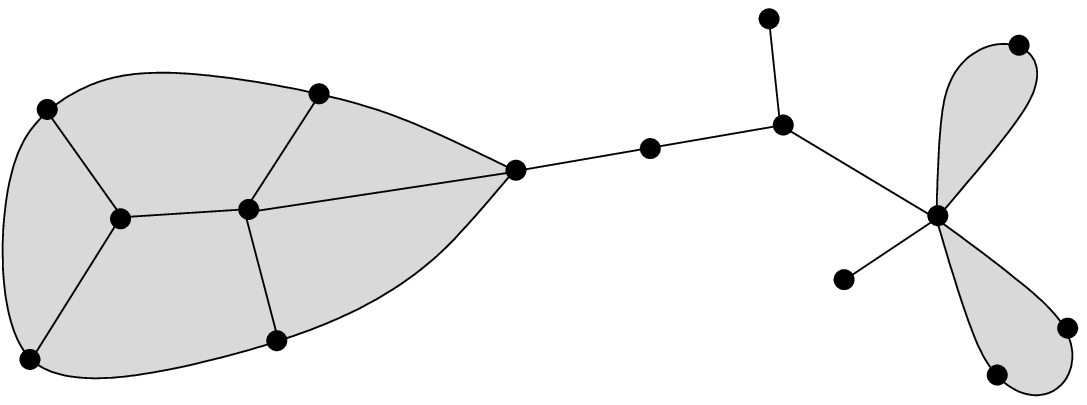}%
\end{picture}%
\setlength{\unitlength}{3947sp}%
\begingroup\makeatletter\ifx\SetFigFont\undefined
\def\x##1##2##3##4##5##6##7\relax{\def\x{##1##2##3##4##5##6}}%
\expandafter\x\fmtname xxxxxx\relax \def\y{splain}%
\ifx\x\y   % LaTeX or SliTeX?
\gdef\SetFigFont##1##2##3{%
  \ifnum ##1<17\tiny\else \ifnum ##1<20\small\else
  \ifnum ##1<24\normalsize\else \ifnum ##1<29\large\else
  \ifnum ##1<34\Large\else \ifnum ##1<41\LARGE\else
     \huge\fi\fi\fi\fi\fi\fi
  \csname ##3\endcsname}%
\else
\gdef\SetFigFont##1##2##3{\begingroup
  \count@##1\relax \ifnum 25<\count@\count@25\fi
  \def\x{\endgroup\@setsize\SetFigFont{##2pt}}%
  \expandafter\x
    \csname \romannumeral\the\count@ pt\expandafter\endcsname
    \csname @\romannumeral\the\count@ pt\endcsname
  \csname ##3\endcsname}%
\fi
\fi\endgroup
\begin{picture}(5221,2187)(420,-1506)
\put(4756,-330){\makebox(0,0)[lb]{\smash{\small{$\pi$}%
}}}
\put(4816,-713){\makebox(0,0)[lb]{\smash{\small{$\pi$}%
}}}
\put(4666,-548){\makebox(0,0)[lb]{\smash{\small{$\pi$}%
}}}
\put(4276,-953){\makebox(0,0)[lb]{\smash{\small{$\pi$}%
}}}
\put(5064,-1448){\makebox(0,0)[lb]{\smash{\small{$\pi$}%
}}}
\put(5641,-938){\makebox(0,0)[lb]{\smash{\small{$\pi$}%
}}}
\put(5311,457){\makebox(0,0)[lb]{\smash{\small{$\pi$}%
}}}
\put(1719,-1305){\makebox(0,0)[lb]{\smash{\small{$\pi$}%
}}}
\put(429,-1373){\makebox(0,0)[lb]{\smash{\small{$\pi$}%
}}}
\put(496,150){\makebox(0,0)[lb]{\smash{\small{$\pi$}%
}}}
\put(1944,232){\makebox(0,0)[lb]{\smash{\small{$\pi$}%
}}}
\put(2836,-120){\makebox(0,0)[lb]{\smash{\small{$\pi$}%
}}}
\put(2866,-495){\makebox(0,0)[lb]{\smash{\small{$\pi$}%
}}}
\put(3511,-383){\makebox(0,0)[lb]{\smash{\small{$\pi$}%
}}}
\put(3481,-23){\makebox(0,0)[lb]{\smash{\small{$\pi$}%
}}}
\put(3999, 22){\makebox(0,0)[lb]{\smash{\small{$\pi$}%
}}}
\put(4269, -8){\makebox(0,0)[lb]{\smash{\small{$\pi$}%
}}}
\put(4111,-240){\makebox(0,0)[lb]{\smash{\small{$\pi$}%
}}}
\put(3954,525){\makebox(0,0)[lb]{\smash{\small{$\pi$}%
}}}
\put(5034,-533){\makebox(0,0)[lb]{\smash{\small{$\pi$}%
}}}
\end{picture}
\end{center} }

\newcommand{\drawBadTriangle}{\begin{center}
\begin{picture}(0,0)%
\includegraphics{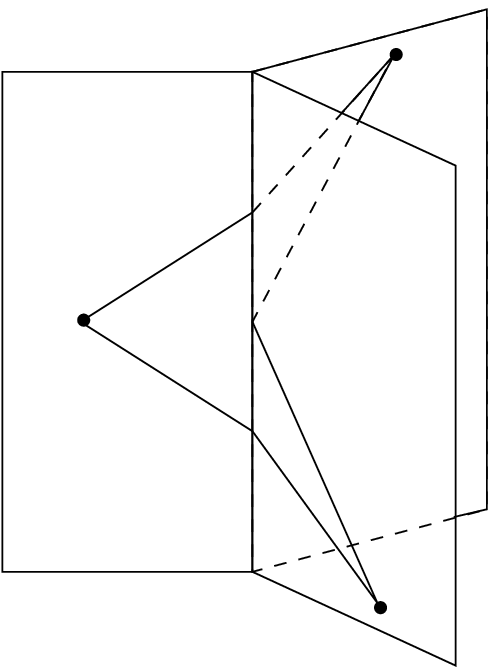}%
\end{picture}%
\setlength{\unitlength}{3947sp}%
\begingroup\makeatletter\ifx\SetFigFont\undefined
\def\x##1##2##3##4##5##6##7\relax{\def\x{##1##2##3##4##5##6}}%
\expandafter\x\fmtname xxxxxx\relax \def\y{splain}%
\ifx\x\y   % LaTeX or SliTeX?
\gdef\SetFigFont##1##2##3{%
  \ifnum ##1<17\tiny\else \ifnum ##1<20\small\else
  \ifnum ##1<24\normalsize\else \ifnum ##1<29\large\else
  \ifnum ##1<34\Large\else \ifnum ##1<41\LARGE\else
     \huge\fi\fi\fi\fi\fi\fi
  \csname ##3\endcsname}%
\else
\gdef\SetFigFont##1##2##3{\begingroup
  \count@##1\relax \ifnum 25<\count@\count@25\fi
  \def\x{\endgroup\@setsize\SetFigFont{##2pt}}%
  \expandafter\x
    \csname \romannumeral\the\count@ pt\expandafter\endcsname
    \csname @\romannumeral\the\count@ pt\endcsname
  \csname ##3\endcsname}%
\fi
\fi\endgroup
\begin{picture}(2349,3174)(1189,-3223)
\put(1373,-1614){\makebox(0,0)[lb]{\smash{\small{$a$}%
}}}
\put(3106,-2986){\makebox(0,0)[lb]{\smash{\small{$b$}%
}}}
\put(3166,-346){\makebox(0,0)[lb]{\smash{\small{$c$}%
}}}
\end{picture}
\end{center}}

\newcommand{\drawSimpleRuffles}{\begin{center}%
\begin{picture}(0,0)%
\includegraphics{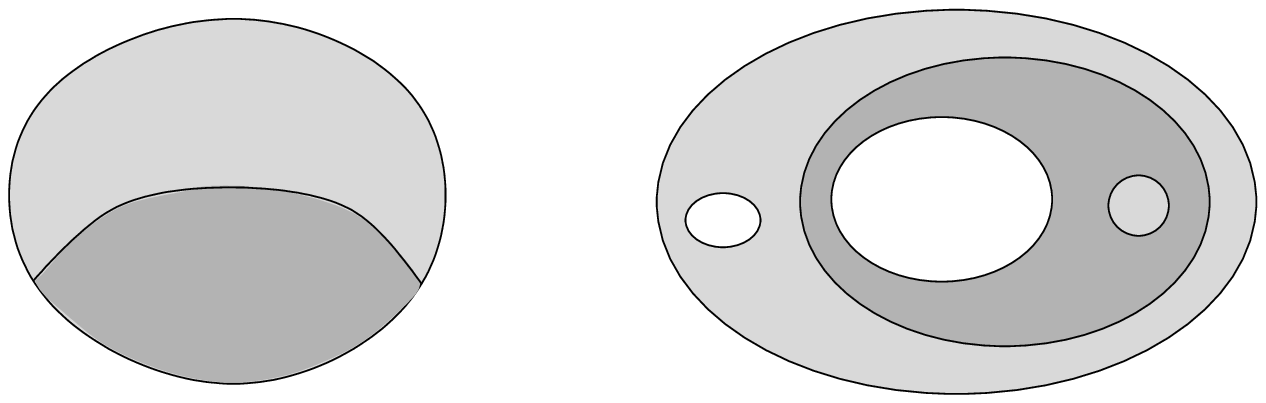}%
\end{picture}%
\setlength{\unitlength}{3947sp}%
\begingroup\makeatletter\ifx\SetFigFont\undefined
\def\x##1##2##3##4##5##6##7\relax{\def\x{##1##2##3##4##5##6}}%
\expandafter\x\fmtname xxxxxx\relax \def\y{splain}%
\ifx\x\y   % LaTeX or SliTeX?
\gdef\SetFigFont##1##2##3{%
  \ifnum ##1<17\tiny\else \ifnum ##1<20\small\else
  \ifnum ##1<24\normalsize\else \ifnum ##1<29\large\else
  \ifnum ##1<34\Large\else \ifnum ##1<41\LARGE\else
     \huge\fi\fi\fi\fi\fi\fi
  \csname ##3\endcsname}%
\else
\gdef\SetFigFont##1##2##3{\begingroup
  \count@##1\relax \ifnum 25<\count@\count@25\fi
  \def\x{\endgroup\@setsize\SetFigFont{##2pt}}%
  \expandafter\x
    \csname \romannumeral\the\count@ pt\expandafter\endcsname
    \csname @\romannumeral\the\count@ pt\endcsname
  \csname ##3\endcsname}%
\fi
\fi\endgroup
\begin{picture}(6048,1941)(451,-1558)
\put(2273,149){\makebox(0,0)[lb]{\smash{\small{$\beta$}%
}}}
\put(4815,-908){\makebox(0,0)[lb]{\smash{\small{$C$}%
}}}
\put(451,239){\makebox(0,0)[lb]{\smash{\small{(a)}%
}}}
\put(3526,239){\makebox(0,0)[lb]{\smash{\small{(b)}%
}}}
\put(2079,-1500){\makebox(0,0)[lb]{\smash{\small{$\alpha$}%
}}}
\put(1756,-458){\makebox(0,0)[lb]{\smash{\small{$\beta'$}%
}}}
\end{picture}
\end{center} }

\newcommand{\drawconcatenations}{\begin{center}%
\begin{picture}(0,0)%
\includegraphics{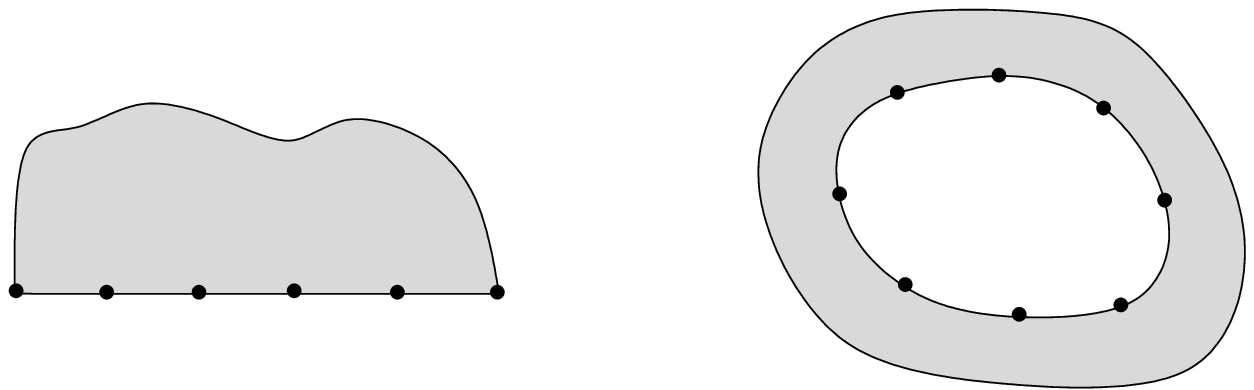}%
\end{picture}%
\setlength{\unitlength}{3947sp}%
\begingroup\makeatletter\ifx\SetFigFont\undefined
\def\x##1##2##3##4##5##6##7\relax{\def\x{##1##2##3##4##5##6}}%
\expandafter\x\fmtname xxxxxx\relax \def\y{splain}%
\ifx\x\y   % LaTeX or SliTeX?
\gdef\SetFigFont##1##2##3{%
  \ifnum ##1<17\tiny\else \ifnum ##1<20\small\else
  \ifnum ##1<24\normalsize\else \ifnum ##1<29\large\else
  \ifnum ##1<34\Large\else \ifnum ##1<41\LARGE\else
     \huge\fi\fi\fi\fi\fi\fi
  \csname ##3\endcsname}%
\else
\gdef\SetFigFont##1##2##3{\begingroup
  \count@##1\relax \ifnum 25<\count@\count@25\fi
  \def\x{\endgroup\@setsize\SetFigFont{##2pt}}%
  \expandafter\x
    \csname \romannumeral\the\count@ pt\expandafter\endcsname
    \csname @\romannumeral\the\count@ pt\endcsname
  \csname ##3\endcsname}%
\fi
\fi\endgroup
\begin{picture}(5992,1838)(226,-1325)
\put(5190,-22){\makebox(0,0)[lb]{\smash{\small{$\alpha_1$}%
}}}
\put(5506,-278){\makebox(0,0)[lb]{\smash{\small{$\alpha_2$}%
}}}
\put(5626,-682){\makebox(0,0)[lb]{\smash{\small{$\alpha_3$}%
}}}
\put(5296,-885){\makebox(0,0)[lb]{\smash{\small{$\alpha_4$}%
}}}
\put(4344,-211){\makebox(0,0)[lb]{\smash{\small{$\alpha_n$}%
}}}
\put(2311,-1073){\makebox(0,0)[lb]{\smash{\small{$\alpha_n$}%
}}}
\put(1336,-1065){\makebox(0,0)[lb]{\smash{\small{$\alpha_3$}%
}}}
\put(894,-1058){\makebox(0,0)[lb]{\smash{\small{$\alpha_2$}%
}}}
\put(450,-1050){\makebox(0,0)[lb]{\smash{\small{$\alpha_1$}%
}}}
\put(226,239){\makebox(0,0)[lb]{\smash{\small{(a)}%
}}}
\put(3676,239){\makebox(0,0)[lb]{\smash{\small{(b)}%
}}}
\put(1201, 89){\makebox(0,0)[lb]{\smash{\small{$\beta$}%
}}}
\put(4726,-61){\makebox(0,0)[lb]{\smash{\small{$\beta$}%
}}}
\end{picture}
\end{center} }

\newcommand{\drawruffles}{\begin{center}
\begin{picture}(0,0)%
\includegraphics{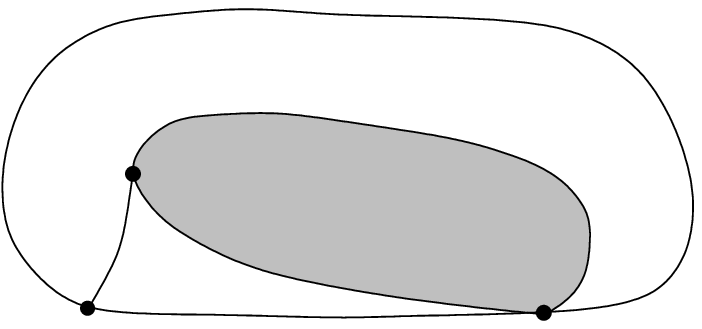}%
\end{picture}%
\setlength{\unitlength}{3947sp}%
\begingroup\makeatletter\ifx\SetFigFont\undefined
\def\x##1##2##3##4##5##6##7\relax{\def\x{##1##2##3##4##5##6}}%
\expandafter\x\fmtname xxxxxx\relax \def\y{splain}%
\ifx\x\y   % LaTeX or SliTeX?
\gdef\SetFigFont##1##2##3{%
  \ifnum ##1<17\tiny\else \ifnum ##1<20\small\else
  \ifnum ##1<24\normalsize\else \ifnum ##1<29\large\else
  \ifnum ##1<34\Large\else \ifnum ##1<41\LARGE\else
     \huge\fi\fi\fi\fi\fi\fi
  \csname ##3\endcsname}%
\else
\gdef\SetFigFont##1##2##3{\begingroup
  \count@##1\relax \ifnum 25<\count@\count@25\fi
  \def\x{\endgroup\@setsize\SetFigFont{##2pt}}%
  \expandafter\x
    \csname \romannumeral\the\count@ pt\expandafter\endcsname
    \csname @\romannumeral\the\count@ pt\endcsname
  \csname ##3\endcsname}%
\fi
\fi\endgroup
\begin{picture}(3339,1738)(832,-1222)
\put(1921,-1164){\makebox(0,0)[lb]{\smash{\small{$\alpha$}%
}}}
\put(2041,-698){\makebox(0,0)[lb]{\smash{\small{$\alpha'$}%
}}}
\put(1209,-639){\makebox(0,0)[lb]{\smash{\small{$\beta$}%
}}}
\put(2663,284){\makebox(0,0)[lb]{\smash{\small{$\gamma$}%
}}}
\put(2948,-316){\makebox(0,0)[lb]{\smash{\small{$\gamma'$}%
}}}
\put(2431,-428){\makebox(0,0)[lb]{\smash{\small{$D'$}%
}}}
\end{picture}
\end{center}}

\newcommand{\drawbroom}{\begin{center}
\begin{picture}(0,0)%
\includegraphics{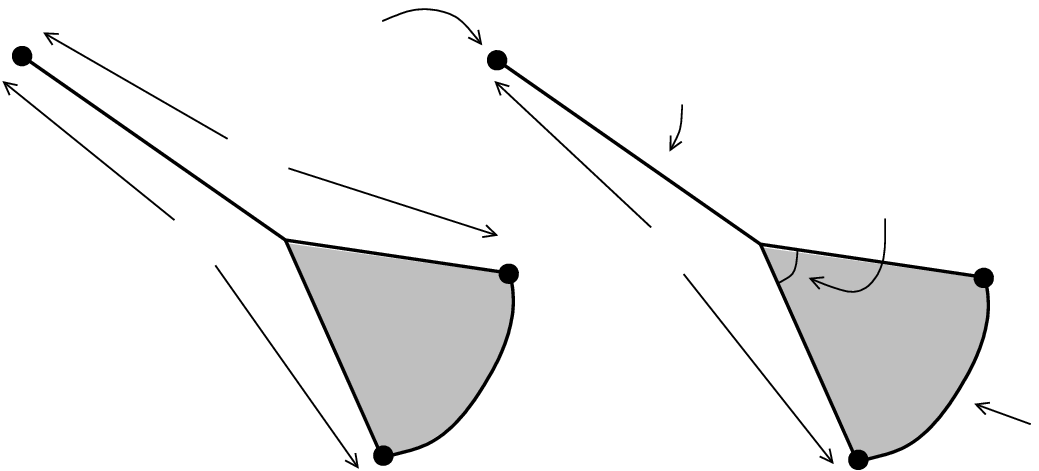}%
\end{picture}%
\setlength{\unitlength}{3947sp}%
\begingroup\makeatletter\ifx\SetFigFont\undefined
\def\x##1##2##3##4##5##6##7\relax{\def\x{##1##2##3##4##5##6}}%
\expandafter\x\fmtname xxxxxx\relax \def\y{splain}%
\ifx\x\y   % LaTeX or SliTeX?
\gdef\SetFigFont##1##2##3{%
  \ifnum ##1<17\tiny\else \ifnum ##1<20\small\else
  \ifnum ##1<24\normalsize\else \ifnum ##1<29\large\else
  \ifnum ##1<34\Large\else \ifnum ##1<41\LARGE\else
     \huge\fi\fi\fi\fi\fi\fi
  \csname ##3\endcsname}%
\else
\gdef\SetFigFont##1##2##3{\begingroup
  \count@##1\relax \ifnum 25<\count@\count@25\fi
  \def\x{\endgroup\@setsize\SetFigFont{##2pt}}%
  \expandafter\x
    \csname \romannumeral\the\count@ pt\expandafter\endcsname
    \csname @\romannumeral\the\count@ pt\endcsname
  \csname ##3\endcsname}%
\fi
\fi\endgroup
\begin{picture}(5270,2229)(131,-1768)
\put(2513,-1515){\makebox(0,0)[lb]{\smash{\small{$\gamma$}%
}}}
\put(1290,-315){\makebox(0,0)[lb]{\smash{\small{$\beta$}%
}}}
\put(1005,-728){\makebox(0,0)[lb]{\smash{\small{$\alpha$}%
}}}
\put(5401,-1261){\makebox(0,0)[lb]{\smash{\small{ }%
}}}
\put(3157, 50){\makebox(0,0)[lb]{\smash{\small{handle}%
}}}
\put(3953,-287){\makebox(0,0)[lb]{\smash{\small{branching}%
}}}
\put(4088,-459){\makebox(0,0)[lb]{\smash{\small{angle}%
}}}
\put(1756,276){\makebox(0,0)[lb]{\smash{\small{tip}%
}}}
\put(5086,-1490){\makebox(0,0)[lb]{\smash{\small{outer}%
}}}
\put(5131,-1692){\makebox(0,0)[lb]{\smash{\small{path}%
}}}
\put(3150,-759){\makebox(0,0)[lb]{\smash{\small{height}%
}}}
\end{picture}
\end{center}}

\newcommand{\drawquadraticdiv}{\begin{center}
\begin{picture}(0,0)%
\includegraphics{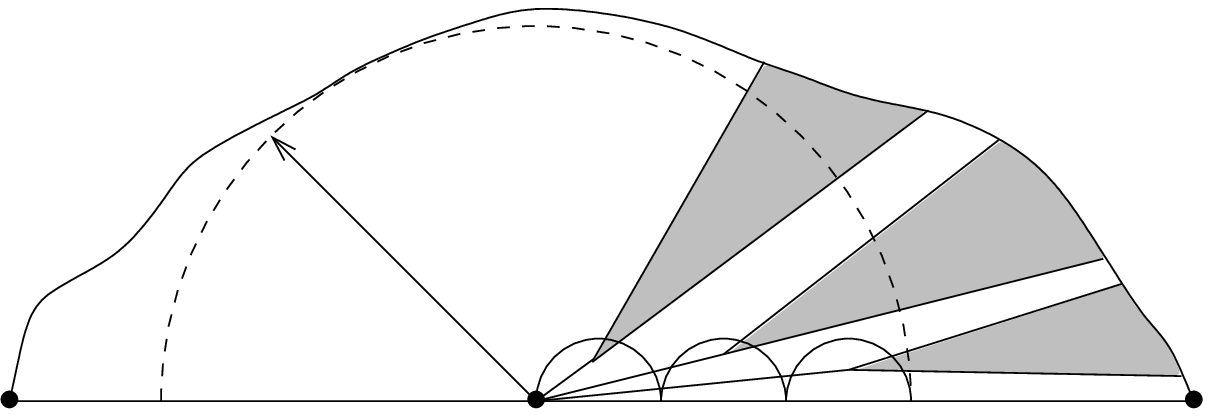}%
\end{picture}%
\setlength{\unitlength}{3947sp}%
\begingroup\makeatletter\ifx\SetFigFont\undefined
\def\x##1##2##3##4##5##6##7\relax{\def\x{##1##2##3##4##5##6}}%
\expandafter\x\fmtname xxxxxx\relax \def\y{splain}%
\ifx\x\y   % LaTeX or SliTeX?
\gdef\SetFigFont##1##2##3{%
  \ifnum ##1<17\tiny\else \ifnum ##1<20\small\else
  \ifnum ##1<24\normalsize\else \ifnum ##1<29\large\else
  \ifnum ##1<34\Large\else \ifnum ##1<41\LARGE\else
     \huge\fi\fi\fi\fi\fi\fi
  \csname ##3\endcsname}%
\else
\gdef\SetFigFont##1##2##3{\begingroup
  \count@##1\relax \ifnum 25<\count@\count@25\fi
  \def\x{\endgroup\@setsize\SetFigFont{##2pt}}%
  \expandafter\x
    \csname \romannumeral\the\count@ pt\expandafter\endcsname
    \csname @\romannumeral\the\count@ pt\endcsname
  \csname ##3\endcsname}%
\fi
\fi\endgroup
\begin{picture}(5777,2390)(202,-1867)
\put(2731,-1809){\makebox(0,0)[lb]{\smash{\small{$p$}%
}}}
\put(1576,-1793){\makebox(0,0)[lb]{\smash{\small{$\gamma$}%
}}}
\put(2026,-736){\makebox(0,0)[lb]{\smash{\small{$r$}%
}}}
\put(2551,367){\makebox(0,0)[lb]{\smash{\small{$\alpha$}%
}}}
\end{picture}
\end{center}}

\newcommand{\drawshadows}{\begin{center}
\begin{picture}(0,0)%
\includegraphics{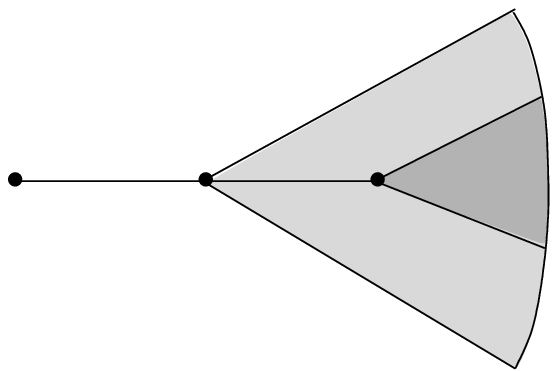}%
\end{picture}%
\setlength{\unitlength}{3947sp}%
\begingroup\makeatletter\ifx\SetFigFont\undefined
\def\x##1##2##3##4##5##6##7\relax{\def\x{##1##2##3##4##5##6}}%
\expandafter\x\fmtname xxxxxx\relax \def\y{splain}%
\ifx\x\y   % LaTeX or SliTeX?
\gdef\SetFigFont##1##2##3{%
  \ifnum ##1<17\tiny\else \ifnum ##1<20\small\else
  \ifnum ##1<24\normalsize\else \ifnum ##1<29\large\else
  \ifnum ##1<34\Large\else \ifnum ##1<41\LARGE\else
     \huge\fi\fi\fi\fi\fi\fi
  \csname ##3\endcsname}%
\else
\gdef\SetFigFont##1##2##3{\begingroup
  \count@##1\relax \ifnum 25<\count@\count@25\fi
  \def\x{\endgroup\@setsize\SetFigFont{##2pt}}%
  \expandafter\x
    \csname \romannumeral\the\count@ pt\expandafter\endcsname
    \csname @\romannumeral\the\count@ pt\endcsname
  \csname ##3\endcsname}%
\fi
\fi\endgroup
\begin{picture}(2789,1749)(384,-1273)
\put(384,-398){\makebox(0,0)[lb]{\smash{\small{$p$}%
}}}
\put(1411,-579){\makebox(0,0)[lb]{\smash{\small{$y_j$}%
}}}
\put(2273,-563){\makebox(0,0)[lb]{\smash{\small{$y_i$}%
}}}
\put(2439,-68){\makebox(0,0)[lb]{\smash{\small{$C_j$}%
}}}
\put(2806,-383){\makebox(0,0)[lb]{\smash{\small{$B_i$}%
}}}
\put(2551,-870){\makebox(0,0)[lb]{\smash{\small{$C'_j$}%
}}}
\end{picture}
\end{center}}

\newcommand{\drawruffledftp}{\begin{center}
\begin{picture}(0,0)%
\includegraphics{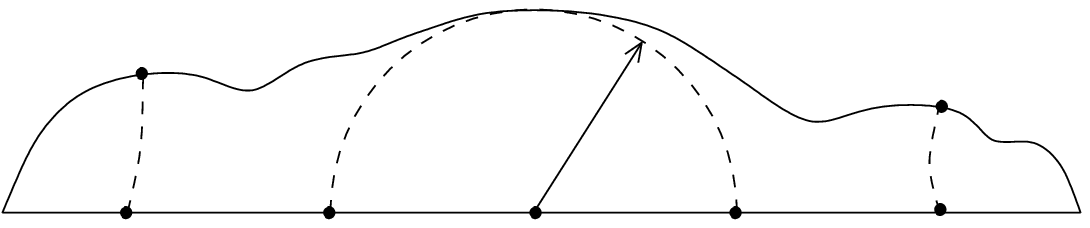}%
\end{picture}%
\setlength{\unitlength}{3947sp}%
\begingroup\makeatletter\ifx\SetFigFont\undefined
\def\x##1##2##3##4##5##6##7\relax{\def\x{##1##2##3##4##5##6}}%
\expandafter\x\fmtname xxxxxx\relax \def\y{splain}%
\ifx\x\y   % LaTeX or SliTeX?
\gdef\SetFigFont##1##2##3{%
  \ifnum ##1<17\tiny\else \ifnum ##1<20\small\else
  \ifnum ##1<24\normalsize\else \ifnum ##1<29\large\else
  \ifnum ##1<34\Large\else \ifnum ##1<41\LARGE\else
     \huge\fi\fi\fi\fi\fi\fi
  \csname ##3\endcsname}%
\else
\gdef\SetFigFont##1##2##3{\begingroup
  \count@##1\relax \ifnum 25<\count@\count@25\fi
  \def\x{\endgroup\@setsize\SetFigFont{##2pt}}%
  \expandafter\x
    \csname \romannumeral\the\count@ pt\expandafter\endcsname
    \csname @\romannumeral\the\count@ pt\endcsname
  \csname ##3\endcsname}%
\fi
\fi\endgroup
\begin{picture}(5199,1312)(214,-1889)
\put(818,-789){\makebox(0,0)[lb]{\smash{\small{$u$}%
}}}
\put(931,-1305){\makebox(0,0)[lb]{\smash{\small{$\le r$}%
}}}
\put(758,-1809){\makebox(0,0)[lb]{\smash{\small{$y'$}%
}}}
\put(1749,-1831){\makebox(0,0)[lb]{\smash{\small{$y$}%
}}}
\put(2739,-1816){\makebox(0,0)[lb]{\smash{\small{$x$}%
}}}
\put(3705,-1815){\makebox(0,0)[lb]{\smash{\small{$z$}%
}}}
\put(4688,-946){\makebox(0,0)[lb]{\smash{\small{$v$}%
}}}
\put(4710,-1381){\makebox(0,0)[lb]{\smash{\small{$\le r$}%
}}}
\put(4673,-1786){\makebox(0,0)[lb]{\smash{\small{$z'$}%
}}}
\put(3083,-1246){\makebox(0,0)[lb]{\smash{\small{$r$}%
}}}
\end{picture}
\end{center}}

\newcommand{\drawpreflats}{\begin{center}%
                           \includegraphics{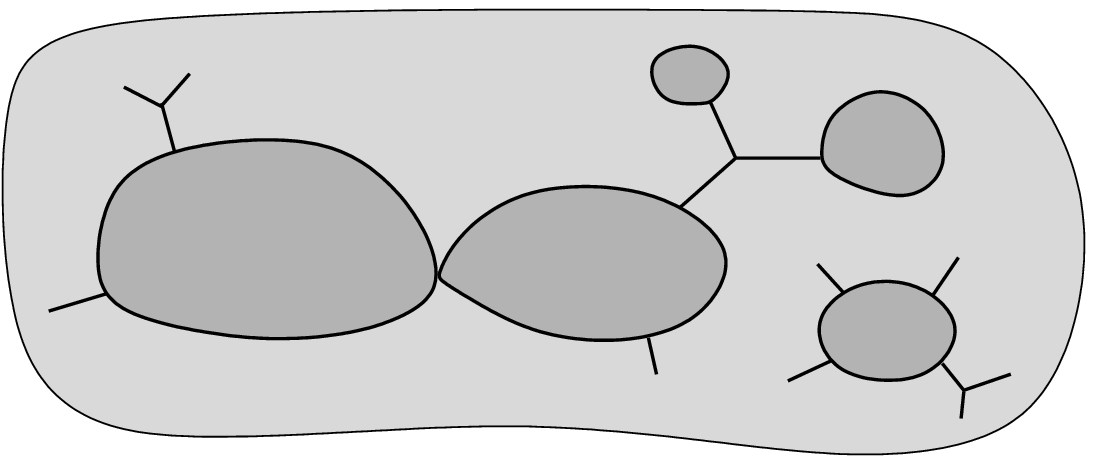}%
                           \end{center} }

\newcommand{\drawflatdisc}{\begin{center}
\begin{picture}(0,0)%
\includegraphics{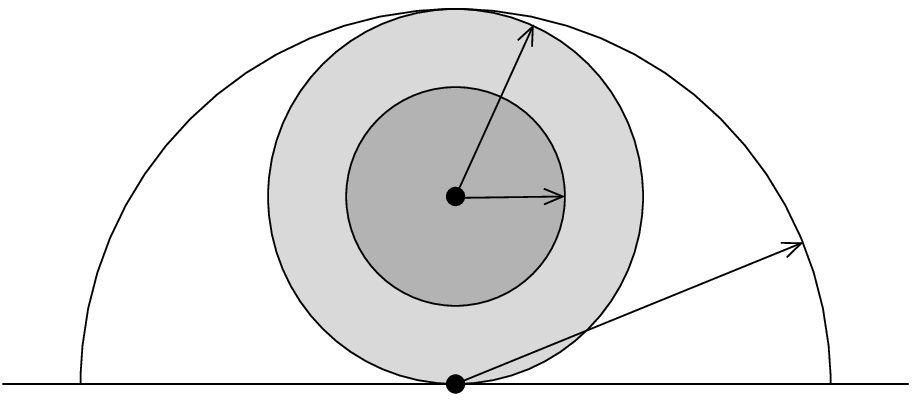}%
\end{picture}%
\setlength{\unitlength}{3947sp}%
\begingroup\makeatletter\ifx\SetFigFont\undefined
\def\x##1##2##3##4##5##6##7\relax{\def\x{##1##2##3##4##5##6}}%
\expandafter\x\fmtname xxxxxx\relax \def\y{splain}%
\ifx\x\y   % LaTeX or SliTeX?
\gdef\SetFigFont##1##2##3{%
  \ifnum ##1<17\tiny\else \ifnum ##1<20\small\else
  \ifnum ##1<24\normalsize\else \ifnum ##1<29\large\else
  \ifnum ##1<34\Large\else \ifnum ##1<41\LARGE\else
     \huge\fi\fi\fi\fi\fi\fi
  \csname ##3\endcsname}%
\else
\gdef\SetFigFont##1##2##3{\begingroup
  \count@##1\relax \ifnum 25<\count@\count@25\fi
  \def\x{\endgroup\@setsize\SetFigFont{##2pt}}%
  \expandafter\x
    \csname \romannumeral\the\count@ pt\expandafter\endcsname
    \csname @\romannumeral\the\count@ pt\endcsname
  \csname ##3\endcsname}%
\fi
\fi\endgroup
\begin{picture}(4374,2092)(214,-1544)
\put(2476,-586){\makebox(0,0)[lb]{\smash{\small{$2M_0$}%
}}}
\put(2476,239){\makebox(0,0)[lb]{\smash{\small{$R$}%
}}}
\put(3451,-1036){\makebox(0,0)[lb]{\smash{\small{$2R$}%
}}}
\put(2251,-1486){\makebox(0,0)[lb]{\smash{\small{$\gamma(t)$}%
}}}
\put(3301,-1486){\makebox(0,0)[lb]{\smash{\small{$\gamma$}%
}}}
\end{picture}
\end{center}}

\newcommand{\drawdefiance}{\begin{center}%
\begin{picture}(0,0)%
\includegraphics{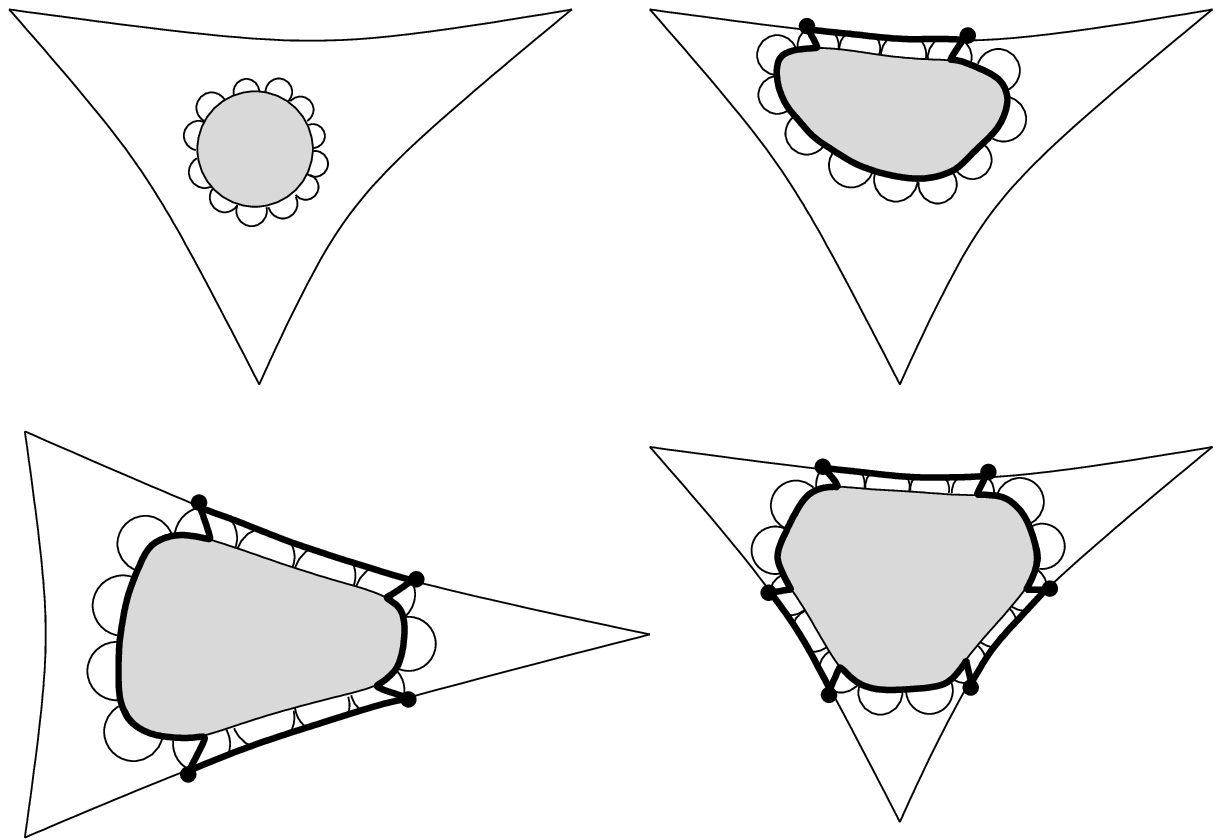}%
\end{picture}%
\setlength{\unitlength}{3947sp}%
\begingroup\makeatletter\ifx\SetFigFont\undefined
\def\x##1##2##3##4##5##6##7\relax{\def\x{##1##2##3##4##5##6}}%
\expandafter\x\fmtname xxxxxx\relax \def\y{splain}%
\ifx\x\y   % LaTeX or SliTeX?
\gdef\SetFigFont##1##2##3{%
  \ifnum ##1<17\tiny\else \ifnum ##1<20\small\else
  \ifnum ##1<24\normalsize\else \ifnum ##1<29\large\else
  \ifnum ##1<34\Large\else \ifnum ##1<41\LARGE\else
     \huge\fi\fi\fi\fi\fi\fi
  \csname ##3\endcsname}%
\else
\gdef\SetFigFont##1##2##3{\begingroup
  \count@##1\relax \ifnum 25<\count@\count@25\fi
  \def\x{\endgroup\@setsize\SetFigFont{##2pt}}%
  \expandafter\x
    \csname \romannumeral\the\count@ pt\expandafter\endcsname
    \csname @\romannumeral\the\count@ pt\endcsname
  \csname ##3\endcsname}%
\fi
\fi\endgroup
\begin{picture}(5937,4079)(151,-3748)
\put(301,-1186){\makebox(0,0)[lb]{\smash{\small{(a)}%
}}}
\put(3451,-1186){\makebox(0,0)[lb]{\smash{\small{(b)}%
}}}
\put(3451,-3361){\makebox(0,0)[lb]{\smash{\small{(d)}%
}}}
\put(151,-3361){\makebox(0,0)[lb]{\smash{\small{(c)}%
}}}
\put(4491,-491){\makebox(0,0)[lb]{\smash{\small{$\xi$}%
}}}
\put(4436,199){\makebox(0,0)[lb]{\smash{\small{$s'$}%
}}}
\put(1646,-3421){\makebox(0,0)[lb]{\smash{\small{$t'$}%
}}}
\put(1641,-2251){\makebox(0,0)[lb]{\smash{\small{$s'$}%
}}}
\end{picture}
\end{center} }%

\newcommand{\drawPointPlaneHull}{\begin{center}%
\begin{picture}(0,0)%
\includegraphics{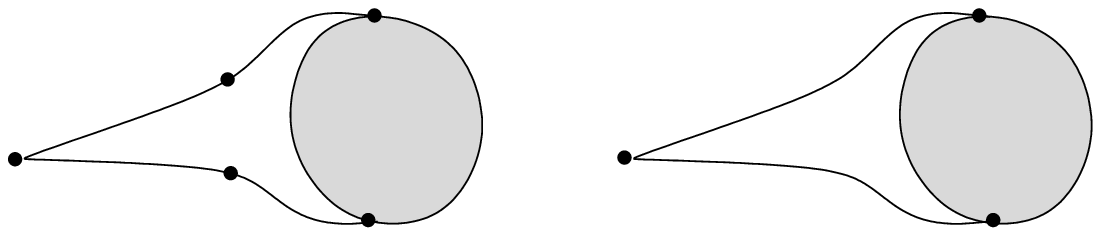}%
\end{picture}%
\setlength{\unitlength}{3947sp}%
\begingroup\makeatletter\ifx\SetFigFont\undefined
\def\x##1##2##3##4##5##6##7\relax{\def\x{##1##2##3##4##5##6}}%
\expandafter\x\fmtname xxxxxx\relax \def\y{splain}%
\ifx\x\y   % LaTeX or SliTeX?
\gdef\SetFigFont##1##2##3{%
  \ifnum ##1<17\tiny\else \ifnum ##1<20\small\else
  \ifnum ##1<24\normalsize\else \ifnum ##1<29\large\else
  \ifnum ##1<34\Large\else \ifnum ##1<41\LARGE\else
     \huge\fi\fi\fi\fi\fi\fi
  \csname ##3\endcsname}%
\else
\gdef\SetFigFont##1##2##3{\begingroup
  \count@##1\relax \ifnum 25<\count@\count@25\fi
  \def\x{\endgroup\@setsize\SetFigFont{##2pt}}%
  \expandafter\x
    \csname \romannumeral\the\count@ pt\expandafter\endcsname
    \csname @\romannumeral\the\count@ pt\endcsname
  \csname ##3\endcsname}%
\fi
\fi\endgroup
\begin{picture}(5392,1384)(84,-1349)
\put(601,-608){\makebox(0,0)[lb]{\smash{\small{$\iota_0$}%
}}}
\put(601,-1074){\makebox(0,0)[lb]{\smash{\small{$\iota_1$}%
}}}
\put(1470,-1291){\makebox(0,0)[lb]{\smash{\small{$\upsilon$}%
}}}
\put(1486,-121){\makebox(0,0)[lb]{\smash{\small{$\omega$}%
}}}
\put(2536,-706){\makebox(0,0)[lb]{\smash{\small{$\eta$}%
}}}
\put(1966,-706){\makebox(0,0)[lb]{\smash{\small{$P$}%
}}}
\put(3924,-1118){\makebox(0,0)[lb]{\smash{\small{$\alpha_1$}%
}}}
\put(3871,-473){\makebox(0,0)[lb]{\smash{\small{$\alpha_0$}%
}}}
\put(4583,-721){\makebox(0,0)[lb]{\smash{\small{$\xi$}%
}}}
\put( 84,-781){\makebox(0,0)[lb]{\smash{\small{$x$}%
}}}
\put(151,-136){\makebox(0,0)[lb]{\smash{\small{(a)}%
}}}
\put(3076,-136){\makebox(0,0)[lb]{\smash{\small{(b)}%
}}}
\put(5476,-736){\makebox(0,0)[lb]{\smash{\small{$\eta$}%
}}}
\end{picture}
\end{center} }

\newcommand{\drawdisjoint}{\begin{center}
\begin{picture}(0,0)%
\includegraphics{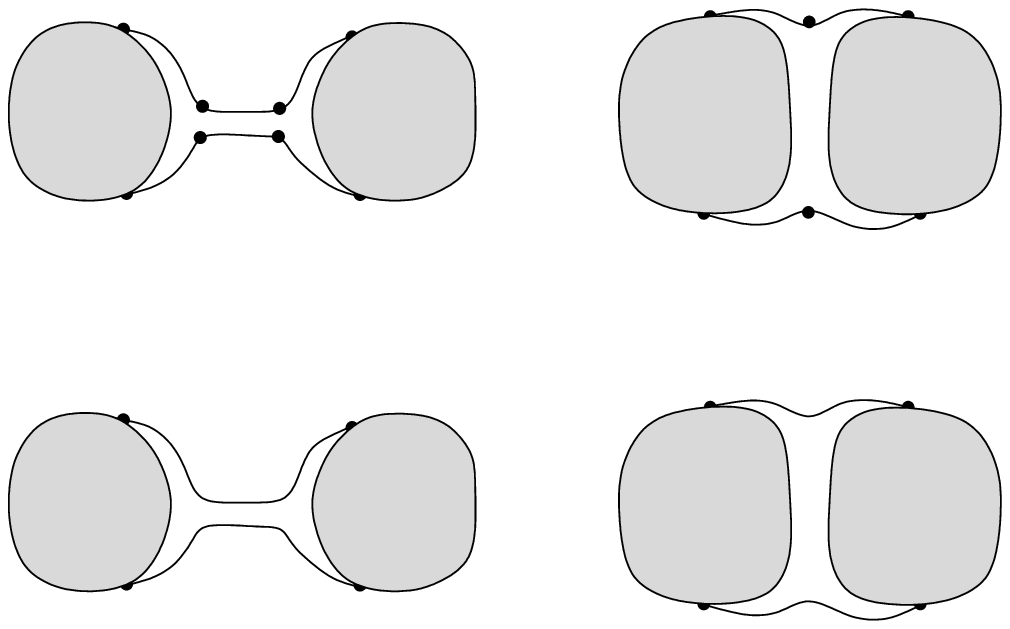}%
\end{picture}%
\setlength{\unitlength}{3947sp}%
\begingroup\makeatletter\ifx\SetFigFont\undefined
\def\x##1##2##3##4##5##6##7\relax{\def\x{##1##2##3##4##5##6}}%
\expandafter\x\fmtname xxxxxx\relax \def\y{splain}%
\ifx\x\y   % LaTeX or SliTeX?
\gdef\SetFigFont##1##2##3{%
  \ifnum ##1<17\tiny\else \ifnum ##1<20\small\else
  \ifnum ##1<24\normalsize\else \ifnum ##1<29\large\else
  \ifnum ##1<34\Large\else \ifnum ##1<41\LARGE\else
     \huge\fi\fi\fi\fi\fi\fi
  \csname ##3\endcsname}%
\else
\gdef\SetFigFont##1##2##3{\begingroup
  \count@##1\relax \ifnum 25<\count@\count@25\fi
  \def\x{\endgroup\@setsize\SetFigFont{##2pt}}%
  \expandafter\x
    \csname \romannumeral\the\count@ pt\expandafter\endcsname
    \csname @\romannumeral\the\count@ pt\endcsname
  \csname ##3\endcsname}%
\fi
\fi\endgroup
\begin{picture}(5130,3348)(144,-2744)
\put(691,-171){\makebox(0,0)[lb]{\smash{\small{$P_0$}%
}}}
\put(2255,-165){\makebox(0,0)[lb]{\smash{\small{$P_1$}%
}}}
\put(1280,107){\makebox(0,0)[lb]{\smash{\small{$\upsilon_0$}%
}}}
\put(1498,-70){\makebox(0,0)[lb]{\smash{\small{$\iota_0$}%
}}}
\put(1264,-569){\makebox(0,0)[lb]{\smash{\small{$\omega_0$}%
}}}
\put(1758,-593){\makebox(0,0)[lb]{\smash{\small{$\upsilon_1$}%
}}}
\put(1522,-407){\makebox(0,0)[lb]{\smash{\small{$\iota_1$}%
}}}
\put(1627,120){\makebox(0,0)[lb]{\smash{\small{$\omega_1$}%
}}}
\put(188,-164){\makebox(0,0)[lb]{\smash{\small{$\eta_0$}%
}}}
\put(2728,-177){\makebox(0,0)[lb]{\smash{\small{$\eta_1$}%
}}}
\put(3715,-187){\makebox(0,0)[lb]{\smash{\small{$P_0$}%
}}}
\put(4251,-543){\makebox(0,0)[lb]{\smash{\small{$\iota_1$}%
}}}
\put(4702,-207){\makebox(0,0)[lb]{\smash{\small{$P_1$}%
}}}
\put(5238,-207){\makebox(0,0)[lb]{\smash{\small{$\eta_1$}%
}}}
\put(4026,448){\makebox(0,0)[lb]{\smash{\small{$\upsilon_0$}%
}}}
\put(4501,444){\makebox(0,0)[lb]{\smash{\small{$\omega_1$}%
}}}
\put(3144,-184){\makebox(0,0)[lb]{\smash{\small{$\eta_0$}%
}}}
\put(4236,107){\makebox(0,0)[lb]{\smash{\small{$\iota_0$}%
}}}
\put(3984,-872){\makebox(0,0)[lb]{\smash{\small{$\omega_0$}%
}}}
\put(4503,-883){\makebox(0,0)[lb]{\smash{\small{$\upsilon_1$}%
}}}
\put(144,329){\makebox(0,0)[lb]{\smash{\small{(a)}%
}}}
\put(3151,314){\makebox(0,0)[lb]{\smash{\small{(b)}%
}}}
\put(196,-2056){\makebox(0,0)[lb]{\smash{\small{$\eta_0$}%
}}}
\put(2723,-2063){\makebox(0,0)[lb]{\smash{\small{$\eta_1$}%
}}}
\put(1501,-1929){\makebox(0,0)[lb]{\smash{\small{$\alpha_0$}%
}}}
\put(1028,-2071){\makebox(0,0)[lb]{\smash{\small{$\xi_0$}%
}}}
\put(1936,-2078){\makebox(0,0)[lb]{\smash{\small{$\xi_1$}%
}}}
\put(3151,-2064){\makebox(0,0)[lb]{\smash{\small{$\eta_0$}%
}}}
\put(4006,-2064){\makebox(0,0)[lb]{\smash{\small{$\xi_0$}%
}}}
\put(4418,-2071){\makebox(0,0)[lb]{\smash{\small{$\xi_1$}%
}}}
\put(4238,-1486){\makebox(0,0)[lb]{\smash{\small{$\alpha_0$}%
}}}
\put(5274,-2071){\makebox(0,0)[lb]{\smash{\small{$\eta_1$}%
}}}
\put(144,-1546){\makebox(0,0)[lb]{\smash{\small{(c)}%
}}}
\put(3151,-1561){\makebox(0,0)[lb]{\smash{\small{(d)}%
}}}
\put(1509,-2304){\makebox(0,0)[lb]{\smash{\small{$\alpha_1$}%
}}}
\put(4246,-2686){\makebox(0,0)[lb]{\smash{\small{$\alpha_1$}%
}}}
\end{picture}
\end{center}}

\newcommand{\drawBetweenPreflats}{\begin{center}
\begin{picture}(0,0)%
\includegraphics{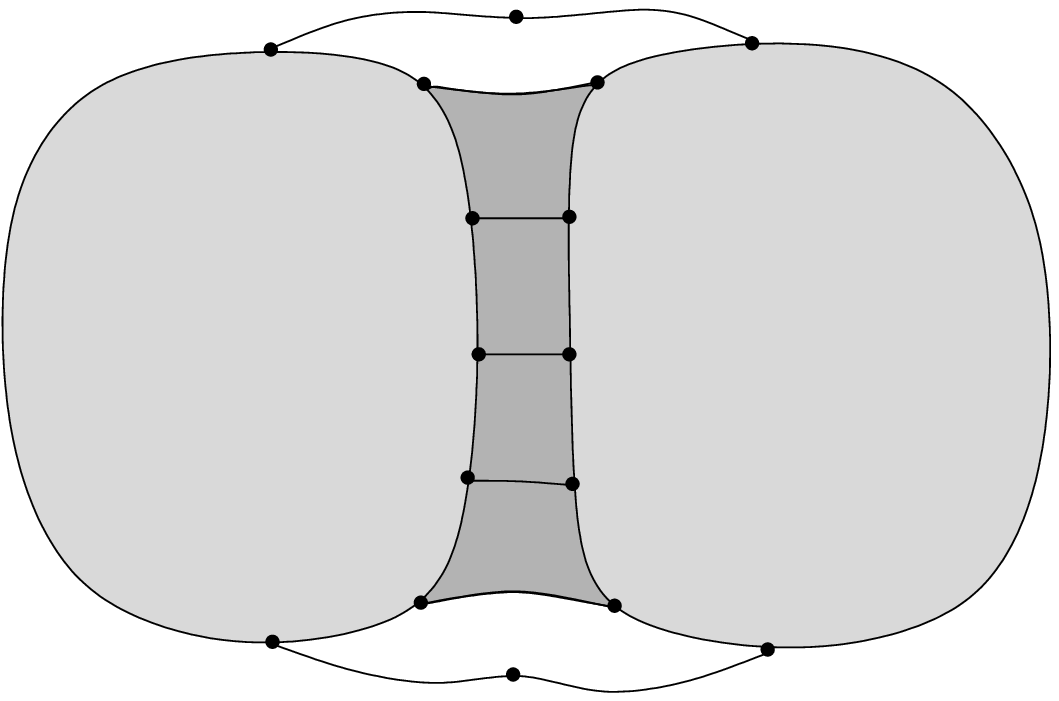}%
\end{picture}%
\setlength{\unitlength}{3947sp}%
\begingroup\makeatletter\ifx\SetFigFont\undefined
\def\x##1##2##3##4##5##6##7\relax{\def\x{##1##2##3##4##5##6}}%
\expandafter\x\fmtname xxxxxx\relax \def\y{splain}%
\ifx\x\y   % LaTeX or SliTeX?
\gdef\SetFigFont##1##2##3{%
  \ifnum ##1<17\tiny\else \ifnum ##1<20\small\else
  \ifnum ##1<24\normalsize\else \ifnum ##1<29\large\else
  \ifnum ##1<34\Large\else \ifnum ##1<41\LARGE\else
     \huge\fi\fi\fi\fi\fi\fi
  \csname ##3\endcsname}%
\else
\gdef\SetFigFont##1##2##3{\begingroup
  \count@##1\relax \ifnum 25<\count@\count@25\fi
  \def\x{\endgroup\@setsize\SetFigFont{##2pt}}%
  \expandafter\x
    \csname \romannumeral\the\count@ pt\expandafter\endcsname
    \csname @\romannumeral\the\count@ pt\endcsname
  \csname ##3\endcsname}%
\fi
\fi\endgroup
\begin{picture}(5053,3704)(426,-3176)
\put(2824,-443){\makebox(0,0)[lb]{\smash{\small{$Q_1$}%
}}}
\put(2859,-1083){\makebox(0,0)[lb]{\smash{\small{$Q_2$}%
}}}
\put(2819,-2288){\makebox(0,0)[lb]{\smash{\small{$Q_\ell$}%
}}}
\put(2829,-3118){\makebox(0,0)[lb]{\smash{\small{$x_1$}%
}}}
\put(2839,372){\makebox(0,0)[lb]{\smash{\small{$x_0$}%
}}}
\put(3239,-2028){\makebox(0,0)[lb]{\smash{\small{$b_{\ell-1}$}%
}}}
\put(3409,-2513){\makebox(0,0)[lb]{\smash{\small{$b_\ell = q_1$}%
}}}
\put(3244,-1413){\makebox(0,0)[lb]{\smash{\small{$b_2$}%
}}}
\put(3249,-763){\makebox(0,0)[lb]{\smash{\small{$b_1$}%
}}}
\put(3364,-173){\makebox(0,0)[lb]{\smash{\small{$b_0 = p_1$}%
}}}
\put(1784,-158){\makebox(0,0)[lb]{\smash{\small{$a_0 = p_0$}%
}}}
\put(2444,-753){\makebox(0,0)[lb]{\smash{\small{$a_1$}%
}}}
\put(2474,-1408){\makebox(0,0)[lb]{\smash{\small{$a_2$}%
}}}
\put(2264,-1993){\makebox(0,0)[lb]{\smash{\small{$a_{\ell-1}$}%
}}}
\put(2901,-1728){\makebox(0,0)[lb]{\smash{\small{$\vdots$}%
}}}
\put(1742,-2490){\makebox(0,0)[lb]{\smash{\small{$a_\ell = q_0$}%
}}}
\put(1463,-1346){\makebox(0,0)[lb]{\smash{\small{$P_0$}%
}}}
\put(4141,-1361){\makebox(0,0)[lb]{\smash{\small{$P_1$}%
}}}
\end{picture}
\end{center} }

\newcommand{\drawintersectingone}{\begin{center}
\begin{picture}(0,0)%
\includegraphics{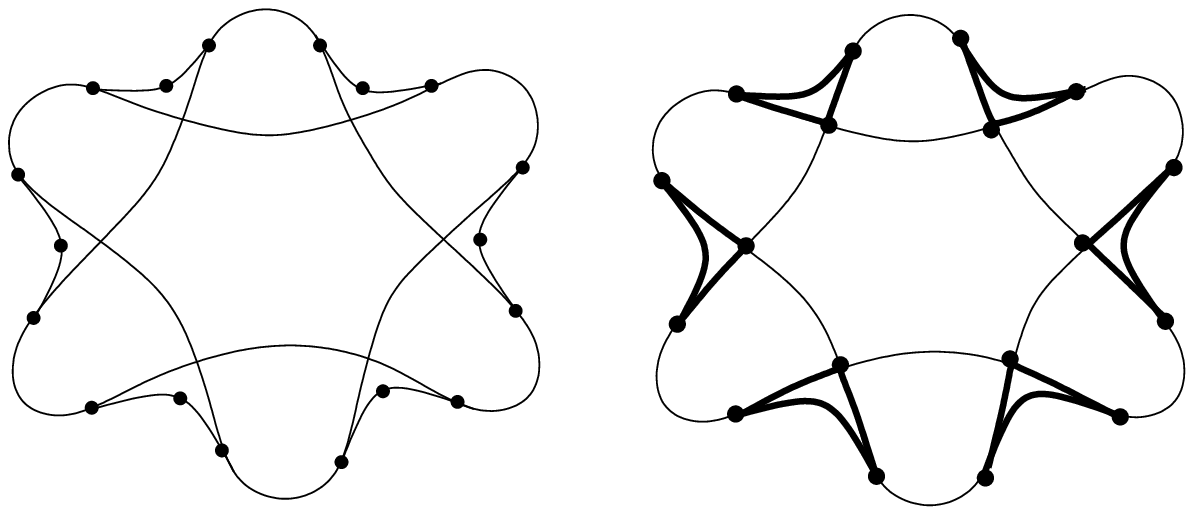}%
\end{picture}%
\setlength{\unitlength}{3947sp}%
\begingroup\makeatletter\ifx\SetFigFont\undefined
\def\x##1##2##3##4##5##6##7\relax{\def\x{##1##2##3##4##5##6}}%
\expandafter\x\fmtname xxxxxx\relax \def\y{splain}%
\ifx\x\y   % LaTeX or SliTeX?
\gdef\SetFigFont##1##2##3{%
  \ifnum ##1<17\tiny\else \ifnum ##1<20\small\else
  \ifnum ##1<24\normalsize\else \ifnum ##1<29\large\else
  \ifnum ##1<34\Large\else \ifnum ##1<41\LARGE\else
     \huge\fi\fi\fi\fi\fi\fi
  \csname ##3\endcsname}%
\else
\gdef\SetFigFont##1##2##3{\begingroup
  \count@##1\relax \ifnum 25<\count@\count@25\fi
  \def\x{\endgroup\@setsize\SetFigFont{##2pt}}%
  \expandafter\x
    \csname \romannumeral\the\count@ pt\expandafter\endcsname
    \csname @\romannumeral\the\count@ pt\endcsname
  \csname ##3\endcsname}%
\fi
\fi\endgroup
\begin{picture}(5774,2777)(266,-2467)
\put(303,140){\makebox(0,0)[lb]{\smash{\small{(a)}%
}}}
\put(2738,-1209){\makebox(0,0)[lb]{\smash{\small{$\omega_3$}%
}}}
\put(2760,-919){\makebox(0,0)[lb]{\smash{\small{$\upsilon_2$}%
}}}
\put(1545,-237){\makebox(0,0)[lb]{\smash{\small{$P_1$}%
}}}
\put(1658,-2409){\makebox(0,0)[lb]{\smash{\small{$\eta_4$}%
}}}
\put(2091,-1994){\makebox(0,0)[lb]{\smash{\small{$\omega_4$}%
}}}
\put(2318,-1886){\makebox(0,0)[lb]{\smash{\small{$\upsilon_3$}%
}}}
\put(2914,-1857){\makebox(0,0)[lb]{\smash{\small{$\eta_3$}%
}}}
\put(2913,-305){\makebox(0,0)[lb]{\smash{\small{$\eta_2$}%
}}}
\put(1971,-181){\makebox(0,0)[lb]{\smash{\small{$\upsilon_1$}%
}}}
\put(2186,-249){\makebox(0,0)[lb]{\smash{\small{$\omega_2$}%
}}}
\put(2430,-595){\makebox(0,0)[lb]{\smash{\small{$P_0$}%
}}}
\put(1402,-1062){\makebox(0,0)[lb]{\smash{\small{$P_0 \cap P_1$}%
}}}
\put(2477,-1527){\makebox(0,0)[lb]{\smash{\small{$P_1$}%
}}}
\put(1596,-1925){\makebox(0,0)[lb]{\smash{\small{$P_0$}%
}}}
\put(1181,-2005){\makebox(0,0)[lb]{\smash{\small{$\upsilon_4$}%
}}}
\put(947,-1926){\makebox(0,0)[lb]{\smash{\small{$\omega_5$}%
}}}
\put(732,-1510){\makebox(0,0)[lb]{\smash{\small{$P_1$}%
}}}
\put(407,-1226){\makebox(0,0)[lb]{\smash{\small{$\upsilon_5$}%
}}}
\put(367,-959){\makebox(0,0)[lb]{\smash{\small{$\omega_0$}%
}}}
\put(266,-368){\makebox(0,0)[lb]{\smash{\small{$\eta_0$}%
}}}
\put(1533,154){\makebox(0,0)[lb]{\smash{\small{$\eta_1$}%
}}}
\put(856,-261){\makebox(0,0)[lb]{\smash{\small{$\upsilon_0$}%
}}}
\put(1096,-157){\makebox(0,0)[lb]{\smash{\small{$\omega_1$}%
}}}
\put(652,-651){\makebox(0,0)[lb]{\smash{\small{$P_0$}%
}}}
\put(339,-1902){\makebox(0,0)[lb]{\smash{\small{$\eta_5$}%
}}}
\put(3499, 23){\makebox(0,0)[lb]{\smash{\small{(b)}%
}}}
\end{picture}
\end{center}}

\newcommand{\drawintersectingtwo}{\begin{center}
\begin{picture}(0,0)%
\includegraphics{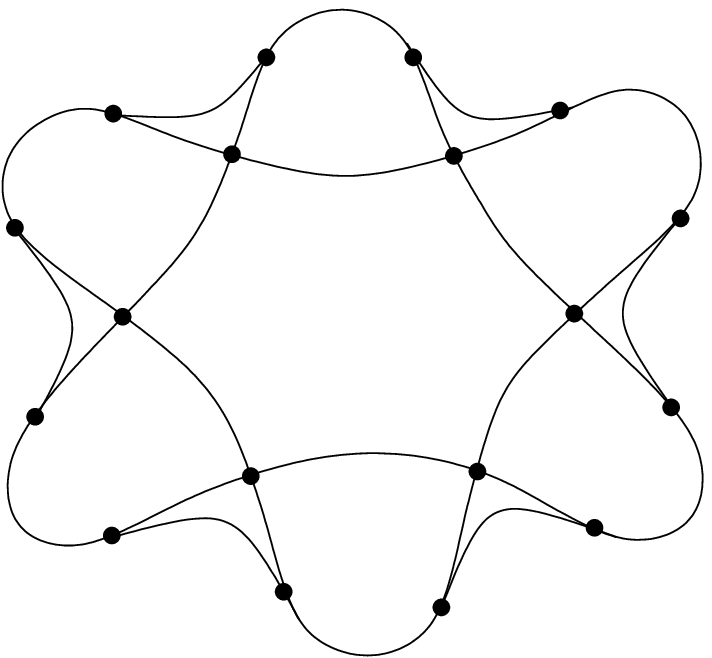}%
\end{picture}%
\setlength{\unitlength}{3947sp}%
\begingroup\makeatletter\ifx\SetFigFont\undefined
\def\x##1##2##3##4##5##6##7\relax{\def\x{##1##2##3##4##5##6}}%
\expandafter\x\fmtname xxxxxx\relax \def\y{splain}%
\ifx\x\y   % LaTeX or SliTeX?
\gdef\SetFigFont##1##2##3{%
  \ifnum ##1<17\tiny\else \ifnum ##1<20\small\else
  \ifnum ##1<24\normalsize\else \ifnum ##1<29\large\else
  \ifnum ##1<34\Large\else \ifnum ##1<41\LARGE\else
     \huge\fi\fi\fi\fi\fi\fi
  \csname ##3\endcsname}%
\else
\gdef\SetFigFont##1##2##3{\begingroup
  \count@##1\relax \ifnum 25<\count@\count@25\fi
  \def\x{\endgroup\@setsize\SetFigFont{##2pt}}%
  \expandafter\x
    \csname \romannumeral\the\count@ pt\expandafter\endcsname
    \csname @\romannumeral\the\count@ pt\endcsname
  \csname ##3\endcsname}%
\fi
\fi\endgroup
\begin{picture}(3384,3124)(1100,-2501)
\put(1786,-173){\makebox(0,0)[lb]{\smash{\small{$\upsilon'_0$}%
}}}
\put(1936,194){\makebox(0,0)[lb]{\smash{\small{$\alpha_0$}%
}}}
\put(2333, 44){\makebox(0,0)[lb]{\smash{\small{$\omega'_1$}%
}}}
\put(2948, 37){\makebox(0,0)[lb]{\smash{\small{$\upsilon'_1$}%
}}}
\put(3353,164){\makebox(0,0)[lb]{\smash{\small{$\alpha_1$}%
}}}
\put(3571,-159){\makebox(0,0)[lb]{\smash{\small{$\omega'_2$}%
}}}
\put(4163,-901){\makebox(0,0)[lb]{\smash{\small{$\alpha_2$}%
}}}
\put(3624,-1659){\makebox(0,0)[lb]{\smash{\small{$\upsilon'_3$}%
}}}
\put(3481,-2011){\makebox(0,0)[lb]{\smash{\small{$\alpha_3$}%
}}}
\put(3045,-1891){\makebox(0,0)[lb]{\smash{\small{$\omega'_4$}%
}}}
\put(2438,-1898){\makebox(0,0)[lb]{\smash{\small{$\upsilon'_4$}%
}}}
\put(1823,-1643){\makebox(0,0)[lb]{\smash{\small{$\omega'_5$}%
}}}
\put(1523,-1217){\makebox(0,0)[lb]{\smash{\small{$\upsilon'_5$}%
}}}
\put(1163,-938){\makebox(0,0)[lb]{\smash{\small{$\alpha_5$}%
}}}
\put(1471,-631){\makebox(0,0)[lb]{\smash{\small{$\omega'_0$}%
}}}
\put(3878,-1239){\makebox(0,0)[lb]{\smash{\small{$\omega'_3$}%
}}}
\put(2056,-2048){\makebox(0,0)[lb]{\smash{\small{$\alpha_4$}%
}}}
\put(3893,-555){\makebox(0,0)[lb]{\smash{\small{$\upsilon'_2$}%
}}}
\end{picture}
\end{center}}

\newcommand{\drawFlatClosureOne}{\begin{center}%
\begin{picture}(0,0)%
\includegraphics{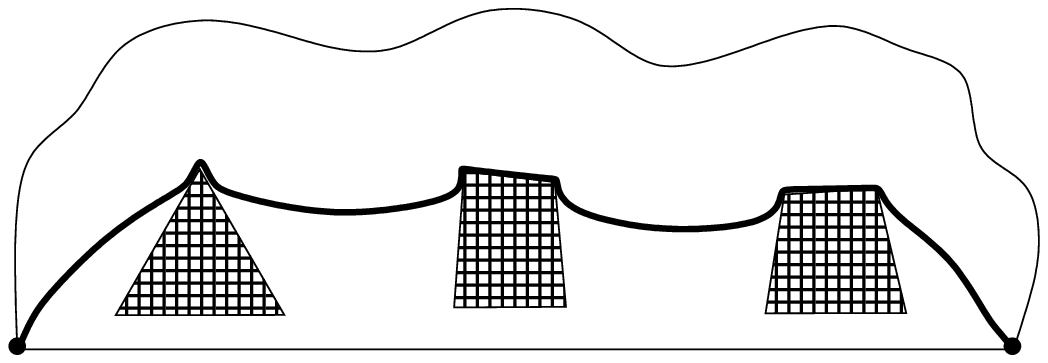}%
\end{picture}%
\setlength{\unitlength}{3947sp}%
\begingroup\makeatletter\ifx\SetFigFont\undefined
\def\x##1##2##3##4##5##6##7\relax{\def\x{##1##2##3##4##5##6}}%
\expandafter\x\fmtname xxxxxx\relax \def\y{splain}%
\ifx\x\y   % LaTeX or SliTeX?
\gdef\SetFigFont##1##2##3{%
  \ifnum ##1<17\tiny\else \ifnum ##1<20\small\else
  \ifnum ##1<24\normalsize\else \ifnum ##1<29\large\else
  \ifnum ##1<34\Large\else \ifnum ##1<41\LARGE\else
     \huge\fi\fi\fi\fi\fi\fi
  \csname ##3\endcsname}%
\else
\gdef\SetFigFont##1##2##3{\begingroup
  \count@##1\relax \ifnum 25<\count@\count@25\fi
  \def\x{\endgroup\@setsize\SetFigFont{##2pt}}%
  \expandafter\x
    \csname \romannumeral\the\count@ pt\expandafter\endcsname
    \csname @\romannumeral\the\count@ pt\endcsname
  \csname ##3\endcsname}%
\fi
\fi\endgroup
\begin{picture}(4962,1906)(555,-1896)
\put(3594,-1838){\makebox(0,0)[lb]{\smash{\small{$\gamma$}%
}}}
\put(3705,-188){\makebox(0,0)[lb]{\smash{\small{$\alpha$}%
}}}
\put(3640,-953){\makebox(0,0)[lb]{\smash{\small{$\beta$}%
}}}
\put(2036,-623){\makebox(0,0)[lb]{\smash{\small{$D_\alpha$}%
}}}
\put(2125,-1328){\makebox(0,0)[lb]{\smash{\small{$D_\gamma$}%
}}}
\end{picture}
\end{center} }

\newcommand{\drawflatclosure}{\begin{center}
\begin{picture}(0,0)%
\includegraphics{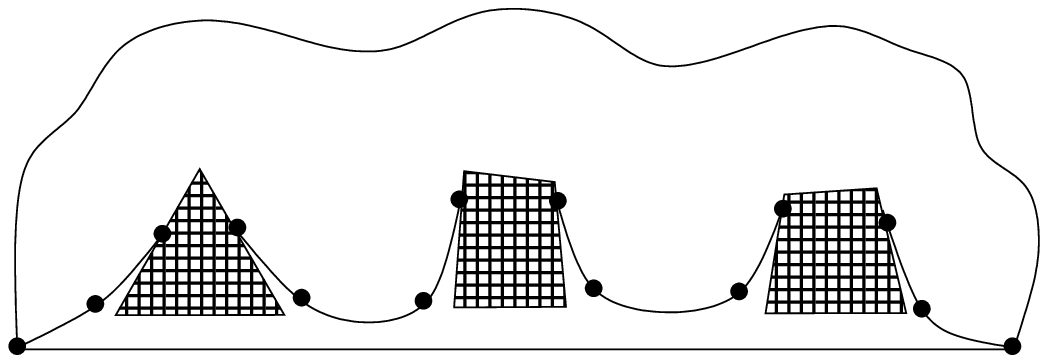}%
\end{picture}%
\setlength{\unitlength}{3947sp}%
\begingroup\makeatletter\ifx\SetFigFont\undefined
\def\x##1##2##3##4##5##6##7\relax{\def\x{##1##2##3##4##5##6}}%
\expandafter\x\fmtname xxxxxx\relax \def\y{splain}%
\ifx\x\y   % LaTeX or SliTeX?
\gdef\SetFigFont##1##2##3{%
  \ifnum ##1<17\tiny\else \ifnum ##1<20\small\else
  \ifnum ##1<24\normalsize\else \ifnum ##1<29\large\else
  \ifnum ##1<34\Large\else \ifnum ##1<41\LARGE\else
     \huge\fi\fi\fi\fi\fi\fi
  \csname ##3\endcsname}%
\else
\gdef\SetFigFont##1##2##3{\begingroup
  \count@##1\relax \ifnum 25<\count@\count@25\fi
  \def\x{\endgroup\@setsize\SetFigFont{##2pt}}%
  \expandafter\x
    \csname \romannumeral\the\count@ pt\expandafter\endcsname
    \csname @\romannumeral\the\count@ pt\endcsname
  \csname ##3\endcsname}%
\fi
\fi\endgroup
\begin{picture}(4962,1906)(555,-1896)
\put(631,-1433){\makebox(0,0)[lb]{\smash{\small{$\tau_0$}%
}}}
\put(923,-1208){\makebox(0,0)[lb]{\smash{\small{$\pi_1$}%
}}}
\put(1545,-818){\makebox(0,0)[lb]{\smash{\small{$\rho_1$}%
}}}
\put(2146,-1425){\makebox(0,0)[lb]{\smash{\small{$\tau_1$}%
}}}
\put(1801,-1216){\makebox(0,0)[lb]{\smash{\small{$\sigma_1$}%
}}}
\put(2934,-713){\makebox(0,0)[lb]{\smash{\small{$\rho_2$}%
}}}
\put(3661,-1359){\makebox(0,0)[lb]{\smash{\small{$\tau_2$}%
}}}
\put(4433,-774){\makebox(0,0)[lb]{\smash{\small{$\rho_3$}%
}}}
\put(4891,-1231){\makebox(0,0)[lb]{\smash{\small{$\sigma_3$}%
}}}
\put(5109,-1493){\makebox(0,0)[lb]{\smash{\small{$\tau_3$}%
}}}
\put(3594,-1838){\makebox(0,0)[lb]{\smash{\small{$\gamma$}%
}}}
\put(3705,-188){\makebox(0,0)[lb]{\smash{\small{$\alpha$}%
}}}
\put(2393,-1193){\makebox(0,0)[lb]{\smash{\small{$\pi_2$}%
}}}
\put(3314,-1165){\makebox(0,0)[lb]{\smash{\small{$\sigma_2$}%
}}}
\put(3969,-1163){\makebox(0,0)[lb]{\smash{\small{$\pi_3$}%
}}}
\end{picture}
\end{center}}

\newcommand{\drawFCstagezero}{\begin{center}
\begin{picture}(0,0)%
\includegraphics{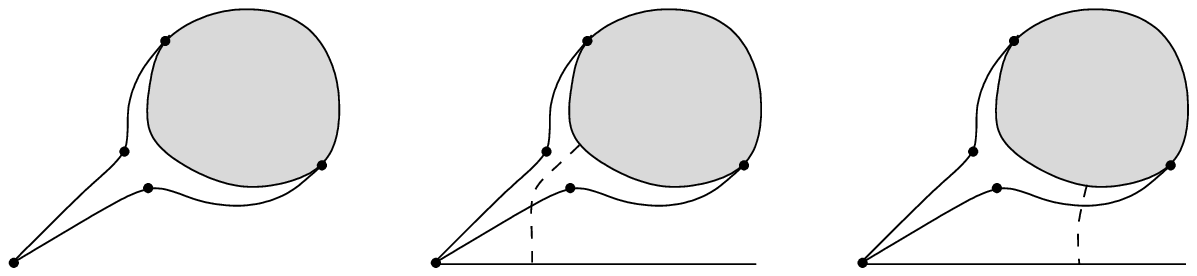}%
\end{picture}%
\setlength{\unitlength}{3947sp}%
\begingroup\makeatletter\ifx\SetFigFont\undefined
\def\x##1##2##3##4##5##6##7\relax{\def\x{##1##2##3##4##5##6}}%
\expandafter\x\fmtname xxxxxx\relax \def\y{splain}%
\ifx\x\y   % LaTeX or SliTeX?
\gdef\SetFigFont##1##2##3{%
  \ifnum ##1<17\tiny\else \ifnum ##1<20\small\else
  \ifnum ##1<24\normalsize\else \ifnum ##1<29\large\else
  \ifnum ##1<34\Large\else \ifnum ##1<41\LARGE\else
     \huge\fi\fi\fi\fi\fi\fi
  \csname ##3\endcsname}%
\else
\gdef\SetFigFont##1##2##3{\begingroup
  \count@##1\relax \ifnum 25<\count@\count@25\fi
  \def\x{\endgroup\@setsize\SetFigFont{##2pt}}%
  \expandafter\x
    \csname \romannumeral\the\count@ pt\expandafter\endcsname
    \csname @\romannumeral\the\count@ pt\endcsname
  \csname ##3\endcsname}%
\fi
\fi\endgroup
\begin{picture}(5859,1594)(127,-1555)
\put(427,-906){\makebox(0,0)[lb]{\smash{\small{$\iota_0$}%
}}}
\put(676,-1227){\makebox(0,0)[lb]{\smash{\small{$\iota_1$}%
}}}
\put(698,-467){\makebox(0,0)[lb]{\smash{\small{$\omega$}%
}}}
\put(1334,-1177){\makebox(0,0)[lb]{\smash{\small{$\upsilon$}%
}}}
\put(1737,-93){\makebox(0,0)[lb]{\smash{\small{$\eta$}%
}}}
\put(3405,-1475){\makebox(0,0)[lb]{\smash{\small{$\gamma$}%
}}}
\put(2862,-1190){\makebox(0,0)[lb]{\smash{\small{$\xi$}%
}}}
\put(5525,-1190){\makebox(0,0)[lb]{\smash{\small{$\xi$}%
}}}
\put(5115,-1468){\makebox(0,0)[lb]{\smash{\small{$\gamma$}%
}}}
\put(2665,-722){\makebox(0,0)[lb]{\smash{\small{$x_0$}%
}}}
\put(4728,-708){\makebox(0,0)[lb]{\smash{\small{$x_0$}%
}}}
\put(5013,-1110){\makebox(0,0)[lb]{\smash{\small{$x_1$}%
}}}
\put(4392,-204){\makebox(0,0)[lb]{\smash{\small{(c)}%
}}}
\put(2344,-204){\makebox(0,0)[lb]{\smash{\small{(b)}%
}}}
\put(222,-204){\makebox(0,0)[lb]{\smash{\small{(a)}%
}}}
\put(1349,-540){\makebox(0,0)[lb]{\smash{\small{$P_1$}%
}}}
\put(127,-1497){\makebox(0,0)[lb]{\smash{\small{$\beta(0)$}%
}}}
\end{picture}
\end{center}}

\newcommand{\drawFCstagezerob}{\begin{center}
\begin{picture}(0,0)%
\includegraphics{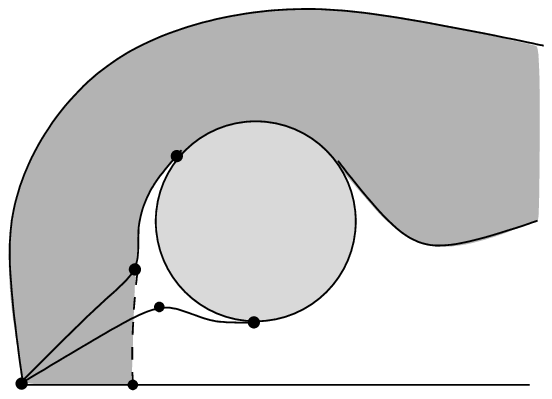}%
\end{picture}%
\setlength{\unitlength}{3947sp}%
\begingroup\makeatletter\ifx\SetFigFont\undefined
\def\x##1##2##3##4##5##6##7\relax{\def\x{##1##2##3##4##5##6}}%
\expandafter\x\fmtname xxxxxx\relax \def\y{splain}%
\ifx\x\y   % LaTeX or SliTeX?
\gdef\SetFigFont##1##2##3{%
  \ifnum ##1<17\tiny\else \ifnum ##1<20\small\else
  \ifnum ##1<24\normalsize\else \ifnum ##1<29\large\else
  \ifnum ##1<34\Large\else \ifnum ##1<41\LARGE\else
     \huge\fi\fi\fi\fi\fi\fi
  \csname ##3\endcsname}%
\else
\gdef\SetFigFont##1##2##3{\begingroup
  \count@##1\relax \ifnum 25<\count@\count@25\fi
  \def\x{\endgroup\@setsize\SetFigFont{##2pt}}%
  \expandafter\x
    \csname \romannumeral\the\count@ pt\expandafter\endcsname
    \csname @\romannumeral\the\count@ pt\endcsname
  \csname ##3\endcsname}%
\fi
\fi\endgroup
\begin{picture}(2885,2115)(855,-1845)
\put(1718,-1778){\makebox(0,0)[lb]{\smash{\small{$x_2$}%
}}}
\put(1463,-1778){\makebox(0,0)[lb]{\smash{\small{$\hat{\gamma}$}%
}}}
\put(3526,112){\makebox(0,0)[lb]{\smash{\small{$\alpha$}%
}}}
\put(3526,-713){\makebox(0,0)[lb]{\smash{\small{$\beta$}%
}}}
\put(3503,-1471){\makebox(0,0)[lb]{\smash{\small{$\gamma$}%
}}}
\put(855,-1787){\makebox(0,0)[lb]{\smash{\small{$\beta(0)$}%
}}}
\put(1529,-944){\makebox(0,0)[lb]{\smash{\small{$x_0$}%
}}}
\put(2566,-76){\makebox(0,0)[lb]{\smash{\small{$\hat D$}%
}}}
\end{picture}
\end{center}}

\newcommand{\drawFCstageia}{\begin{center}
\begin{picture}(0,0)%
\includegraphics{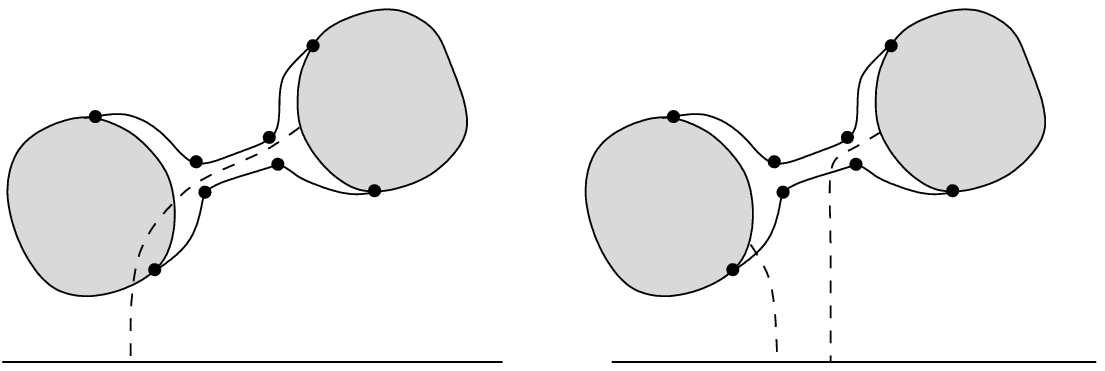}%
\end{picture}%
\setlength{\unitlength}{3947sp}%
\begingroup\makeatletter\ifx\SetFigFont\undefined
\def\x##1##2##3##4##5##6##7\relax{\def\x{##1##2##3##4##5##6}}%
\expandafter\x\fmtname xxxxxx\relax \def\y{splain}%
\ifx\x\y   % LaTeX or SliTeX?
\gdef\SetFigFont##1##2##3{%
  \ifnum ##1<17\tiny\else \ifnum ##1<20\small\else
  \ifnum ##1<24\normalsize\else \ifnum ##1<29\large\else
  \ifnum ##1<34\Large\else \ifnum ##1<41\LARGE\else
     \huge\fi\fi\fi\fi\fi\fi
  \csname ##3\endcsname}%
\else
\gdef\SetFigFont##1##2##3{\begingroup
  \count@##1\relax \ifnum 25<\count@\count@25\fi
  \def\x{\endgroup\@setsize\SetFigFont{##2pt}}%
  \expandafter\x
    \csname \romannumeral\the\count@ pt\expandafter\endcsname
    \csname @\romannumeral\the\count@ pt\endcsname
  \csname ##3\endcsname}%
\fi
\fi\endgroup
\begin{picture}(5274,1988)(364,-1619)
\put(3864,-1208){\makebox(0,0)[lb]{\smash{\small{$\xi_0$}%
}}}
\put(4403,-1029){\makebox(0,0)[lb]{\smash{\small{$\xi_1$}%
}}}
\put(4006,-279){\makebox(0,0)[lb]{\smash{\small{$x_0$}%
}}}
\put(4261,-196){\makebox(0,0)[lb]{\smash{\small{$x_1$}%
}}}
\put(3226,164){\makebox(0,0)[lb]{\smash{\small{(b)}%
}}}
\put(2026,-1561){\makebox(0,0)[lb]{\smash{\small{$\gamma$}%
}}}
\put(4876,-1561){\makebox(0,0)[lb]{\smash{\small{$\gamma$}%
}}}
\put(1404,-264){\makebox(0,0)[lb]{\smash{\small{$\iota_0$}%
}}}
\put(376, 89){\makebox(0,0)[lb]{\smash{\small{(a)}%
}}}
\put(1036,-1186){\makebox(0,0)[lb]{\smash{\small{$\xi_1$}%
}}}
\put(2041,-114){\makebox(0,0)[lb]{\smash{\small{$P_{i+1}$}%
}}}
\put(706,-616){\makebox(0,0)[lb]{\smash{\small{$P_i$}%
}}}
\put(3504,-638){\makebox(0,0)[lb]{\smash{\small{$P_i$}%
}}}
\put(4809,-151){\makebox(0,0)[lb]{\smash{\small{$P_{i+1}$}%
}}}
\end{picture}
\end{center}}

\newcommand{\drawFCstageib}{\begin{center}
\begin{picture}(0,0)%
\includegraphics{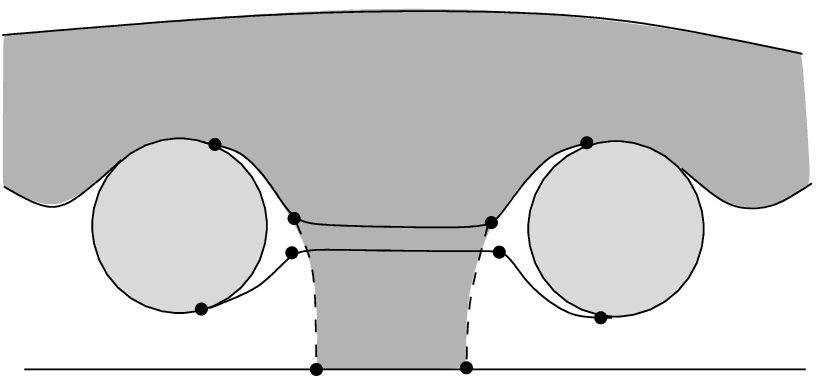}%
\end{picture}%
\setlength{\unitlength}{3947sp}%
\begingroup\makeatletter\ifx\SetFigFont\undefined
\def\x##1##2##3##4##5##6##7\relax{\def\x{##1##2##3##4##5##6}}%
\expandafter\x\fmtname xxxxxx\relax \def\y{splain}%
\ifx\x\y   % LaTeX or SliTeX?
\gdef\SetFigFont##1##2##3{%
  \ifnum ##1<17\tiny\else \ifnum ##1<20\small\else
  \ifnum ##1<24\normalsize\else \ifnum ##1<29\large\else
  \ifnum ##1<34\Large\else \ifnum ##1<41\LARGE\else
     \huge\fi\fi\fi\fi\fi\fi
  \csname ##3\endcsname}%
\else
\gdef\SetFigFont##1##2##3{\begingroup
  \count@##1\relax \ifnum 25<\count@\count@25\fi
  \def\x{\endgroup\@setsize\SetFigFont{##2pt}}%
  \expandafter\x
    \csname \romannumeral\the\count@ pt\expandafter\endcsname
    \csname @\romannumeral\the\count@ pt\endcsname
  \csname ##3\endcsname}%
\fi
\fi\endgroup
\begin{picture}(3909,2104)(1384,-1941)
\put(2791,-856){\makebox(0,0)[lb]{\smash{\small{$x_0$}%
}}}
\put(3601,-849){\makebox(0,0)[lb]{\smash{\small{$x_1$}%
}}}
\put(3578,-1853){\makebox(0,0)[lb]{\smash{\small{$y_1$}%
}}}
\put(2858,-1876){\makebox(0,0)[lb]{\smash{\small{$y_0$}%
}}}
\put(3211,-1883){\makebox(0,0)[lb]{\smash{\small{$\hat{\gamma}$}%
}}}
\put(5086,  7){\makebox(0,0)[lb]{\smash{\small{$\alpha$}%
}}}
\put(5131,-683){\makebox(0,0)[lb]{\smash{\small{$\beta$}%
}}}
\put(5094,-1546){\makebox(0,0)[lb]{\smash{\small{$\gamma$}%
}}}
\put(3271,-534){\makebox(0,0)[lb]{\smash{\small{$\hat{D}$}%
}}}
\end{picture}
\end{center}}

\newcommand{\drawFCstageic}{\begin{center}
\begin{picture}(0,0)%
\includegraphics{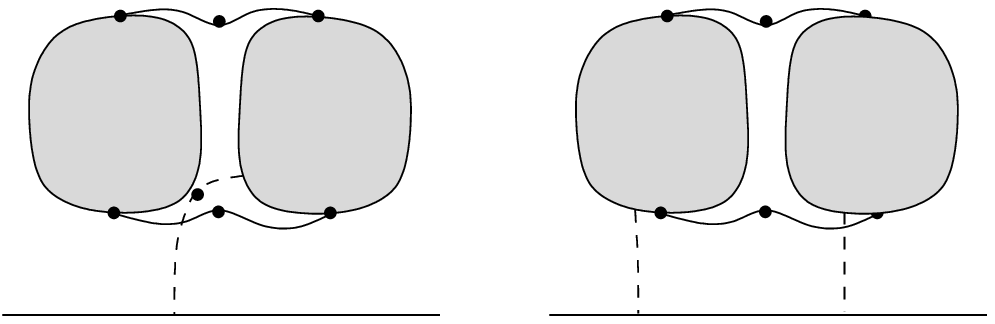}%
\end{picture}%
\setlength{\unitlength}{3947sp}%
\begingroup\makeatletter\ifx\SetFigFont\undefined
\def\x##1##2##3##4##5##6##7\relax{\def\x{##1##2##3##4##5##6}}%
\expandafter\x\fmtname xxxxxx\relax \def\y{splain}%
\ifx\x\y   % LaTeX or SliTeX?
\gdef\SetFigFont##1##2##3{%
  \ifnum ##1<17\tiny\else \ifnum ##1<20\small\else
  \ifnum ##1<24\normalsize\else \ifnum ##1<29\large\else
  \ifnum ##1<34\Large\else \ifnum ##1<41\LARGE\else
     \huge\fi\fi\fi\fi\fi\fi
  \csname ##3\endcsname}%
\else
\gdef\SetFigFont##1##2##3{\begingroup
  \count@##1\relax \ifnum 25<\count@\count@25\fi
  \def\x{\endgroup\@setsize\SetFigFont{##2pt}}%
  \expandafter\x
    \csname \romannumeral\the\count@ pt\expandafter\endcsname
    \csname @\romannumeral\the\count@ pt\endcsname
  \csname ##3\endcsname}%
\fi
\fi\endgroup
\begin{picture}(4749,2002)(739,-1544)
\put(2626,-1486){\makebox(0,0)[lb]{\smash{\small{$\gamma$}%
}}}
\put(5251,-1486){\makebox(0,0)[lb]{\smash{\small{$\gamma$}%
}}}
\put(4831,-1148){\makebox(0,0)[lb]{\smash{\small{$\xi_1$}%
}}}
\put(3841,-1156){\makebox(0,0)[lb]{\smash{\small{$\xi_0$}%
}}}
\put(4365,-961){\makebox(0,0)[lb]{\smash{\small{$x_1$}%
}}}
\put(4358,239){\makebox(0,0)[lb]{\smash{\small{$x_0$}%
}}}
\put(1733,239){\makebox(0,0)[lb]{\smash{\small{$x_0$}%
}}}
\put(1613,-1126){\makebox(0,0)[lb]{\smash{\small{$\xi_1$}%
}}}
\put(751,314){\makebox(0,0)[lb]{\smash{\small{(a)}%
}}}
\put(3451,314){\makebox(0,0)[lb]{\smash{\small{(b)}%
}}}
\put(1553,-631){\makebox(0,0)[lb]{\smash{\small{$y$}%
}}}
\put(1224,-343){\makebox(0,0)[lb]{\smash{\small{$P_i$}%
}}}
\put(2211,-363){\makebox(0,0)[lb]{\smash{\small{$P_{i+1}$}%
}}}
\put(3849,-343){\makebox(0,0)[lb]{\smash{\small{$P_i$}%
}}}
\put(4836,-363){\makebox(0,0)[lb]{\smash{\small{$P_{i+1}$}%
}}}
\end{picture}
\end{center}}

\newcommand{\drawFCclaimtwoa}{\begin{center}
\begin{picture}(0,0)%
\includegraphics{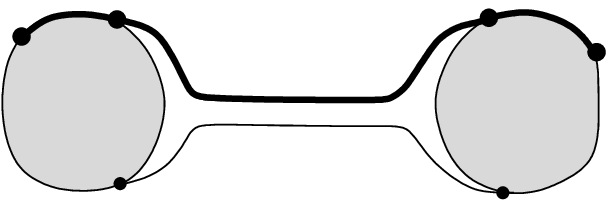}%
\end{picture}%
\setlength{\unitlength}{3947sp}%
\begingroup\makeatletter\ifx\SetFigFont\undefined
\def\x##1##2##3##4##5##6##7\relax{\def\x{##1##2##3##4##5##6}}%
\expandafter\x\fmtname xxxxxx\relax \def\y{splain}%
\ifx\x\y   % LaTeX or SliTeX?
\gdef\SetFigFont##1##2##3{%
  \ifnum ##1<17\tiny\else \ifnum ##1<20\small\else
  \ifnum ##1<24\normalsize\else \ifnum ##1<29\large\else
  \ifnum ##1<34\Large\else \ifnum ##1<41\LARGE\else
     \huge\fi\fi\fi\fi\fi\fi
  \csname ##3\endcsname}%
\else
\gdef\SetFigFont##1##2##3{\begingroup
  \count@##1\relax \ifnum 25<\count@\count@25\fi
  \def\x{\endgroup\@setsize\SetFigFont{##2pt}}%
  \expandafter\x
    \csname \romannumeral\the\count@ pt\expandafter\endcsname
    \csname @\romannumeral\the\count@ pt\endcsname
  \csname ##3\endcsname}%
\fi
\fi\endgroup
\begin{picture}(2912,1128)(1854,-1423)
\put(3249,-856){\makebox(0,0)[lb]{\smash{\small{$\xi$}%
}}}
\put(2116,-451){\makebox(0,0)[lb]{\smash{\small{$\rho_{r-1}$}%
}}}
\put(4486,-466){\makebox(0,0)[lb]{\smash{\small{$\rho_{s+1}$}%
}}}
\put(2183,-1006){\makebox(0,0)[lb]{\smash{\small{$P_{r-1}$}%
}}}
\put(4283,-1014){\makebox(0,0)[lb]{\smash{\small{$P_{s+1}$}%
}}}
\end{picture}
\end{center}}

\newcommand{\drawFCclaimtwob}{\begin{center}
\begin{picture}(0,0)%
\includegraphics{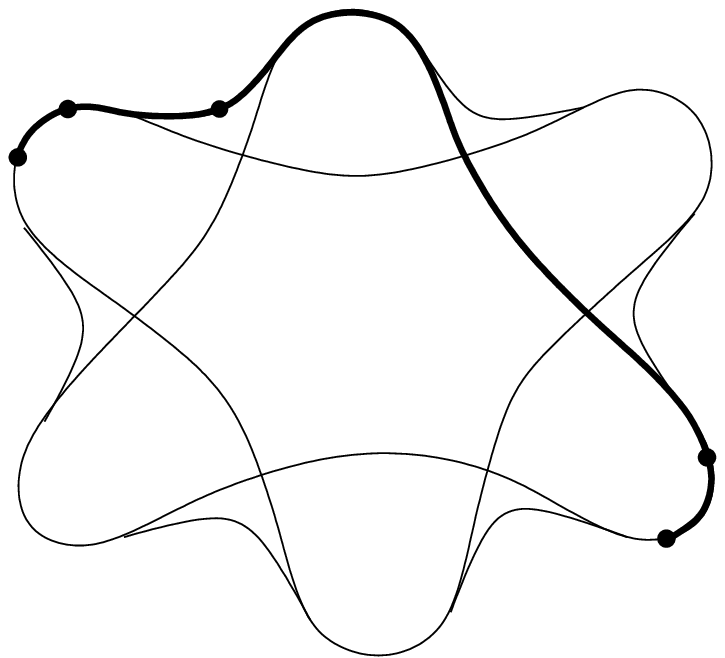}%
\end{picture}%
\setlength{\unitlength}{3947sp}%
\begingroup\makeatletter\ifx\SetFigFont\undefined
\def\x##1##2##3##4##5##6##7\relax{\def\x{##1##2##3##4##5##6}}%
\expandafter\x\fmtname xxxxxx\relax \def\y{splain}%
\ifx\x\y   % LaTeX or SliTeX?
\gdef\SetFigFont##1##2##3{%
  \ifnum ##1<17\tiny\else \ifnum ##1<20\small\else
  \ifnum ##1<24\normalsize\else \ifnum ##1<29\large\else
  \ifnum ##1<34\Large\else \ifnum ##1<41\LARGE\else
     \huge\fi\fi\fi\fi\fi\fi
  \csname ##3\endcsname}%
\else
\gdef\SetFigFont##1##2##3{\begingroup
  \count@##1\relax \ifnum 25<\count@\count@25\fi
  \def\x{\endgroup\@setsize\SetFigFont{##2pt}}%
  \expandafter\x
    \csname \romannumeral\the\count@ pt\expandafter\endcsname
    \csname @\romannumeral\the\count@ pt\endcsname
  \csname ##3\endcsname}%
\fi
\fi\endgroup
\begin{picture}(3616,3131)(879,-2501)
\put(3586,-465){\makebox(0,0)[lb]{\smash{\small{$\pi$}%
}}}
\put(4471,-1921){\makebox(0,0)[lb]{\smash{\small{$\rho_{s+1}$}%
}}}
\put(879,120){\makebox(0,0)[lb]{\smash{\small{$\rho_{r-1}$}%
}}}
\put(1658,171){\makebox(0,0)[lb]{\smash{\small{$\sigma$}%
}}}
\put(2078,-46){\makebox(0,0)[lb]{\smash{\small{$\tau$}%
}}}
\end{picture}
\end{center}}

\begin{document}

\title{Nonpositively curved 2-complexes with isolated flats}

\author{G Christopher Hruska}
\address{Department of Mathematics, University of Chicago\\
5734 S University Ave, Chicago, IL 60637, USA}
\email{chruska@math.uchicago.edu}

\primaryclass{%  AMS 2000 Math subject classifications
20F67} % Hyperbolic groups and nonpositively curved groups
\secondaryclass{%
20F06, % Cancellation theory; application of van Kampen diagrams
57M20} % Two-dimensional complexes

\keywords{Word hyperbolic, nonpositive curvature, thin triangles,
quasigeodesics, isolated flats}

\begin{abstract}
We introduce the class of nonpositively curved $2$--complexes
with the Isolated Flats Property.  These $2$--complexes are, in a sense,
hyperbolic relative to their flats.
More precisely, we show that several important properties of Gromov-hyperbolic
spaces hold ``relative to flats'' in
nonpositively curved $2$--complexes with the Isolated Flats Property.

We introduce the Relatively Thin Triangle Property, which states roughly
that the fat part of a geodesic triangle lies near a single flat.
We also introduce the Relative Fellow Traveller Property,
which states that pairs of quasigeodesics with common endpoints
fellow travel relative to flats, in a suitable sense.
The main result of this paper states that in the setting of $\CAT(0)$
$2$--complexes, the Isolated Flats Property is equivalent to the Relatively
Thin Triangle Property and is also equivalent to the Relative Fellow
Traveller Property.
\end{abstract}

\asciiabstract{%
We introduce the class of nonpositively curved 2-complexes
with the Isolated Flats Property.  These 2-complexes are, in a sense,
hyperbolic relative to their flats.
More precisely, we show that several important properties of Gromov-hyperbolic
spaces hold `relative to flats' in
nonpositively curved 2-complexes with the Isolated Flats Property.

We introduce the Relatively Thin Triangle Property, which states roughly
that the fat part of a geodesic triangle lies near a single flat.
We also introduce the Relative Fellow Traveller Property,
which states that pairs of quasigeodesics with common endpoints
fellow travel relative to flats, in a suitable sense.
The main result of this paper states that in the setting of CAT(0)
2-complexes, the Isolated Flats Property is equivalent to the Relatively
Thin Triangle Property and is also equivalent to the Relative Fellow
Traveller Property.}

\maketitlepage

%%%%%%%%%%%%%%%%%%%%%%%%%%%%%%%%%%%%%%%%%%%%%%%%%%%%%%%%%%%%%%%%%%%%%%%%%%%
\section{Introduction}
\label{sec:Introduction}
%%%%%%%%%%%%%%%%%%%%%%%%%%%%%%%%%%%%%%%%%%%%%%%%%%%%%%%%%%%%%%%%%%%%%%%%%%%

The theory of $\delta$--hyperbolic spaces has been enormously fruitful since
it was first introduced by Gromov in his seminal article \cite{Gromov87}.
In that article, Gromov establishes two facts about
$\delta$--hyperbolic spaces which are used heavily in the proofs
of many results in the theory.
The first is that geodesic triangles
are thin, and the second is that quasigeodesics with common endpoints
(asynchronously) fellow travel.

In this article, we introduce nonpositively curved spaces
with the \emph{Isolated Flats Property}, generalizing $\delta$--hyperbolic
spaces.
Spaces with the Isolated Flats Property can be studied using techniques
analogous to those used in the study of $\delta$--hyperbolic spaces.
In fact, many results about $\delta$--hyperbolic spaces have natural
extensions to this new setting.
The resulting theory shares much of the robust
character of the $\delta$--hyperbolic setting.
In contrast, few methods are currently known for
extending results from $\delta$--hyperbolic spaces to arbitrary
nonpositively curved spaces.

The main result of this article provides a starting point for the
process of generalizing
results from the $\delta$--hyperbolic setting to the isolated flats
setting.  We introduce the Relatively Thin Triangle Property,
which extends the notion of thin triangles, and
the Relative Fellow Traveller Property, which generalizes the
fellow travelling of quasigeodesics.
The main theorem shows that these ``relative'' properties
are each equivalent to the Isolated Flats Property
in the $2$--dimensional setting.

\begin{thm}\label{thm:2dEquivalent}
Let $X$ be a proper, cocompact piecewise Euclidean $\CAT(0)$ $2$--complex.
The following are equivalent.
\begin{enumerate}
  \item $X$ has the Isolated Flats Property.
  \item $X$ has the Relatively Thin Triangle Property.
  \item $X$ has the Relative Fellow Traveller Property.
\end{enumerate}
\end{thm}

The implications (2)~$\implies$~(1) and (3)~$\implies$~(1) are
fairly straightforward.  The converse implications
(1)~$\implies$~(2) and (1)~$\implies$~(3) are more difficult
and also have thus far provided more applications.
Typically, the Isolated Flats Property is easier to detect
than the Relatively Thin Triangle Property and
the Relative Fellow Traveller Property.
However, the latter two properties are quite useful in applications.
For instance, one can start with an argument in the $\delta$--hyperbolic
setting which uses the thinness of triangles
or the fellow travelling of quasigeodesics
and try to convert it to an argument in the isolated flats setting
which uses the corresponding
relative properties.
In Subsection~\ref{subsec:Applications}, we discuss several applications
of this nature.

In the course of the proof of Theorem~\ref{thm:2dEquivalent},
we prove a ``quadratic divergence'' theorem for geodesic rays
in a $2$--complex with the Isolated Flats Property
(Proposition~\ref{prop:QuadraticDivergence}).  This theorem
can be interpreted as saying that given a pair of geodesic rays
in such a $2$--complex neither of which lingers very long near a single
flat, the given rays diverge from each other at a rate which is at least
quadratic.  This divergence theorem is of independent
interest because of its
similarity to the exponential divergence theorem for $\delta$--hyperbolic
spaces proved by Cooper and Mihalik in
\cite[Theorem~2.19]{ABC91}.

The exponential divergence theorem is a key ingredient of
Lustig and Mihalik's proof that quasigeodesics track close
to geodesics in hyperbolic spaces \cite[Proposition~3.3]{ABC91}.
Similarly, Proposition~\ref{prop:QuadraticDivergence} implies that
given any geodesic segment
which does not linger very long near a single flat, that geodesic tracks
close to any quasigeodesic connecting its endpoints (see
Section~\ref{sec:Divergence}).

It seems likely that the hypothesis that $X$ is $2$--dimensional can be
dropped from Theorem~\ref{thm:2dEquivalent}.
Some evidence for this conjecture is presented below in
Subsection~\ref{subsec:History}.  However,
the present methods are specific to the $2$--dimensional setting.
Our techniques depend heavily on the observation that a van~Kampen
diagram over a $\CAT(0)$ $2$-complex is itself a $\CAT(0)$ space.
In higher dimensions, completely new techniques would seem to be
necessary.

%%%%%%%%%%%%%%%%%%%%%%%%%%%%%%%%%%%%%%%%%%%%%%%%%%%%%%%%%%%%%%%%%%%%%%%%%%%
\subsection{The Isolated Flats Property}
\label{subsec:IsolatedFlats}
%%%%%%%%%%%%%%%%%%%%%%%%%%%%%%%%%%%%%%%%%%%%%%%%%%%%%%%%%%%%%%%%%%%%%%%%%%%

The Isolated Flats Property is defined precisely in
Section~\ref{sec:FlatTriplaneTheorem}.
Roughly speaking, a $\CAT(0)$ space
has the Isolated Flats Property if its
isometrically embedded flat Euclidean subspaces diverge from each other
in all directions, in the sense that their corresponding boundary spheres
at infinity are disjoint.  Note that the Isolated Flats Property
is vacuously satisfied in any $\delta$--hyperbolic $\CAT(0)$ space.

The prototypical example is the universal cover of a truncated
finite volume hyperbolic manifold~$M$.  Such a space,
called a \emph{neutered space} or \emph{core}, is obtained from
hyperbolic space~$\Hyp^n$ by removing a family of disjoint open horoballs
corresponding to the cusps of~$M$.  The neutered space (with the induced
path metric) is a $\CAT(0)$ space whose only flat subspaces are the
boundaries of the deleted horoballs, which are isolated.

The idea of studying spaces with isolated flats is implicit in work
of Kapovich--Leeb \cite{KapovichLeeb95} and of Wise
\cite{Wise96,WiseFigure8}, and has also been studied by Kleiner
(personal communication).
In unpublished work, Wise has proved a Flat Triplane Theorem,
which states that in the $2$--dimensional setting the Isolated
Flats Property is equivalent to an absence of isometrically embedded
triplanes.  A \emph{triplane} is a space obtained by gluing three
Euclidean half-planes together along their boundary lines.
We provide Wise's proof of the Flat Triplane Theorem
in Section~\ref{sec:FlatTriplaneTheorem}.

Wise observed in \cite[\S4.0]{Wise96} that if~$X$ is any compact
nonpositively curved $2$--complex
each of whose $2$--cells is isometric to a regular Euclidean hexagon
then the universal cover of~$X$ has the Isolated Flats Property.
The reason is that a triplane cannot be built out of regular hexagons.
In \cite{BallmannBrin94}, Ballmann and Brin give explicit techniques for
constructing $\CAT(0)$ hexagonal $2$--complexes with arbitrary local
data.  The following theorem due to Moussong
indicates some of the richness of this class of $2$--complexes.
The theorem is a special case of the main result of Moussong's
thesis (\cite{Moussong88}, see also Haglund \cite{Haglund91} and
Benakli \cite{Benakli94}).

\begin{thm}[Moussong]\label{thm:Moussong}
For any simplicial graph~$L$, there is a $\CAT(0)$ hexagonal
$2$--complex~$X$ such that the link of each vertex is isomorphic to the
graph~$L$.
If $L$ is finite, then the Coxeter group~$W$ defined by the graph~$L$
\textup{(}labeled with a $3$ on each edge\textup{)}
acts properly, cocompactly, and cellularly on~$X$ with a compact quotient.
\end{thm}

In the preceding theorem,
we use the convention that a simplicial graph~$L$
(with all edges labeled by the number~$3$)
defines a Coxeter system with one generator~$s_i$ of order two
for each vertex~$v_i$, and a relation $s_i s_j s_i = s_j s_i s_j$
whenever two vertices $v_i$ and $v_j$ are connected by an edge.
The $2$--complex arising in Theorem~\ref{thm:Moussong} is the
Davis--Moussong geometric realization of the Coxeter system given
by~$L$.

In addition to hexagonal complexes,
there are also many squared complexes with isolated flats.
For instance, if $L$ is a hyperbolic, prime, alternating link then
$\pi_1(S^3 \setminus L)$ is the fundamental group of a
nonpositively curved squared $2$--complex~$X$
whose universal cover has the Isolated Flats Property.
This squared $2$--complex was constructed by Dehn
\cite{Dehn87} and was shown to be nonpositively curved by
Weinbaum \cite{Weinbaum71}.  The Isolated Flats Property for these
$2$-complexes follows from \cite{HruskaGeometric}.

More generally, Aitchison has shown that every finite volume
cusped hyperbolic $3$--manifold deformation retracts onto a
compact $2$--complex that admits a
piecewise Euclidean metric with nonpositive curvature
\cite{Aitchison}.
The universal cover of this $2$--complex has the Isolated Flats Property by
\cite{HruskaGeometric}.
The $2$--complexes arising in Aitchison's
construction typically have irregularly shaped cells.

Wise encountered spaces with isolated flats while investigating
the $\Z \times \Z$ conjecture for $\CAT(0)$ groups.  This conjecture
states that, if a group acts properly and cocompactly
by isometries on a $\CAT(0)$ space, then
either the group is word hyperbolic or it contains a $\Z \times \Z$
subgroup.  The $\Z\times\Z$ conjecture has been proved by Bangert--Schroeder
in the case that the $\CAT(0)$ space is a real analytic manifold
\cite{BangertSchroeder91},
however the general conjecture seems quite difficult.
In fact, a theorem of Kari--Papasoglu \cite{KariPapasoglu99}
strongly indicates that the $\Z\times\Z$ conjecture may be false
even for $\CAT(0)$ squared complexes.

Wise noticed that, in the presence of the Isolated Flats Property,
the situation is much simpler, since
one can then prove that all flats are periodic.
Consequently the $\Z\times\Z$ conjecture is true in this setting
(\cite[Proposition~4.0.4]{Wise96}, for a complete proof see
\cite{HruskaGeometric}).

\begin{thm}[Wise]
Let $G$ act properly, cocompactly, and isometrically on a $\CAT(0)$
space
with the Isolated Flats Property.  Then either $G$ is word hyperbolic,
or $G$ contains a subgroup isomorphic to $\Z \times \Z$.
\end{thm}

The $\Z\times\Z$ conjecture was previously established by Ballmann--Brin
in the special case that the $\CAT(0)$ space is a hexagonal $2$-complex
\cite{BallmannBrin94}.

Sela has conjectured that the Isolated Flats Property is closely
related to limit groups (or fully residually free groups),
which arise in the study of the elementary theory of free groups.
In particular, he has conjectured that each limit group acts properly
and cocompactly on a $\CAT(0)$ space with the Isolated Flats Property
\cite{SelaProblems}.

%%%%%%%%%%%%%%%%%%%%%%%%%%%%%%%%%%%%%%%%%%%%%%%%%%%%%%%%%%%%%%%%%%%%%%%%%%%
\subsection{Hyperbolicity relative to flats}
\label{subsec:HyperbolicityRelFlats}
%%%%%%%%%%%%%%%%%%%%%%%%%%%%%%%%%%%%%%%%%%%%%%%%%%%%%%%%%%%%%%%%%%%%%%%%%%%

The Relatively Thin Triangle Property and
the Relative Fellow Traveller Property arise from an intuitive notion
that spaces with the Isolated Flats Property
are $\delta$--hyperbolic ``relative to flats.''

\begin{figure}[ht!]
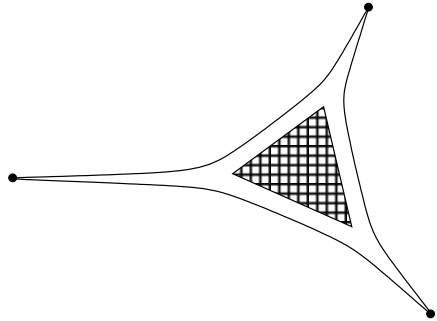

\drawthintriangle
\caption{A triangle which is $\delta$--thin relative to a flat}
\label{fig:ThinTriangle}
\end{figure}

For instance, a triangle is \emph{$\delta$--thin}
if each side lies in a $\delta$--neighborhood of the union of the other two
sides.
A geodesic space is \emph{$\delta$--hyperbolic} if
every geodesic triangle in the space is $\delta$--thin.
In Section~\ref{sec:HyperbolicityRelFlats}, we
introduce the notion of a geodesic triangle being
\emph{$\delta$--thin relative to a flat}.
The idea is that each side of the triangle lies in a $\delta$--neighborhood
of the union of the other two sides and some flat, as illustrated in
Figure~\ref{fig:ThinTriangle}.
A space has the Relatively Thin Triangle Property if there is
a constant~$\delta$ such that every geodesic triangle is either
$\delta$--thin in the standard sense or $\delta$--thin relative
to some flat.

\begin{figure}[ht!]
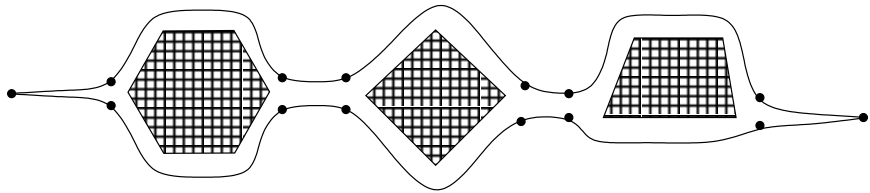

\drawrftp
\caption{A pair of paths which fellow travel relative to flats}
\label{fig:RFTP}
\end{figure}

A lemma due to Morse states that in the hyperbolic plane
any pair of quasigeodesics with common endpoints asynchronously fellow travel.
Theorems of Gromov and Masur--Minsky (\cite[Proposition~7.2.A]{Gromov87}
and \cite[Lemma~7.2]{MasurMinsky99})
together show that this fellow traveller property
is equivalent to $\delta$--hyperbolicity.
In Section~\ref{sec:HyperbolicityRelFlats}, we introduce the more general
notion of a pair of quasigeodesics that
\emph{fellow travel relative to flats}.
The idea is that the curves alternate between two types of behavior:
``tracking'' close together and travelling near a common flat,
as illustrated in Figure~\ref{fig:RFTP}.
A nonpositively curved space has the \emph{Relative Fellow Traveller
Property} if every pair of quasigeodesics with common endpoints fellow
travels relative to flats.

%%%%%%%%%%%%%%%%%%%%%%%%%%%%%%%%%%%%%%%%%%%%%%%%%%%%%%%%%%%%%%%%%%%%%%%%%%%
\subsection{Historical background}
\label{subsec:History}
%%%%%%%%%%%%%%%%%%%%%%%%%%%%%%%%%%%%%%%%%%%%%%%%%%%%%%%%%%%%%%%%%%%%%%%%%%%

Several results analogous to Theorem~\ref{thm:2dEquivalent}
have previously been shown for certain nonpositively curved manifolds.
In \cite{KapovichLeeb95}, Kapovich and Leeb considered a class
of nonpositively curved manifolds with isolated flats in which the flat
subspaces are separated by regions of strict negative curvature.
In that setting they proved that the Relatively Thin Triangle
Property holds.

Epstein showed that an analogue of the Relative Fellow Traveller Property
holds for real
hyperbolic space $\Hyp^n$ with a disjoint family of open horoballs
removed \cite[Theorem~11.3.1]{ECHLPT92}.
Lang generalized Epstein's result to manifolds with pinched sectional
curvature $-a^2 \le \kappa \le -1$ for $1 \le a < 2$
in \cite{Lang96}.
Farb proved a result in the same spirit using a slightly different metric
on the neutered space \cite[Lemma~4.5]{Farb98}.
Additionally, we remark that the Relative Fellow Traveller Property
is similar to Farb's bounded coset penetration property.
This similarity is exploited in \cite{HruskaRelHyp}
in the proof of the relative hyperbolicity
theorem for groups acting on spaces with isolated flats.

Kapovich--Leeb and Epstein's results mentioned above provide additional
examples of spaces with the Isolated Flats Property
which also have either the Relatively Thin Triangle Property
or the Relative Fellow Traveller Property.  In light of these results,
it seems likely that the hypothesis that $X$ is $2$--dimensional
can be dropped from Theorem~\ref{thm:2dEquivalent}.

%%%%%%%%%%%%%%%%%%%%%%%%%%%%%%%%%%%%%%%%%%%%%%%%%%%%%%%%%%%%%%%%%%%%%%%
\subsection{Applications of Theorem~\ref{thm:2dEquivalent}}
\label{subsec:Applications}
%%%%%%%%%%%%%%%%%%%%%%%%%%%%%%%%%%%%%%%%%%%%%%%%%%%%%%%%%%%%%%%%%%%%%%%

Thus far, the main application of the Relatively Thin Triangle Property
and the Relative Fellow Traveller Property has been to extend
results from the $\delta$--hyperbolic setting to the setting of $\CAT(0)$
spaces with the Isolated Flats Property.  By way of example,
we list below several immediate consequences
of Theorem~\ref{thm:2dEquivalent} whose proofs make use of
either the Relatively Thin Triangle Property
or the Relative Fellow Traveller Property.

Suppose a group $G$ acts properly and cocompactly by isometries on
a $\CAT(0)$ space~$X$.
In \cite{HruskaGeometric}, the author shows that in the presence of
the Isolated Flats Property and the Relative Fellow Traveller Property
the boundary at infinity of~$X$ is an invariant of the group~$G$.
Recall that the \emph{boundary} $\boundary X$ of~$X$
is the space of geodesic rays emanating from a fixed basepoint
with the compact-open topology.
Together with Theorem~\ref{thm:2dEquivalent} we have the following
consequence in the $2$--dimensional setting.

\begin{thm}[Boundary is well-defined]\label{thm:BoundaryWellDefined}
Suppose a group~$G$ acts properly and cocompactly by isometries
on two $\CAT(0)$ spaces $X$ and~$Y$.
Suppose further that $X$ is a piecewise Euclidean $2$--complex
with the Isolated Flats Property.
Then any equivariant quasi-isometry $X \to Y$ induces an
equivariant homeomorphism $\boundary X \to \boundary Y$.
\end{thm}

The previous theorem was established in the word hyperbolic setting by
Gromov in \cite{Gromov87}.  Croke and Kleiner showed in \cite{CrokeKleiner00}
that this theorem does not extend to the general $\CAT(0)$ setting
by constructing two homeomorphic nonpositively curved
$2$--complexes whose universal covers have nonhomeomorphic boundaries.
Julia Wilson has since shown that Croke and Kleiner's
construction actually provides
a continuous family of homeomorphic $2$--complexes whose universal covers
are pairwise nonhomeomorphic \cite{WilsonBoundary}.%
\footnote{The author has been informed that Kleiner has unpublished
work related to the article \cite{CrokeKleiner02}
in which he proves Theorem~\ref{thm:BoundaryWellDefined}
without the $2$--dimensional hypothesis.
}

A second consequence of Theorem~\ref{thm:2dEquivalent} deals with
quasiconvex subgroups.
A subspace $Z$ of~$X$ is \emph{quasiconvex} if there is a constant~$\epsilon$
so that every geodesic in~$X$ connecting two points of~$Z$ lies inside
an $\epsilon$--neighborhood of~$Z$.  If $\rho$ is a proper, cocompact
action of a group $G$ by isometries on a $\CAT(0)$ space~$X$, then
a subgroup $H \le G$ is \emph{quasiconvex} with respect to~$\rho$
if for some $x \in X$, the orbit $Hx$ is a quasiconvex subspace of~$X$.

In \cite{HruskaGeometric} the author shows that in the presence of
the Isolated Flats Property and the Relative Fellow Traveller Property
the notion of a subgroup $H\le G$ being quasiconvex
does not depend on the choice of
action~$\rho$ or on the choice of $\CAT(0)$ space~$X$.
Together with Theorem~\ref{thm:2dEquivalent}, we have the following
consequence for $2$--dimensional complexes.

\begin{thm}[Quasiconvexity is well-defined]\label{thm:QCWellDefined}
Let $\rho$ and~$\sigma$ be two proper, cocompact actions
of a group $G$ by isometries on $\CAT(0)$ spaces $X$ and~$Y$.
Suppose further that $X$ is a piecewise Euclidean $2$--complex
with the Isolated Flats Property.
For each subgroup $H\le G$, the following are equivalent.
\begin{enumerate}
  \item $H$ is quasiconvex with respect to~$\rho$.
  \item $H$ is quasiconvex with respect to~$\sigma$.
  \item The inclusion $H \inclusion G$ is a quasi-isometric embedding.
\end{enumerate}
\end{thm}

The previous theorem was established by Short
for word hyperbolic groups
in \cite{Short91}.  The author shows in \cite{HruskaGeometric}
that this result does not extend to the general $\CAT(0)$ setting.

In a subsequent article \cite{HruskaRelHyp}, the author will
use Theorem~\ref{thm:2dEquivalent}
to prove the following result,
which provides a precise group theoretic manifestation of the intuitive
notion that
a space with isolated flats is hyperbolic ``relative to flats.''

\begin{thm}[Isolated Flats $\implies$ Relatively Hyperbolic]
\label{thm:RelHyp}
Suppose a group $G$ acts properly and cocompactly by isometries on a
$\CAT(0)$ $2$--complex with the Isolated Flats
Property.
Then $G$ is hyperbolic relative to the collection of
maximal virtually abelian subgroups of rank two.
\end{thm}

The previous theorem is stated using Gromov and Bowditch's terminology for
relative hyperbolicity \cite{Gromov87,BowditchRelHyp}.
Using the terminology of Farb \cite{Farb98}, the conclusion is that
the group in question
is relatively hyperbolic with bounded coset penetration.
The proof of Theorem~\ref{thm:RelHyp} uses the fact that
a $2$--complex with isolated flats has the Relatively Thin Triangle
Property.

Theorem~\ref{thm:RelHyp} together with a result of Rebbechi
\cite{Rebbechi01} has the following immediate consequence.

\begin{thm}[Isolated Flats $\implies$ Biautomatic]
\label{thm:Biautomatic}
Suppose $G$ acts properly and cocompactly by isometries on a $\CAT(0)$
$2$--complex with the Isolated Flats Property.
Then $G$ is biautomatic.
\end{thm}

It is unknown whether a group acting properly and cocompactly
on an arbitrary $\CAT(0)$ space is necessarily biautomatic
(or even automatic).
Previously, biautomaticity has been proven only
for $\CAT(0)$ complexes built from a small number of allowed shapes
of polyhedral cells.
For instance, Gersten--Short establish biautomaticity for $2$--complexes
built of squares and three shapes of triangles
(specifically the $2$--complexes of type $A_1 \times A_1$, $A_2$, $B_2$,
and~$G_2$)
in \cite{GerstenShort90Automatic,GerstenShort91Automatic}.
Niblo--Reeves establish biautomaticity for $\CAT(0)$ cube complexes
in \cite{NibloReeves98}.
The main difference between Theorem~\ref{thm:Biautomatic}
and these previous results is that
Theorem~\ref{thm:Biautomatic} allows convex polygonal cells of
arbitrary shapes.

Finally we mention that Theorem~\ref{thm:RelHyp} and a result of Tukia
\cite{Tukia94} show that the Tits Alternative holds
for $2$--complexes with isolated flats.

\begin{thm}[Isolated Flats $\implies$ Tits Alternative]
Suppose $G$ acts properly and cocompactly on a $\CAT(0)$ $2$--complex
with the Isolated Flats Property.  Then $G$ satisfies the Tits Alternative.
In other words, every subgroup of~$G$ is either virtually abelian
or contains a free subgroup of rank two.
\end{thm}

Again it is unknown whether the previous theorem remains true
if the Isolated Flats Property is dropped from the hypothesis.

If the hypothesis that $X$ is $2$--dimensional can be dropped from
Theorem~\ref{thm:2dEquivalent}, as conjectured in the previous subsection,
then the results
listed above will apply to any $\CAT(0)$ space with the Isolated Flats
Property.

%%%%%%%%%%%%%%%%%%%%%%%%%%%%%%%%%%%%%%%%%%%%%%%%%%%%%%%%%%%%%%%%%%%%%%%%%%%
\subsection{Summary of the sections}
\label{subsec:Summary}
%%%%%%%%%%%%%%%%%%%%%%%%%%%%%%%%%%%%%%%%%%%%%%%%%%%%%%%%%%%%%%%%%%%%%%%%%%%

We begin with a few sections of background material.
In Section~\ref{sec:CAT(0)} we define $\CAT(0)$ spaces,
and review several of their important properties.
In Section~\ref{sec:PEComplexes}, we discuss piecewise Euclidean complexes
and their relation to nonpositive curvature by way of the Link Condition.
In Section~\ref{sec:Diagrams}, we review definitions and basic results
about diagrams, reduced maps, and the Combinatorial Gauss--Bonnet theorem,
following the development of McCammond--Wise \cite{McCammondWise02}.

In Section~\ref{sec:FlatTriplaneTheorem}, we give a definition of the
Isolated Flats Property which is catered to the $2$--dimensional setting.
We also give Wise's proof of the Flat Triplane Theorem, which has not
previously appeared in the literature.
In Section~\ref{sec:HyperbolicityRelFlats}, we state the Relatively
Thin Triangle Property and the Relative Fellow Traveller Property.
In the $2$--dimensional setting, we use
the Flat Triplane Theorem to show that
each of these properties implies the Isolated Flats Property,
thus establishing the implications
(2)~$\implies$~(1) and (3)~$\implies$~(1)
of Theorem~\ref{thm:2dEquivalent}.

The reverse implications (1)~$\implies$~(2) and (1)~$\implies$~(3)
are more difficult and require a deeper analysis
of the structure of disc diagrams corresponding to $2$--complexes with isolated flats.
In Section~\ref{sec:RuffledBoundaries}, we introduce the notion of a
diagram which is \emph{ruffled} along a boundary path.  Such a diagram
has the property that every point of the given boundary path is either
close to a vertex with strictly negative curvature or close to some other
part of the boundary.  This notion generalizes the fact that
in the negatively curved setting the fat part of
any reduced disc diagram is filled with
vertices at which the curvature is negative.
We use the Combinatorial Gauss--Bonnet Theorem to relate ruffled diagrams
to the usual notion of $\delta$--thin triangles.

In Section~\ref{sec:Preflats}, we study \emph{preflats} in reduced disc
diagrams.  Preflats in a diagram $D \to X$ correspond to flats in~$X$.
We prove two results about preflats, which play a key role in the proof
of Theorem~\ref{thm:2dEquivalent}.
The first is Proposition~\ref{prop:AlongGeodesic}, which states roughly that
if a geodesic segment occurs in the boundary path of a nonpositively curved
disc diagram, then either the diagram is ruffled along the geodesic or
the geodesic is close to a large preflat.
The second result is Proposition~\ref{prop:AroundFlat}
which states that in the presence of the Isolated Flats Property,
preflats in disc diagrams are surrounded by ruffles.

In Section~\ref{sec:RelThinTriangles}, we use the ideas developed in
the previous two sections to prove the implication
(1)~$\implies$~(2) of Theorem~\ref{thm:2dEquivalent}.
The idea of the proof is that a triangular diagram
either contains a large preflat surrounded by ruffles, or it is
ruffled along one side.  In the first case, the triangle is thin relative
to a flat, and in the second case, the triangle is thin in the standard
sense.

In Section~\ref{sec:Divergence} we turn our attention to the remaining
implication (1)~$\implies$~(3) of Theorem~\ref{thm:2dEquivalent}.
In this section we prove the ``quadratic divergence'' theorem
mentioned above.
Cooper, Lustig, and Mihalik use an analogous exponential divergence
theorem in their proof that quasigeodesics track close to geodesics
in a $\delta$--hyperbolic space \cite[Proposition~3.3]{ABC91}.
The main result of Section~\ref{sec:Divergence} is a similar fellow travelling
result for $2$--complexes with isolated flats in the presence of ruffles.

In order to prove (1)~$\implies$~(3), it will be useful to understand the
structure of the convex hull of a union of preflats in a disc diagram.
In Section~\ref{sec:ConvexHull}, we examine in detail the convex hull
of a union of two preflats and also the convex hull of the union of one
preflat and a point.  We show that in the presence of
the Isolated Flats Property, the convex hull
of the two objects in question is essentially just the union of those objects
together with a path of shortest length connecting them.

Section~\ref{sec:FlatClosure} contains the main part of the proof of
(1)~$\implies$~(3).  In the isolated flats setting, we study disc diagrams
with a geodesic segment along the boundary.  We consider the convex hull of
the union of
that geodesic and all the preflats that come near the geodesic.
We show that the boundary of this convex hull contains a path which fellow
travels the geodesic relative to flats.  We also show that this convex hull
is surrounded by ruffles in the diagram.

In Section~\ref{sec:RelativeFTP} we complete the proof of
Theorem~\ref{thm:2dEquivalent}.  We consider a disc diagram whose boundary
consists of a geodesic and quasigeodesic.  Once one removes the convex hull
constructed in the previous section, one obtains a diagram which
satisfies the hypothesis of the ruffled fellow travelling result from
Section~\ref{sec:Divergence}.  We conclude that
any geodesic and quasigeodesic with common endpoints fellow travel
relative to flats.

Finally, we establish the Relative Fellow Traveller Property
by considering the general case of a pair of quasigeodesics with
common endpoints.  We derive the general case from the special case
by comparing each quasigeodesic with the geodesic connecting its endpoints
and piecing together the two sequences of flats associated to each such
pair.

%%%%%%%%%%%%%%%%%%%%%%%%%%%%%%%%%%%%%%%%%%%%%%%%%%%%%%%%%%%%%%%%%%%%%%%%%%%
\subsection{Acknowledgements}
\label{subsec:Acknowledgements}
%%%%%%%%%%%%%%%%%%%%%%%%%%%%%%%%%%%%%%%%%%%%%%%%%%%%%%%%%%%%%%%%%%%%%%%%%%%

The results in this article were originally published as part of my PhD
dissertation at Cornell University.  That dissertation was prepared under the
guidance of Daniel Wise and Karen Vogtmann.
I would like to thank Karen for supporting my desire to work with Dani on this
project, although he was a postdoc at the time.
I would also like to thank
Dani for pointing me in the direction of such a fruitful research area,
and for devoting such a large amount of time to encouraging this research even
after he moved from Ithaca to Boston and eventually to Montreal.
On many occasions Dani and his family graciously hosted me in their home so that
I could talk mathematics with Dani for long hours at a stretch.
I also benefited from conversations about this research with Martin
Bridson, Marshall Cohen, Ilya Kapovich,
Bruce Kleiner, Jon McCammond, John Meier and surely
others I have forgotten to mention.
Additionally, I would like to thank the referee for numerous helpful
suggestions that I hope have improved the exposition of this
article.

This research was partially supported by a grant from the National Science
Foundation.

%%%%%%%%%%%%%%%%%%%%%%%%%%%%%%%%%%%%%%%%%%%%%%%%%%%%%%%%%%%%%%%%%%%%%%%%%%%
\section{$\CAT(0)$ spaces}
\label{sec:CAT(0)}
%%%%%%%%%%%%%%%%%%%%%%%%%%%%%%%%%%%%%%%%%%%%%%%%%%%%%%%%%%%%%%%%%%%%%%%%%%%

In this section, we review some definitions and several well-known results
about $\CAT(0)$ spaces.
We refer the reader to Bridson--Haefliger \cite{BH99}
for proofs of the results listed in this section, as well as for
historical information about the origins of these ideas.
We give precise theorem numbers for the corresponding statements in \cite{BH99}
when they do not have a common name that can be found in the index of that
book.

A \emph{geodesic} in a metric space~$X$ is an isometric embedding of
an interval of~$\R$ into~$X$.
A metric space~$X$ is \emph{geodesic} if every pair of points in~$X$
is connected by at least one geodesic.
If $p$ and~$q$ are points in a geodesic space,
we use the notation $[p,q]$ to denote a particular choice of
geodesic segment connecting the points $p$ and~$q$.

\begin{defn}[$\CAT(0)$]\label{def:CAT0}
Let $X$ be a geodesic metric space.  A \emph{geodesic triangle}
$\Delta(p,q,r)$ in~$X$
is the union of three geodesic segments $[p,q]$, \ $[q,r]$, and $[p,r]$
in~$X$.
A \emph{comparison triangle} $\bar\Delta
= \bar\Delta(\bar{p}, \bar{q}, \bar{r})$ for $\Delta$ is a triangle
in the Euclidean plane~$\E^2$ such that
\[
   d(p,q) = d(\bar{p},\bar{q}), \quad
   d(q,r) = d(\bar{q},\bar{r}), \quad \text{and} \quad
   d(p,r) = d(\bar{p},\bar{r}),
\]
as illustrated in Figure~\ref{fig:CAT0}.
A point $\bar{x} \in [\bar{p},\bar{q}]$ is a \emph{comparison point}
for $x \in [p,q]$ provided that $d(\bar{p},\bar{x}) = d(p,x)$.
Comparison points are defined similarly for points on the other sides
$[q,r]$ and $[p,r]$.
\begin{figure}[ht!]
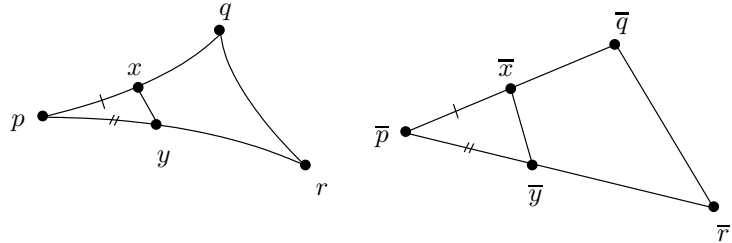

\drawCATzero
\caption[The $\CAT(0)$ inequality.]{A comparison triangle and a pair of
comparison points for the $\CAT(0)$ inequality}
\label{fig:CAT0}
\end{figure}

Let $\Delta$ be a triangle in $X$, and let $\bar\Delta$ be a comparison
triangle for $\Delta$.  We say that $\Delta$ satisfies the \emph{$\CAT(0)$
inequality} if for all points $x,y \in \Delta$ and
comparison points $\bar{x},\bar{y} \in \bar\Delta$ we have
\[
   d(x,y) \le d(\bar{x}, \bar{y}).
\]
If every geodesic triangle in $X$ satisfies the $\CAT(0)$ inequality, then
$X$ is called a \emph{$\CAT(0)$ space}.
The space~$X$ has \emph{nonpositive curvature}
if for each $x \in X$ there is an $\epsilon >0$ so that
$\ball{x}{\epsilon}$ is $\CAT(0)$.
\end{defn}

\begin{defn}[Angles]\label{def:angles}
Let $X$ be a nonpositively curved space,
and let $\gamma \co$\break $ [0,s]\to X$ and
$\gamma'\co [0,s']\to X$ be geodesic segments
with $p=\gamma(0)=\gamma'(0)$.
The \emph{angle} $\angle(\gamma,\gamma')$ between $\gamma$~and~$\gamma'$
is defined by the formula
\[
  \angle(\gamma,\gamma') = \lim_{t\to 0} 2\arcsin \frac{1}{2t}
  d \bigl( \gamma(t), \gamma'(t) \bigr).
\]
The above limit exists by \cite[Proposition~II.3.1]{BH99}.
When $X$ is $\CAT(0)$, then geodesic segments are uniquely determined by their
endpoints.  In this case,
if $x\ne p$ and $y\ne p$, then the angle between the segments
$[p,x]$ and $[p,y]$ will frequently be denoted $\angle_p(x,y)$.
\end{defn}

Note that the preceding definition agrees
with the usual Euclidean angle between geodesic segments
in the Euclidean plane~$\E^2$.

\begin{defn}[Comparison angles]
Let $\Delta(p,q,r)$ be a geodesic triangle in a $\CAT(0)$ space and
let $\bar\Delta(\bar p, \bar q, \bar r)$ be a comparison triangle for~$\Delta$.
The interior angle of $\bar\Delta$ at $\bar{p}$, denoted
${\bar\angle}_p(q,r)$, is called the \emph{comparison angle}
between $q$ and~$r$ at~$p$.
\end{defn}

\begin{thm}[\cite{BH99}, II.1.7(4)]\label{thm:comparisonangle}
Let $X$ be a $\CAT(0)$ space.  The angle between two sides of a geodesic
triangle $\Delta$ in~$X$ with distinct vertices is no greater than the
corresponding comparison angle in $\bar\Delta$.
\end{thm}

\begin{thm}[Convexity of the $\CAT(0)$ metric]\label{thm:CAT0convexity}
Let $\gamma$~and~$\gamma'$ be two geodesic segments in a $\CAT(0)$
space~$X$, each parametrized from $0$ to~$1$
proportional to arclength.  Then for each $t\in [0,1]$ we have
\[
  d \bigl( \gamma(t), \gamma'(t) \bigr) \le
   (1-t)\, d \bigl( \gamma(0), \gamma'(0) \bigr)
   + t \, d \bigl( \gamma(1),\gamma'(1) \bigr) .
\]
\end{thm}

\begin{thm}[Cartan--Hadamard Theorem]\label{thm:CartanHadamard}
Let $X$ be a complete metric space which is connected and simply connected.
If $X$ has nonpositive curvature, then $X$ is a $\CAT(0)$ space.
\end{thm}

A map $f\co Y \to X$ between metric spaces is
a \emph{local isometry} if for every $y\in Y$
there is an $\epsilon > 0$
such that the restriction of~$f$ to $\ball{y}{\epsilon}$ is an isometry
onto its image.  A \emph{local geodesic} in a metric space~$X$
is a local isometry from an interval of~$\R$ into~$X$.

The following two results are consequences of the Cartan--Hadamard Theorem.

\begin{thm}[\cite{BH99}, II.4.13]
Fix two points $x_0$ and~$x_1$ in a complete, nonpositively curved
metric space~$X$.
Then any homotopy class of paths between
$x_0$ and~$x_1$ in~$X$ is represented by a unique local geodesic.
\end{thm}

\begin{thm}[\cite{BH99}, II.4.14]\label{thm:localisometry}
Let $f\co Y\to X$ be a local isometry between two complete, connected
nonpositively curved spaces.
Then any lift of\/~$f$ to a map
$\tilde{f} \co \tilde{Y} \to \tilde{X}$
between their universal covers
is an isometric embedding.
\end{thm}

%%%%%%%%%%%%%%%%%%%%%%%%%%%%%%%%%%%%%%%%%%%%%%%%%%%%%%%%%%%%%%%%%%%%%%%%%%
\section{Piecewise Euclidean $2$--complexes}
\label{sec:PEComplexes}
%%%%%%%%%%%%%%%%%%%%%%%%%%%%%%%%%%%%%%%%%%%%%%%%%%%%%%%%%%%%%%%%%%%%%%%%%%

In this section, we collect definitions and basic results about
piecewise Euclidean $2$--complexes.  For a more thorough
development of these ideas see \cite{BH99}.

\begin{defn}
A \emph{convex Euclidean polyhedron}~$P$ is the convex hull
of a finite set of
points in Euclidean space~$\E^n$.  The \emph{dimension} of~$P$ is
the minimal
dimension of an affine subspace $E\subseteq \E^n$ containing~$P$.
If $P$ is contained in one of the
closed half-spaces determined by some hyperplane $H \subset \E^n$,
then $H\cap P$ is called a \emph{face} of~$P$.  We also
consider $P$ itself to be a face of~$P$. Note that a face~$F$ of~$P$ is itself
a
convex polyhedron.
The \emph{vertices} of~$P$ are its $0$--dimensional faces.
\end{defn}

\begin{defn}
A \emph{piecewise Euclidean complex} is a
complex~$X$ formed from a disjoint union of convex
Euclidean polyhedra by gluing isometric faces by isometries.
A \emph{metric graph} is a one dimensional piecewise Euclidean complex.
\end{defn}

Note that a piecewise Euclidean complex $X$ has a natural cell structure
whose cells are the polyhedra of~$X$.
A piecewise Euclidean complex has a natural pseudometric
where the distance between two points is the infimum of the lengths of
paths connecting them.  According to the following theorem due to Bridson,
in many complexes of interest this infimum is actually realized by
a geodesic path.

\begin{thm}{\rm\cite{Bridson91}}\qua
If a connected, piecewise Euclidean complex~$X$
has only finitely many isometry types of cells, then $X$~is a complete
geodesic metric space.
\end{thm}

\begin{defn}[Link]
Let $X$~be a piecewise Euclidean $2$--complex, and let $v$~be a vertex of~$X$.
The \emph{link} $\Lk(v,X)$ of~$v$ in~$X$ is the space of all germs
of geodesics
in~$X$ based at~$v$.
The \emph{radial projection} onto $\Lk(v,X)$ is the function
\[
   \pi \co X \setminus \{v\} \to \Lk(v,X)
\]
sending each point~$x$
to the germ of the geodesic $[v,x]$.
If $X$ has only finitely many isometry types of cells, then
by \cite[Theorem~I.7.39]{BH99}, the function~$d$ given by
\[
   d \bigl( \pi(x), \pi(y) ) = \angle_v(x,y)
\]
is a metric on $\Lk(v,X)$.
Under this metric each link $\Lk(v,X)$ has a natural structure as a metric
graph with one edge for each corner of a $2$--cell
incident at~$v$.  The length of each edge is equal to the
angle of the corresponding corner.
\end{defn}

\begin{defn}
A \emph{locally geodesic loop of length~$L$} in a metric space~$X$
is a local isometry $C\to X$ where $C$ is a metric circle
of length~$L$.
\end{defn}

\begin{defn}[Link condition]\label{def:LinkCondition}
A piecewise Euclidean $2$--complex $X$ satisfies the \emph{link condition}
if for every vertex $v\in X^{(0)}$, every locally geodesic
loop $C\to \Lk(v,X)$ has length at least~$2\pi$.
\end{defn}

The following theorem is due to Gromov \cite[$\S$4.2]{Gromov87}.
Ballmann provided a proof in the locally finite case
\cite{Ballmann90}, and Bridson in the general case \cite{Bridson91}.

\begin{thm}\label{thm:LinkCondition}
Let $X$ be a piecewise Euclidean
$2$--complex with finitely many isometry
types of cells.
Then $X$ has nonpositive curvature if and only if
it satisfies the link condition.
\end{thm}

\begin{cor}
Let $X$~be a nonpositively curved, piecewise Euclidean $2$--complex
with finitely many isometry types of cells.  Then any subcomplex~$Y$
of~$X$ is also nonpositively curved
\textup{(}using the induced path metric on~$Y$\textup{)}.
\end{cor}

\begin{cor}\label{cor:diameterpi}
Let $X$~be a piecewise Euclidean $\CAT(0)$ $2$--complex,
and let $Y$~be a connected subcomplex.
Suppose that, for every vertex $v\in Y^{(0)}$, the link $\Lk(v,Y)$ has
diameter at most~$\pi$.  Then the inclusion $Y\inclusion X$ is an isometric
embedding \textup{(}using the induced path metric on~$Y$\textup{)}.
\end{cor}

%%%%%%%%%%%%%%%%%%%%%%%%%%%%%%%%%%%%%%%%%%%%%%%%%%%%%%%%%%%%%%%%%%%%%%%%%%%%
\section{Diagrams and curvature}
\label{sec:Diagrams}
%%%%%%%%%%%%%%%%%%%%%%%%%%%%%%%%%%%%%%%%%%%%%%%%%%%%%%%%%%%%%%%%%%%%%%%%%%%%

In this section we review some background on diagrams and curvature, including
the Combinatorial Gauss--Bonnet Theorem.  Diagrams are planar $2$--complexes
which
often allow one to more clearly visualize and reason about combinatorial
homotopies of paths.  The early sections of
\cite{McCammondWise02}
contain a concise and elegant development of many of the ideas in this section.

\begin{defn}[Combinatorial maps and complexes]
A map $Y\to X$ between CW~complexes is \emph{combinatorial}
if its restriction to each open cell of~$Y$ is a homeomorphism onto an open
cell of~$X$.  A CW~complex is itself \emph{combinatorial} if the attaching
map of each cell is a combinatorial map (possibly after subdividing the
given cell structures).
\end{defn}

Notice that a convex Euclidean polyhedron is a combinatorial complex.
Consequently a piecewise Euclidean complex is also combinatorial.

\begin{conv}
\label{conv:CocompactComplex}
A piecewise Euclidean $2$--complex is \emph{cocompact} if it is
cocompact as a combinatorial $2$--complex.  More precisely,
the $2$--complex~$X$ is cocompact if the group of combinatorial isometries
acts on~$X$ with a compact quotient.

Note that cocompactness of a
piecewise Euclidean $2$--complex is a stronger notion than cocompactness of the
underlying metric space.
For instance, the Euclidean plane equipped with a Penrose tiling
is not a cocompact $2$--complex even though the underlying metric space~$\E^2$ is
cocompact.
\end{conv}

\begin{defn}
A \emph{diagram}~$D$
is a finite connected
combinatorial $2$--complex equipped
with a fixed combinatorial embedding in the $2$--sphere which
misses at least one point.
A simply connected diagram is a \emph{disc diagram}.
\begin{figure}[ht!]
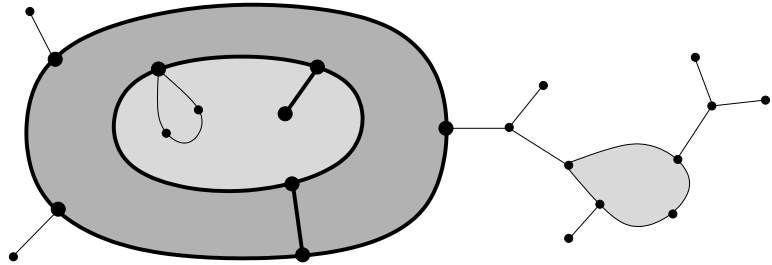

\drawdiagram
\caption[A disc diagram and a subdiagram.]{A disc diagram and a darkened
subdiagram which is not simply connected}
\label{fig:diagram}
\end{figure}

A \emph{subdiagram} $E$ of~$D$ is a subcomplex such that the embedding
$E \inclusion S^2$ factors as $E \inclusion D \inclusion S^2$
where $D \inclusion S^2$ is the given embedding of~$D$ into the sphere.
Figure~\ref{fig:diagram} shows a disc diagram and a subdiagram which is not
simply connected.
\end{defn}

Disc diagrams are often called
\emph{van~Kampen diagrams}, particularly in the context of
group presentations.
Some authors reserve the term \emph{disc diagram} for a diagram
which is homeomorphic to a disc and use \emph{singular disc diagram}
for a contractible diagram.
The author prefers the present terminology, as contractible diagrams
occur with much higher frequency in applications and seem to be
more natural objects than diagrams with the topology of a disc.

In practice, we often suppress mention of the embedding of a diagram
into the sphere.  However the boundary cycles of a diagram
(defined below) will typically be different for different embeddings of the
diagram into the sphere.  The issue is that the diagram may have cut points.

\begin{defn}[Boundary cycles]
Let $D \inclusion S^2$ be a diagram with its given combinatorial
embedding in the sphere.
Since $D$ is connected, its complement in $S^2$ has a
finite number of components $R_1,\dots,R_k$, each
homeomorphic to an open disc.
Without loss of generality, we may assume that $S^2$
has a cell structure with a single $2$--cell~$e_i$ for each region~$R_i$.
A choice of orientation for~$S^2$ determines a collection of
\emph{boundary cycles} of~$D$ so that the boundary cycle corresponding
to the region~$R_i$ is the attaching
map of the $2$--cell~$e_i$.
If $D$ contains at least one $1$--cell, then each boundary cycle is a
combinatorial map $C_i \to D$ where $C_i$ is a subdivided circle.
\end{defn}

\begin{figure}[ht!]
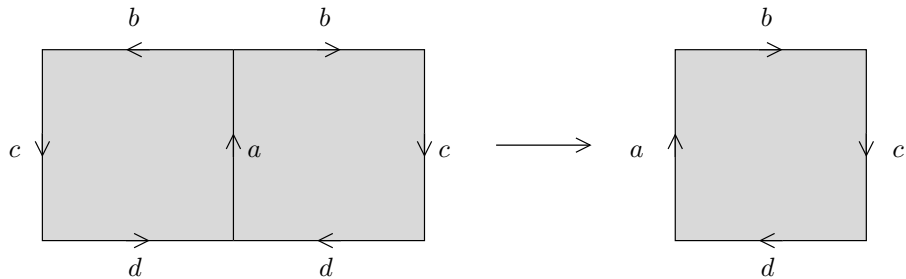

\drawcancelable
\caption[A cancelable pair of $2$--cells.]{A combinatorial map which is not
reduced because it contains a cancelable pair of $2$--cells}
\label{fig:cancelable}
\end{figure}

\begin{defn}[Reduced maps]
Let $D$~be a disc diagram and $\phi\co D\to X$ a combinatorial map.
A pair of (not necessarily distinct) $2$--cells $C_1$~and~$C_2$ in~$D$
which meet along a $1$--cell~$e$
is a \emph{cancelable pair} with respect to~$\phi$ if the boundary cycles
of $C_1$ and~$C_2$ beginning with~$e$ (in the same direction)
are not identical in~$D$, but
are sent to identical paths in~$X$ by~$\phi$.
Figure~\ref{fig:cancelable} shows a cancelable pair of $2$--cells.
The map~$\phi$ is \emph{reduced} if $D$~does not contain a cancelable pair
of $2$--cells.
We will often refer to a reduced map $D\to X$
with domain a disc diagram as a \emph{reduced disc diagram}.
\end{defn}

The following theorem concerning the existence of reduced disc diagrams
was discovered by van~Kampen \cite{vanKampen33} and independently by
Lyndon \cite{Lyndon66}.
For a complete proof, see for instance \cite[Lemma~2.17]{McCammondWise02}.

\begin{thm}\label{thm:ReducedDiagram}
If $X$~is a combinatorial $2$--complex and $P\to X$ is a combinatorial
closed path which is nullhomotopic in~$X$, then there exists a reduced disc
diagram $D\to X$ so that $P\to X$ is the boundary
path of~$D$.
\end{thm}

We will often refer to the reduced disc diagram $D\to X$
obtained in the preceding theorem as a \emph{reduced disc diagram for $P$}.

\begin{figure}[ht!]
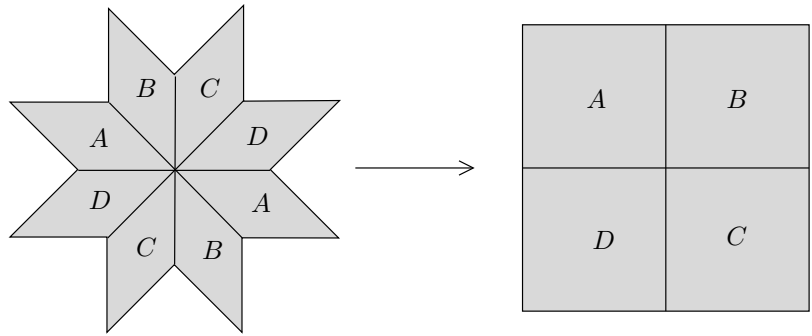

\drawBranchedCover
\caption[The induced $\CAT(0)$ structure of a reduced disc diagram.]%
{A reduced map from a disc diagram to a nonpositively curved squared
$2$--complex. The diagram on the left inherits a $\CAT(0)$ metric in which the
eight $2$--cells are Euclidean squares.}
\label{fig:BranchedCover}
\end{figure}

In the following theorem we show that nonpositive curvature pulls back under
reduced disc diagrams.  We emphasize that the metric constructed on the disc
diagram is unrelated to any metric inherited from the given embedding of the
diagram in the plane~$\R^2$ (see Figure~\ref{fig:BranchedCover}).

\begin{prop}\label{prop:CAT0diagram}
Let $\phi\co D \to X$ be a reduced disc diagram, such that $X$ is a
nonpositively curved piecewise Euclidean $2$--complex.
Then $D$ is a $\CAT(0)$ space.
\end{prop}

\begin{proof}
Since $D$~is connected and simply connected, we just need to show that
$D$~is
nonpositively curved.  Since $\phi$~is combinatorial, $D$~has a natural
piecewise
Euclidean structure obtained by pulling back the Euclidean metrics
on each cell of the image.

We will now show that for each vertex $v\in D^{(0)}$, the induced map
\[
   \phi_*\co \Lk(v,D) \to \Lk \bigl( \phi(v),X \bigr)
\]
is a local isometry.
If not, then $\phi_*$ \emph{folds} a pair of edges together.
In other words, there is a pair of distinct oriented
edges $e_1$ and~$e_2$ in $\Lk(v,D)$ with the same initial vertex which map
to the same oriented edge under~$\phi_*$.
But such a pair of edges in
$\Lk(v,D)$ corresponds to two distinct corners of
of a cancelable pair of $2$--cells
of~$D$, contradicting the fact that $\phi$~is reduced.  Hence for every vertex
$v\in D^{(0)}$ the map $\phi_*$ is a local isometry.

Thus if $C\to\Lk(v,D)$ is a locally geodesic loop,
then the composition
$C \to \Lk(v,D) \to \Lk \bigl( \phi(v),X \bigr)$
is a locally geodesic loop of the same length.
Hence $C$~has length at least $2\pi$, and we see that $D$~has
nonpositive curvature.
\end{proof}

\begin{figure}[ht!]
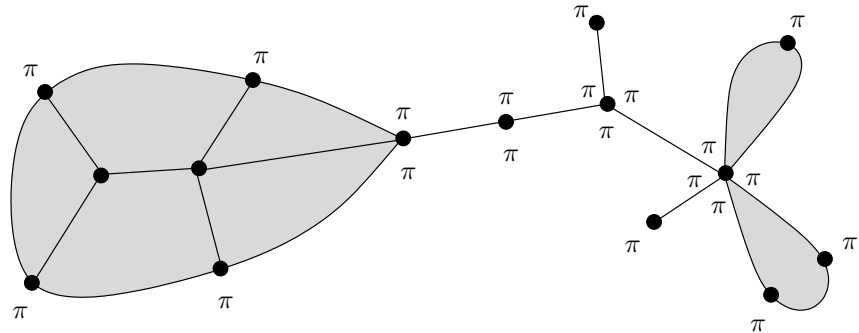

\drawcurvature
\caption[The curvature at a vertex.]%
{Formally place an angle of size~$\pi$ at each corner of $S^2 \setminus D$.
The curvature at a vertex is $2\pi$ minus the angle sum of all corners at~$v$
including these ``exterior'' corners.}
\label{fig:Curvature}
\end{figure}

\begin{defn}[Curvature]
Let $D$ be a piecewise Euclidean diagram.
The \emph{curvature at a vertex}~$v$ of~$D$, denoted $\kappa(v)$,
is defined by the formula
\[
  \kappa(v) = 2\pi - \pi\,\chi \bigl( \Lk(v,D) \bigr) \, - \!
      \sum_{e\in \operatorname{Edges} (\Lk(v,D) )} \!\!\norm{e}\text{,}
\]
where $\chi$ denotes Euler characteristic
and $\norm{e}$ denotes the length of the edge~$e$.
We occasionally use the notation $\kappa_D(v)$ to emphasize
the specific diagram~$D$ in which the curvature is measured.
An alternate way to interpret curvature is to formally place an angle
of size~$\pi$ at each corner of the complement of~$D$ in~$S^2$
as illustrated in Figure~\ref{fig:Curvature}.
Then the curvature
at~$v$ is equal to $2\pi$ minus the angle sum of all the corners at~$v$
including these ``exterior'' corners.
\end{defn}

We now have two distinct notions of nonpositive curvature
in a piecewise Euclidean diagram.  The following lemma gives a simple
correspondence between these two ideas.

\begin{lem}\label{lem:NPCisNPC}
Let $D$~be a piecewise Euclidean diagram.  Then $D$~has nonpositive
curvature
if and only if the curvature at each interior vertex of~$D$
is nonpositive.
\end{lem}

\begin{proof}
Notice that the link of an interior vertex is always a circle.
The proof now follows immediately from the definition of curvature
at a vertex
and the Link Condition for~$D$.
\end{proof}

The Combinatorial Gauss--Bonnet Theorem is based on Lyndon's Curvature Formula
\cite{LS77},
and was proved for piecewise Euclidean
disc diagrams independently by Gersten and Pride \cite{Gersten87,Pride88}.
It has since been generalized to
arbitrary combinatorial $2$--complexes by McCammond--Wise
\cite{McCammondWise02}.

\begin{thm}[Combinatorial Gauss--Bonnet]\label{thm:CGB}
Let $D$~be a piecewise Euclidean diagram.  Then
\[
  \sum_{v\in D^{(0)}} \kappa(v) = 2\pi\,\chi(D).
\]
\end{thm}

%%%%%%%%%%%%%%%%%%%%%%%%%%%%%%%%%%%%%%%%%%%%%%%%%%%%%%%%%%%%%%%%%%%%%%%%%%%%
\section{Isolated Flats and the Flat Triplane Theorem}
\label{sec:FlatTriplaneTheorem}
%%%%%%%%%%%%%%%%%%%%%%%%%%%%%%%%%%%%%%%%%%%%%%%%%%%%%%%%%%%%%%%%%%%%%%%%%%%%

\begin{defn}\label{def:IsolatedFlats}
A $\CAT(0)$ $2$--complex~$X$ has the \emph{Isolated Flats Property}
if there is a function $\psi \co \R_+ \to \R_+$
such that for every pair of distinct flat planes $F_1 \ne F_2$ in~$X$
and for every $k\ge 0$, the intersection $\nbd{F_1}{k} \cap \nbd{F_2}{k}$
of $k$--neighborhoods of $F_1$ and~$F_2$ has diameter at most $\psi(k)$.
\end{defn}

The definition of the Isolated Flats Property given above is
catered to the two dimensional setting, and as such, is
slightly
simpler than the definition used in \cite{HruskaGeometric}
and \cite{HruskaRelHyp}.  This simplification is due to the fact that
several phenomema
present in higher dimensions are absent in the
$2$--dimensional setting.  For instance, in general a space with
the Isolated Flats Property could contain subspaces of the form
$\E^k \times K$ where $K$ is a nontrivial compact set.

We note the following immediate consequence of the Isolated Flats Property,
which will be useful in the sequel.

\begin{prop}\label{prop:IFPEquivalent}
Let $X$ be a $\CAT(0)$ space with the Isolated Flats Property.
Then for every $k \ge 0$, each flat disc in~$X$ of radius
at least $\psi(k)$ lies in a $k$--neighborhood of
at most one flat.
\end{prop}

The Flat Plane Theorem states that a proper, cocompact $\CAT(0)$ space is
$\delta$--hyperbolic if and only if it does not contain an isometrically
embedded flat plane.%
\footnote{The Flat Plane Theorem was proved for Riemannian manifolds
  by Eberlein \cite{Eberlein73}.
  Gromov stated the theorem for general $\CAT(0)$ spaces in
  \cite[\S 4.1]{Gromov87}.
  Heber and Bridson independently provided proofs in this general setting
  \cite{Heber87,Bridson95}.
}
Wise has proved an analogous result for two-complexes with the
Isolated Flats Property.
His Flat Triplane Theorem shows that in the two-dimensional case,
the Isolated Flats Property is equivalent to an absence of triplanes.
Since this useful result has not appeared in the literature,
we provide Wise's proof in this section for the sake of completeness.
For related results regarding the presence of triplanes in
nonpositively curved $2$--complexes, see \cite{Wise96,WiseFigure8}.

\begin{defn}\label{def:triplane}
A \emph{triplane} is the space formed from three closed half-planes
by gluing their boundary lines (by isometries) to a common line.
Notice that a triplane is a piecewise Euclidean $\CAT(0)$ $2$--complex.

Fix a basepoint $x$ on the singular line of a triplane~$T$.
Let $T_r$ denote the closed ball in~$T$ of radius~$r$ around~$X$.
Notice that $T_r$ looks like three half discs of radius~$r$ with their
straight boundary sides glued together.
\end{defn}

\begin{thm}[Flat Triplane Theorem]\label{thm:FlatTriplane}
Let $X$ be a proper, cocompact piecewise Euclidean $2$--complex.
The following are equivalent:
\begin{enumerate}
   \item \label{item:TriplaneIFP}
   $X$ has the Isolated Flats Property.
   \item \label{item:TriplaneUBI}
   There is a universal bound~$L$ on the diameter of the intersection
   of any two distinct flat planes in~$X$.
   \item \label{item:TriplaneNTP}
   $X$ does not contain an isometrically embedded triplane.
\end{enumerate}
\end{thm}

The proof of Theorem~\ref{thm:FlatTriplane} uses the following
variant of the Arzel\`{a}--Ascoli theorem, which is proved in
\cite[Lemma~II.9.34]{BH99}.
The statement given here is slightly stronger than the
one given by Bridson and Haefliger. This additional strength follows immediately
from their proof.

\begin{lem}\label{lem:embedball}
Let $Y$ be a separable metric space with basepoint $y_0$, and let
$X$ be a proper, cocompact metric space.  If for each $n \in \N$ there
is an isometric embedding $\phi_n\co \ball{y_0}{n} \inclusion X$
then there exists an isometric embedding $\phi\co Y \inclusion X$.

Furthermore, if we assume the existence of a compact set $K$ such that
for each~$n$ the point $\phi_n(y_0)$ lies in $K$, then we may take $\phi$ to
be a pointwise limit of a subsequence $\{\phi_{n_i}\}$ of the original
sequence of embeddings.
\end{lem}

The proof of Theorem~\ref{thm:FlatTriplane} also uses the Flat Strip
Theorem \cite{BH99},
which states the following.

\begin{thm}[Flat Strip Theorem]\label{thm:FlatStrip}
Let $X$ be a $\CAT(0)$ space, and let $\gamma\co\R\to X$ and
$\gamma'\co\R\to X$
be geodesic lines in $X$.  If there is a constant~$K$ such that
$d \bigl( \gamma(t), \gamma'(t) \bigr) \le K$ for all $t\in\R$,
then the convex hull of\/
$\gamma(\R) \cup \gamma'(\R)$ is isometric to a flat strip
$\R \times [0,L] \subset \E^2$.
\end{thm}

\begin{proof}[Proof of Theorem~\ref{thm:FlatTriplane}]
The implications
$\text{(\ref{item:TriplaneIFP})} \implies \text{(\ref{item:TriplaneUBI})}
  \implies \text{(\ref{item:TriplaneNTP})}$
are immediate.
We now show that (\ref{item:TriplaneNTP}) implies~(\ref{item:TriplaneUBI}).
Suppose for each $k>0$ there is a pair of distinct flat planes $E_k\ne F_k$
whose intersection~$I_k$ has diameter at least~$k$.
Since $E_k$~and~$F_k$ are convex, their intersection~$I_k$ is isometric to
a convex polygonal region (possibly unbounded) in the Euclidean plane.
Since $X$ is a cocompact $2$--complex,
it has only finitely many isometry types of cells (see
Convention~\ref{conv:CocompactComplex}).
So there are finitely many possible turning angles on the boundary of~$I_k$.
Thus there is a positive number~$\theta$ (independent of~$k$)
so that every positive turning angle on~$\boundary I_k$ is
at least~$\theta$.
It follows that there are at most $2\pi/\theta$ turns on the boundary
of~$I_k$.  Thus for each $m>0$, there is an $n>0$ so that $\boundary I_n$
contains a straight line segment~$\gamma_m$ of length~$m$.
This segment lies in two distinct flat discs $D_m$ and~$D'_m$ of
diameter~$m$
which intersect in a half disc.
Notice that $D_m \cup D'_m$ with its induced path metric
is isometric to the space $T_m$ (in the notation of
Definition~\ref{def:triplane}).
But for every point $v\in T_m$, the link $\Lk(v,T_m)$
has diameter~$\pi$, so $T_m\inclusion X$ is an isometric embedding by
Corollary~\ref{cor:diameterpi}.  Since $X$ contains an isometrically
embedded
copy of~$T_m$ for each $m>0$, it follows from Lemma~\ref{lem:embedball}
that $X$ contains an isometrically embedded triplane.

We have shown that (\ref{item:TriplaneUBI}) and~(\ref{item:TriplaneNTP})
are equivalent.
Now we show that (\ref{item:TriplaneUBI}) and~(\ref{item:TriplaneNTP})
together
imply~(\ref{item:TriplaneIFP}).
Suppose $X$ does not have the Isolated Flats Property.
Then there exists a constant
$R > 0$ such that for every $k \in \N$ there are distinct flats
$E_k$ and~$F_k$ in~$X$ with geodesic segments
$\gamma_k \subset E_k$ and $\gamma'_k \subset F_k$,
each with length greater than~$k$ such that
the Hausdorff distance between their images is at most~$R$.

Consider the sequence of segments~$\gamma_k$.
By composing with a suitable isometry, we may assume that the midpoint~$x_k$
of each segment~$\gamma_k$ lies inside a given compact set~$K$.
Then Lemma~\ref{lem:embedball} implies that there is a subsequence
$\gamma_{k_i}$ that converges pointwise to a geodesic~$\gamma$.
Applying Lemma~\ref{lem:embedball} to the sequence of embedded flat planes
$E_{k_i}$ gives a further subsequence of integers $\{m_i\}$
such that the flat planes $E_{m_i}$ converge to a
flat plane~$E$ containing the geodesic~$\gamma$.
Letting $x'_k$ denote the midpoint of~$\gamma'_k$, we see that
$d(x_k,x'_k)\le 3R$,
since the endpoints of $\gamma_k$ and $\gamma'_k$
are within $3R$ of each other.
So $x'_k$ lies in the closure $K'$ of a $3R$--neighborhood of~$K$, which is
compact, since $X$ is proper.

Continuing as above, we eventually obtain a sequence of integers
$\{n_i\}$ such that $E_{n_i}$ and~$F_{n_i}$ converge to isometrically
embedded
flat planes $E$ and~$F$ containing geodesic lines $\gamma$ and~$\gamma'$.
Furthermore, we may assume that the Hausdorff distance between the images of
$\gamma$ and~$\gamma'$ is at most~$R$.
If $E$ and~$F$ are actually the same flat plane, then the intersections
$E_{n_i}\cap F_{n_i}$ must have arbitrarily large diameter,
contradicting~(\ref{item:TriplaneUBI}).  So we may assume that $E\ne F$.

By Theorem~\ref{thm:FlatStrip}, we conclude that the convex hull~$S$ of
$\Image(\gamma) \cup \Image(\gamma')$ is isometric to a flat strip
$\R \times [0,L] \subset \E^2$.
Since $E$,~$F$, and~$S$ are convex, any intersection of them is also convex.
It now follows that $E\cap S$ and $F\cap S$ are substrips
$S_E = \R \times [0,s]$ and $S_F = \R\times [t,L]$ respectively.
We now have two cases, depending on whether these substrips intersect.

{\bf Case 1}\qua Suppose $s < t$.
Then $E$~and~$F$ are disjoint planes connected by a
flat strip.  Letting $U = E\cup F\cup S$, we see that
for every point $v\in U$ the link $\Lk(v,U)$
has diameter~$\pi$.  So the inclusion
$U\inclusion X$ is an isometric embedding by Corollary~\ref{cor:diameterpi}.
Since $U$~contains an isometrically embedded triplane, we are done.

{\bf Case 2}\qua Suppose $s \ge t$.  Then $E\cap F$ is a closed convex set
of the plane
containing a line.  It follows that $E\cap F$ is either a half-plane, or a
(possibly degenerate) strip $\R \times [0, s - t]$.
Let $U = E\cup F$.
As in the previous case, $U$~is isometrically embedded in~$X$ since
every link $\Lk(v,U)$ has diameter~$\pi$.
Since $U$~contains a triplane, we are done.
\end{proof}

%%%%%%%%%%%%%%%%%%%%%%%%%%%%%%%%%%%%%%%%%%%%%%%%%%%%%%%%%%%%%%%%%%%%%%%%%%
\section{Hyperbolicity relative to flats}
\label{sec:HyperbolicityRelFlats}
%%%%%%%%%%%%%%%%%%%%%%%%%%%%%%%%%%%%%%%%%%%%%%%%%%%%%%%%%%%%%%%%%%%%%%%%%%

In this section we give precise definitions of the Relative Fellow Traveller
Property and the Relatively Thin Triangle Property.  We show that in the
$2$--dimensional
setting, each of these properties implies the Isolated Flats Property,
establishing the implications (2)~$\implies$~(1) and (3)~$\implies$~(1)
of Theorem~\ref{thm:2dEquivalent}.

\begin{figure}[ht!]
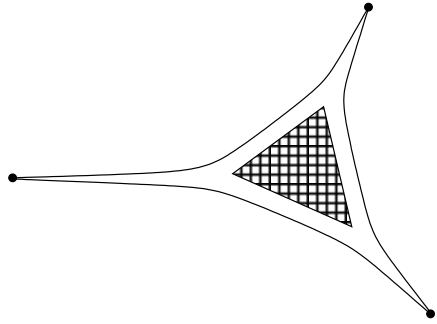

\drawthintriangle
\caption{A triangle which is $\delta$--thin relative to a flat}
\label{fig:ThinTriangleRelFlat}
\end{figure}

\begin{defn}[Relatively Thin Triangle Property]
A geodesic triangle in a space is \emph{$\delta$--thin relative to the flat~$F$}
if each side of the triangle lies in a $\delta$--neighborhood of the other two
sides and the flat~$F$, as illustrated in Figure~\ref{fig:ThinTriangleRelFlat}.
A space~$X$ has the \emph{Relatively Thin Triangle Property}
if there is a constant $\delta$ so that
each triangle in~$X$ is either $\delta$--thin in the usual sense or
$\delta$--thin relative to some flat.
\end{defn}

It is not hard to see that any $2$--complex with the Relatively Thin Triangle
Property must also have the Isolated Flats Property.

\begin{thm}[Relatively Thin Triangle Property
             $\implies$ Isolated Flats Property]
\label{thm:ThinTriangles=>IsolatedFlats}
Let $X$ be a proper, cocompact, piecewise Euclidean $\CAT(0)$ $2$--complex
satisfying the Relatively Thin Triangle Property.
Then $X$ also has the Isolated Flats Property.
\end{thm}

\begin{proof}
Assume by way of contradiction that $X$ does not have the Isolated Flats
Property.
Then by Theorem~\ref{thm:FlatTriplane}, $X$ contains an isometrically
embedded triplane~$T$.
Let us parametrize $T$ as
\[
   T = \bigl\{\, (x,y,i) \bigm|
       \text{$x,y \in \R$; $y \ge 0$; and $i \in \{1,2,3\}$} \bigr\}
         \Big/  (x,0,j) \sim (x,0,k).
\]
Now consider the triangle $\Delta_n$ with vertices
\[
  a = (0,n,1) \qquad  b = (-2n,n,2) \qquad c = (2n,n,3)
\]
illustrated in Figure~\ref{fig:BadTriangle}.
Notice that for any fixed~$\delta$, each side of $\Delta_n$ only intersects
the $\delta$--neighborhood of the union of the other two sides
near the corners of the triangle.  This intersection consists
of two
segments whose lengths are bounded by a constant which does not depend on
the value of~$n$.  So for large values of~$n$, the triangle~$\Delta_n$ is
not $\delta$--thin.
\begin{figure}[ht!]
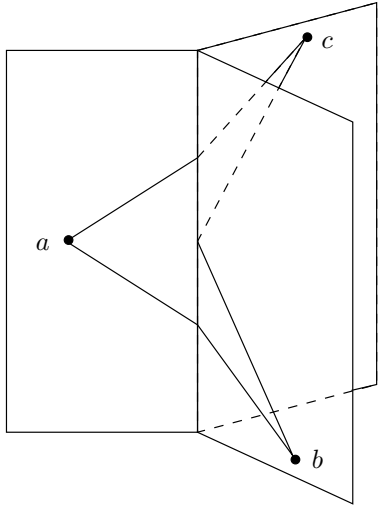

\drawBadTriangle
\caption[A fat triangle which does not lie close to a single flat.]%
{A fat triangle which does not lie close to a single flat.
The two paths from $b$ to~$c$ are quasigeodesics which are only close
together near their endpoints.}
\label{fig:BadTriangle}
\end{figure}

But if~$n$ is sufficiently large, each side of~$\Delta_n$ lies in a
$\delta$--neighborhood of at most one flat in~$X$.  Since
the three sides of~$\Delta_n$ lie in three distinct flats,
it is clear that $\Delta_n$ is not $\delta$--thin relative to any
single flat in~$X$.
So $X$ does not have the Relatively Thin Triangle Property.
\end{proof}

\begin{defn}
A \emph{$(\lambda,\epsilon)$--quasigeodesic} in a metric space~$X$
is a function $\alpha \co [a,b] \to X$ for some real interval $[a,b]$
satisfying
\[
  \frac{1}{\lambda} \, d(s,t) - \epsilon
  \le d \bigl( \alpha(s), \alpha(t) \bigr)
  \le \lambda \, d(s,t) + \epsilon
\]
for all $s,t \in [a,b]$.
A map $\alpha\co [a,b] \to X$ is a \emph{quasigeodesic} if
there exist constants
$\lambda$ and~$\epsilon$ such that $\alpha$ is a
$(\lambda,\epsilon)$--quasigeodesic.
\end{defn}

In hyperbolic geometry, quasigeodesics with common endpoints
satisfy an asynchronous fellow traveller property.  This fact was
established for the hyperbolic plane by Morse \cite{Morse24}
and for $\Hyp^n$ by
Efromovich--Tihomirova
\cite{EfromovichTihomirova63}.
Gromov generalized the fellow traveller property to the
following result about $\delta$--hyperbolic spaces.

\begin{thm}[\cite{Gromov87}, Proposition~7.2.A]
\label{thm:GromovFTP}
Let $\alpha$ and~$\beta$ be a pair of $(\lambda,\epsilon)$--quasigeodesics
with common endpoints in a $\delta$--hyperbolic space~$X$.
Then the Hausdorff distance between $\Image(\alpha)$ and
$\Image(\beta)$ is at most~$L$, where $L$ depends only
on the constants $\delta$, $\lambda$, and~$\epsilon$.
\end{thm}

\begin{figure}[ht!]
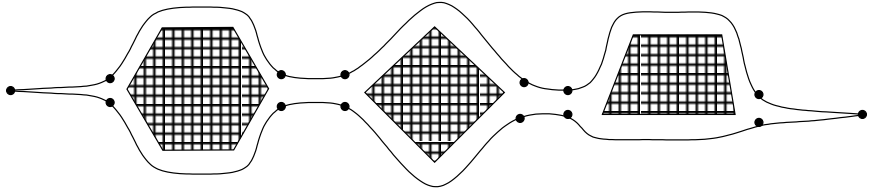

\drawrftp
\caption{A pair of paths which fellow travel relative to flats}
\label{fig:RelativeFTP}
\end{figure}

The concept of a pair of paths which \emph{fellow travel relative to flats}
generalizes the asynchronous fellow travelling described above.
Roughly speaking, the idea is that the two paths alternate between tracking
close together and travelling near a common flat
as illustrated in Figure~\ref{fig:RelativeFTP}.

\begin{defn}[Fellow travelling relative to flats]
\label{def:RelativeFellowTravelling}
A pair of paths
\[
   \alpha\co [0,a]\to X \quad \text{and}
   \quad \alpha'\co [0,a']\to X
\]
in
a $\CAT(0)$ space \emph{$L$--fellow travel relative to a sequence of flats}
$(F_1,\dots,F_n)$
if there are partitions
\[
  0 = t_0 \le s_0 \le t_1 \le s_1 \le \dots \le t_n \le s_n = a
\]
and
\[
  0 = t'_0\le s'_0\le t'_1\le s'_1\le \dots \le t'_n\le s'_n= a'
\]
so that for $0\le i \le n$ the Hausdorff distance between the
sets $\alpha \bigl( [t_i,s_i] \bigr)$ and
$\alpha' \bigl( [t'_i,s'_i] \bigr)$ is at most~$L$,
while for $1\le i \le n$ the sets
$\alpha \bigl( [s_{i-1},t_i] \bigr)$
and $\alpha' \bigl( [s'_{i-1},t'_i] \bigr)$ lie in an $L$--neighborhood of
the flat $F_i$.

We will frequently say that paths \emph{$L$--fellow travel relative
to flats} if they $L$--fellow travel relative to some sequence of flats.
\end{defn}

\begin{defn}\label{def:RelFTP}
A $\CAT(0)$ space $X$ satisfies the
\emph{Relative Fellow Traveller Property}
if for each choice of constants $\lambda$ and~$\epsilon$
there is a constant $L = L(\lambda,\epsilon,X)$ such that
$(\lambda,\epsilon)$--quasigeodesics in~$X$
with common endpoints $L$--fellow travel relative to flats.
\end{defn}

In the $2$--dimensional setting, it is easy to show that the Relative Fellow
Traveller Property implies the Isolated Flats Property using
the same counterexample as the proof of
Theorem~\ref{thm:ThinTriangles=>IsolatedFlats}.
The converse requires substantially more work and occupies
most of the present article.

\begin{thm}[Relative Fellow Traveller Property $\implies$ Isolated Flats
Prop\-erty]
\label{thm:RFTP=>IFP}
Let $X$ be a proper, cocompact piecewise Euclidean $\CAT(0)$ $2$--complex.
If $X$ has the Relative Fellow Traveller Property, then $X$ also has the
Isolated Flats Property.
\end{thm}

\begin{proof}
Suppose $X$ contains a triplane~$T$.  As in the proof of
Theorem~\ref{thm:ThinTriangles=>IsolatedFlats}, consider the triangle
$\Delta_n$ with vertices
\[
  a = (0,n,1) \qquad  b = (-2n,n,2) \qquad c = (2n,n,3)
\]
illustrated in Figure~\ref{fig:BadTriangle}.
Notice that the union $[b,a] \cup [a,c]$ is a quasigeodesic when parametrized
by
arclength.  Furthermore, the associated constants of this quasigeodesic are
independent of~$n$.

Also note that for each~$L$, the geodesic $[b,c]$ intersects the
$L$--neighborhood of the quasigeodesic only near their common endpoints.
But
$[b,a] \cup [a,c]$ does not lie in a $C$--neighborhood of any flat
for any constant~$C$ which does not depend on~$n$.
So $X$ does not have the Relative Fellow Traveller Property.
\end{proof}

%%%%%%%%%%%%%%%%%%%%%%%%%%%%%%%%%%%%%%%%%%%%%%%%%%%%%%%%%%%%%%%%%%%%%%%%%%
\section{Ruffled boundaries and thin triangles}
\label{sec:RuffledBoundaries}
%%%%%%%%%%%%%%%%%%%%%%%%%%%%%%%%%%%%%%%%%%%%%%%%%%%%%%%%%%%%%%%%%%%%%%%%%%

In the previous section, we established the implications
(2)~$\implies$~(1) and (3)~$\implies$~(1) of
Theorem~\ref{thm:2dEquivalent}.  The remainder of this article is devoted
to proving (1)~$\implies$~(2) and (1)~$\implies$~(3), which we will establish
in Theorems \ref{thm:IFP=>ThinTriangles} and~\ref{thm:RelativeFTP}
respectively.

Our objective for the next two sections is to establish some diagrammatic
tools that will be useful in the proofs of these two theorems.
In this section we introduce the notion of a diagram which is
\emph{ruffled} along a certain part of its boundary.  This notion
generalizes the fact that
in the $\delta$--hyperbolic setting, every point of a piecewise Euclidean
diagram
is either close to a point of negative curvature or lies near the boundary of
the diagram.
As a consequence, every $\delta$--hyperbolic diagram
is ``ruffled'' throughout its interior.

When studying the Isolated Flats Property, one frequently encounters
diagrams with the property
that negative curvature is distributed evenly along
a certain portion of the boundary of the diagram
except in places where the
diagram is very thin.  Such a diagram is
\emph{ruffled} along that boundary path.

Diagrams which are ruffled along part of their boundary have many features in
common with $\delta$--hyperbolic diagrams.
For instance, we will see in Proposition~\ref{prop:ThinTriangles} that a
triangular diagram which is ruffled along one side must be $\delta$--thin.
In this section,
we also prove results which describe various ways
that ruffles can be inherited by a subdiagram.

\begin{defn}[Ruffled]\label{def:Ruffled}
Let $D$ be a nonpositively curved piecewise Euclidean diagram with boundary
cycles $C_1, \dots, C_n$.  Suppose $C_1$ is a concatenation
of two paths $\alpha$ and~$\beta$.  Then the pair
$(D, \alpha)$ is
\emph{$(R,\theta)$--ruffled} for positive constants $R$
and~$\theta$ provided that for each point $\alpha(t)$, the open ball
$B = \bigball{\alpha(t)}{R}$ in~$D$
satisfies one of the following two properties:
\begin{description}
   \item[(R-1)] $B$~contains a vertex~$v$ with curvature
      $\kappa_D(v) \le -\theta$, or
   \item[(R-2)] $B$~intersects the image of at least one of the curves
      $C_2,\dots,C_n,$ or $\beta$.
\end{description}

The pair $(D,C_1)$ is \emph{$(R,\theta)$--ruffled} provided that
for each point $C_1(t)$, the open ball
$B= \bigball{C_1(t)}{R}$ in~$D$
satisfied one of the following two properties:
\begin{description}
   \item[(R$\mathbf{'}$-1)] $B$ contains a vertex~$v$ with curvature
      $\kappa_D(v) \le -\theta$, or
   \item[(R$\mathbf{'}$-2)] $B$ intersects the image of at least one
      of the curves $C_2,\dots,C_n$.
\end{description}
\end{defn}

\begin{rem}
Note that if $(D,\alpha)$ or $(D,C_1)$ is $(R_0,\theta_0)$--ruffled, then
it is also $(R,\theta)$--ruffled for any $R\ge R_0$
and any positive $\theta \le \theta_0$.
\end{rem}

The following lemma follows immediately from the definition of ruffled.

\begin{figure}[ht!]
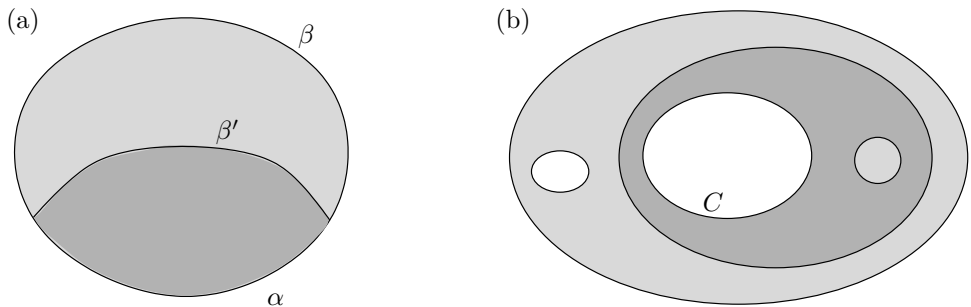

\drawSimpleRuffles
\caption[A diagram and a subdiagram which share a boundary arc.]% 
{(a) A diagram and a darkened subdiagram which share a boundary
arc~$\alpha$\qua (b) A diagram and a darkened subdiagram which share
an entire boundary cycle}
\label{fig:SimpleRuffles}
\end{figure}

\begin{lem}[Subdiagrams]\label{lem:SubdiagramRuffled}
Consider a subdiagram~$D'$ of a nonpositively curved piecewise Euclidean
diagram~$D$.
Suppose there are boundary cycles $C$~and\/~$C'$ of\/ $D$~and\/~$D'$
respectively
so that $C$ is a concatenation $\alpha\beta$
and $C'$ is a concatenation $\alpha\beta'$,
as illustrated in Figure~\ref{fig:SimpleRuffles}\textup{(}a\textup{)}.
If the pair $(D,\alpha)$ is $(R,\theta)$--ruffled
then so is the pair $(D',\alpha)$.

Similarly suppose $C$ is a boundary cycle of
both $D$ and~$D'$ as in Figure~\ref{fig:SimpleRuffles}\textup{(}b\textup{)}.
If $(D,C)$ is $(R,\theta)$--ruffled,
then $(D',C)$ is also $(R,\theta)$--ruffled.
\end{lem}

The following lemma says that under many circumstances
a concatenation of ruffled boundary arcs is also a ruffled
boundary arc.  In general the constants associated to the ruffles
depend on the number of segments being concatenated.
The issue is that one way for a boundary arc to be ruffled
is for the arc to be very short.
If each arc in the concatenation is known to be sufficiently long,
then ruffling constants are obtained which do not depend on the
number of concatenated segments.

\begin{lem}[Concatenations]\label{lem:ConcatenatedRuffles}
Let $C$ be a boundary cycle of a nonpositively curved piecewise Euclidean
diagram~$D$.  Suppose $C$ is a concatenation
$\alpha_1 \dotsm \alpha_n \beta$
such that $(D,\alpha_i)$ is $(R,\theta)$--ruffled for each~$i$
and $\alpha=\alpha_1 \dotsm \alpha_n$ is a local geodesic
in~$D$.

Suppose further that any \textup{(}global\textup{)} geodesic in~$D$ connecting
two points of\/
$\Image(\alpha)$ lies inside $\Image(C)$.
Then
\begin{enumerate}
   \item\label{item:CRa} $(D,\alpha)$ is $(nR,\theta)$--ruffled.
   \item\label{item:CRb} If each $\alpha_i$ has length at least~$2R$, then
      $(D,\alpha)$ is $(2R,\theta)$--ruffled.
   \item\label{item:CRc} If each $\alpha_i$ has length at least~$2R$
      and $\beta$ has image a single point, then $(D,C)$ is
      $(2R,\theta)$--ruffled.
\end{enumerate}
\end{lem}

\begin{figure}[ht!]
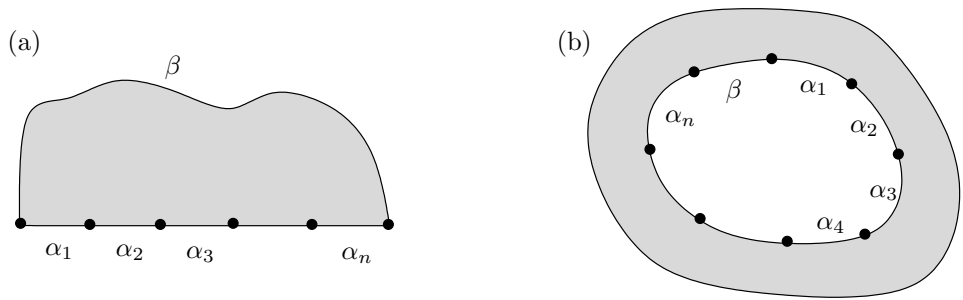

\drawconcatenations
\caption[An illustration of the hypothesis
    of Lemma~\ref{lem:ConcatenatedRuffles}.]%
   {(a)~A disc diagram such that the concatenation
   $\alpha_1 \dotsm \alpha_n$ is a geodesic boundary arc\qua
   (b)~An annular diagram such that the concatenation
$\alpha_1 \dotsm \alpha_n \beta$ is a locally geodesic boundary component}
\label{fig:concatenations}
\end{figure}

The requirement that any geodesic in~$D$ connecting two points of
$\Image(\alpha)$ lies inside $\Image(C)$ is satisfied if
$D$ is a disc diagram and $\alpha$ is a geodesic boundary component,
as shown in Figure~\ref{fig:concatenations}(a),
or alternately if
$D$ is an annular diagram and $C$ is a locally geodesic boundary
component as in Figure~\ref{fig:concatenations}(b).

\begin{proof}
If every $\alpha_i$ has length less than~$2R$, then every point of
$\Image(\alpha)$ is within a distance $nR$ of an endpoint of~$\alpha$,
and hence $(D,\alpha)$ is $(nR,\theta)$--ruffled.
On the other hand, suppose some $\alpha_i$ has length at least~$2R$,
and $x_i$ is a point
on $\Image(\alpha_i)$ at least a distance~$R$ from both ends of~$\alpha_i$.
Then the ball $B=\ball{x_i}{R}$ intersects $\Image(\alpha)$
only in the interior of $\Image(\alpha_i)$.
If $B$ does not contain a vertex with curvature $\le -\theta$,
then $B$ intersects $\Image(\beta)$ or some other boundary component of~$D$.
Every point of $\Image(\alpha)$ is within a distance $(n-1)R$ of such
a point~$x_i$.  Therefore (\ref{item:CRa}) holds.

If each $\alpha_i$ has length at least~$2R$, then each contains
a point~$x_i$ at least a distance~$R$ from both ends of $\Image(\alpha)$.
Therefore, every point of~$\Image(\alpha)$ is within a distance~$R$ of such a
point~$x_i$.  Claim~(\ref{item:CRb}) now follows immediately.
Furthermore, if $\beta$ has image a single point, then that point
is also within
a distance~$R$ of some~$x_i$, giving~(\ref{item:CRc}).
\end{proof}

\begin{defn}[Triangular]\label{def:TriangularDiagram}
A piecewise Euclidean $\CAT(0)$ disc diagram~$\Delta$ is \emph{triangular}
if its boundary cycle is a concatenation of three geodesics
$\alpha$,~$\beta$, and~$\gamma$.  These geodesics are the \emph{sides}
of the triangle, and their endpoints are its \emph{corners}.
\end{defn}

\begin{lem}\label{lem:OneSideThin}
Let $\Delta$ be a geodesic triangle in a $\CAT(0)$ space $X$.
Suppose one of the sides of~$\Delta$ lies in
a $\delta$--neighborhood
of the union of the other two sides.
Then $\Delta$~is $2\delta$--thin.
\end{lem}

\begin{proof}
Let $a$,~$b$, and~$c$ be points lying on sides
$\alpha$,~$\beta$, and~$\gamma$
respectively so that $d(a,b)$ and $d(a,c)$ are each less than~$\delta$.
Then $b$ lies in a $2\delta$--neighborhood of both $\alpha$ and~$\gamma$.
Similarly, $c$ lies in a $2\delta$--neighborhood of both
$\alpha$ and~$\beta$.
The result now follows from the convexity of the $\CAT(0)$ metric in~$X$.
\end{proof}

The following proposition says that a triangular diagram which is ruffled along 
one side
is $\delta$--thin for some~$\delta$.  The proof uses the Combinatorial
Gauss--Bonnet Theorem to bound the total amount of negative curvature inside the
triangle.

\begin{prop}[Thin triangles]\label{prop:ThinTriangles}
Let $\Delta$~be a piecewise Euclidean $\CAT(0)$ triangular diagram with
sides $\alpha$,~$\beta$, and~$\gamma$.
Suppose $R$ and~$\theta$ are positive constants so that
$(\Delta,\alpha)$ is $(R,\theta)$--ruffled.
Then there is a constant $\delta=\delta(R,\theta)$
so that $\Delta$~is $\delta$--thin.
\end{prop}

\begin{proof}
By Theorem~\ref{thm:CGB}, the sum of
the curvatures
at all the vertices of~$\Delta$ is exactly~$2\pi$.
But by Lemma~\ref{lem:NPCisNPC}
the only positive curvature in~$\Delta$ occurs at its three corners since
$\Delta$~is nonpositively curved.
Furthermore, the curvature at each corner is at most~$\pi$.
So the sum of all positive curvatures in~$\Delta$ is at most~$3\pi$.
It follows that the sum of all negative curvatures in~$\Delta$
has magnitude at most~$\pi$.

If $\alpha$ lies in an $R$--neighborhood of the union
of $\beta$~and~$\gamma$, then by Lemma~\ref{lem:OneSideThin}
setting $\delta=2R$ completes the proof.
Otherwise, let $[x,y]$ be the maximal subsegment of~$\alpha$ which lies
outside the open $R$~neighborhood of $\beta \cup \gamma$.
Setting
$m= \bigl\lfloor d(x,y) \big/ 2R \bigr\rfloor$, we can choose points
$x_1,\ldots,x_m$ on $[x,y]$ so that the open balls $B_i=\ball{x_i}{R}$
are disjoint.
But each ball contains negative curvature with magnitude at
least~$\theta$.  So there can be at most $\pi / \theta$
such balls.  Therefore $d(x,y)$ is bounded in terms of
$R$~and~$\theta$.  The conclusion now follows
from Lemma~\ref{lem:OneSideThin}.
\end{proof}

\begin{figure}[ht!]
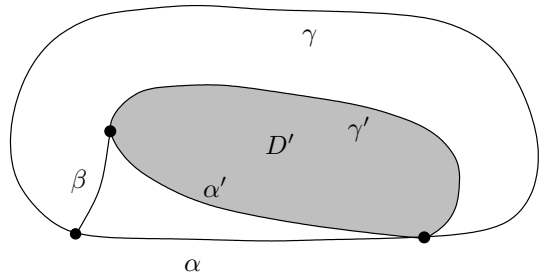

\drawruffles
\caption{If $(D, \alpha)$ is ruffled and $\beta$ is short, then
$(D', \alpha')$ is also ruffled.}
\label{fig:Ruffles}
\end{figure}

Figure~\ref{fig:Ruffles} illustrates a diagram~$D$ which is ruffled
along a boundary arc~$\alpha$
and a subdiagram~$D'$ so that $\alpha$, $\beta$, and $\alpha'$ form a geodesic
triangle.  If the side~$\beta$ is sufficiently short, then the following lemma
states that $D'$ is ruffled along~$\alpha'$.
As in the proof of
Proposition~\ref{prop:ThinTriangles}, the reason is that the Combinatorial
Gauss--Bonnet Theorem provides a bound on the total amount of negative curvature
inside the triangle.

\begin{lem}\label{lem:RuffledTriangle}
Given positive constants $R$,~$\theta$, and~$\delta$, there is a
constant $R'=R'(R,\theta,\delta)$ such that the following property
holds.
Let $D$ be a piecewise Euclidean $\CAT(0)$ disc diagram
whose boundary cycle is a concatenation $\alpha\gamma$.
Let $D'$ and~$\Delta$ be subdiagrams of\/~$D$ with disjoint interiors
such that $\Delta$ is triangular with sides $\alpha$,~$\beta$, and~$\alpha'$
and such that the boundary cycle of\/~$D'$ is a concatenation
$\alpha'\gamma'$ as illustrated in Figure~\ref{fig:Ruffles}.
If\/ $(D,\alpha)$ is $(R,\theta)$--ruffled and $\beta$ has length
at most~$\delta$,
then $(D',\alpha')$ is $(R',\theta/2)$--ruffled.
\end{lem}

\begin{proof}
If $\alpha$ has length at most $\delta$, then we are done since $\alpha'$
then has length at most $2\delta$, and hence $(D',\alpha')$ is
trivially $(\delta,\theta)$--ruffled.  Thus we may assume that
the length of~$\alpha$ is more than~$\delta$.

Let $x$ be the common endpoint of $\alpha$ and~$\alpha'$.
For
\[
   1\le i \le m
     = \bigl\lfloor (\ell(\alpha) - \delta) \big/ 2R \bigr\rfloor,
\]
let $x_i$ be the point on $\Image(\alpha)$ at a distance $(2i-1)R$
from~$x$. Then the open balls
$B_i = \ball{x_i}{R}$ are pairwise disjoint.
Since $\ell(\alpha) \ge 2mR + \delta$ and $\ell(\beta) \le \delta$,
it follows that the image of~$\beta$ does not intersect any~$B_i$.
Therefore, if $B_i$ intersects $\Image(\gamma)$, it must also
intersect $\Image(\gamma')$.
Since $(D,\alpha)$ is $(R,\theta)$--ruffled,
if $B_i$ does not intersect $\Image(\gamma')$, then $B_i$ contains
a vertex~$y_i$ with $\kappa_D(y_i) \le -\theta$.
If $y_i$ is in the interior of~$D'$, then
$\kappa_{D'}(y_i) = \kappa_D(y_i)$.  Similarly,
if $y_i \in \Delta \setminus \Image(\alpha')$, then
$\kappa_{\Delta}(y_i) = \kappa_D(y_i)$.
If $y_i$ lies on the image of~$\alpha'$, it is a vertex in both
$\Delta$ and~$D'$ satisfying
\[
   \kappa_D(y_i) = \kappa_{\Delta}(y_i) + \kappa_{D'}(y_i).
\]
Call $B_i$ \emph{defiant} if it contains a vertex
$y_i\in \Delta$ with $\kappa_{\Delta}(y_i) \le -\theta/2$.
As in the proof of Proposition~\ref{prop:ThinTriangles}, the triangular
diagram~$\Delta$ contains negative curvature with total magnitude
at most~$\pi$.
So at most $2\pi / \theta$ of the $B_i$'s are defiant.
Each nondefiant ball~$B_i$ either intersects $\Image(\gamma')$ or contains
a vertex $y_i \in D'$ with $\kappa_{D'}(y_i) \le -\theta/2$.

Now fix any point~$p'\in\Image(\alpha')$.  Since the endpoints of $\alpha$
and~$\alpha'$ are at most~$\delta$ apart, $p'$ is within a distance~$\delta$
of some point $p\in \Image(\alpha)$.  But $p$ is within $\delta + 3R$ of
the center~$x_i$ of some ball~$B_i$.
This $x_i$ is within a distance $4R\pi / \theta$ of the center~$x_j$
of some nondefiant ball~$B_j$.
Thus $p'$ is within a total distance $2\delta + 4R + (4R\pi / \theta)$
of either the image of~$\gamma'$ or some vertex $y_j\in D'$ with
$\kappa_{D'}(y_j) \le -\theta/2$.
In other words, $(D',\alpha')$ is $(R',\theta/2)$--ruffled,
where $R' = 2\delta + 4R + (4R\pi / \theta)$.
\end{proof}

%%%%%%%%%%%%%%%%%%%%%%%%%%%%%%%%%%%%%%%%%%%%%%%%%%%%%%%%%%%%%%%%%%%%%%%%%%
\section{Preflats in reduced disc diagrams}
\label{sec:Preflats}
%%%%%%%%%%%%%%%%%%%%%%%%%%%%%%%%%%%%%%%%%%%%%%%%%%%%%%%%%%%%%%%%%%%%%%%%%%

Let $\phi\co D \to X$ be a reduced disc diagram
where $X$ is a piecewise Euclidean $\CAT(0)$ $2$--complex.
When $X$ contains flat planes, we will frequently be
interested in examining the preimages under~$\phi$
of these various planes.
These preimages of flats are especially useful in the presence of the
Isolated Flats Property since distinct flats in~$X$ have only a small
intersection.

\begin{figure}[ht!]
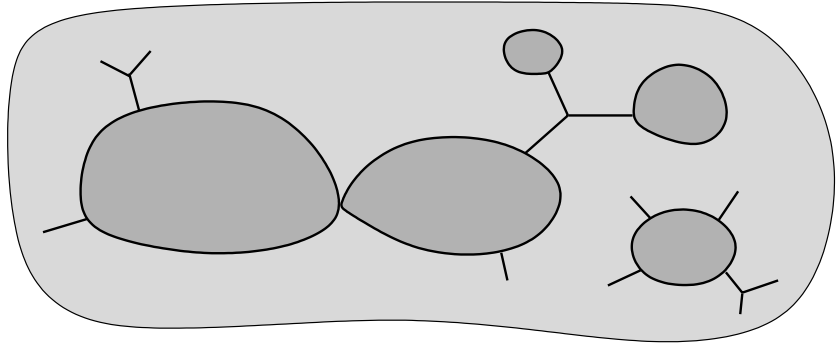

\drawpreflats
\caption[Some preflats in a disc diagram.]%
{A reduced disc diagram $\phi\co D \to X$ and the preimage
under~$\phi$ of some flat plane~$F$ in~$X$.
The preimage shown is disconnected
and contains exactly five preflats.}
\label{fig:preflats}
\end{figure}

Each component of such a preimage is
composed of various topological discs and isolated edges glued together in a
treelike fashion as illustrated in Figure~\ref{fig:preflats}.
It is frequently more convenient to deal with these
discs individually rather than considering the entire preimage as a unit.
We call each such disc a \emph{preflat}.  The following definition
makes this notion more precise.

\begin{defn}[Preflat]\label{def:Preflat}
Let $\phi\co D\to X$ be a reduced disc diagram where $X$~is a $\CAT(0)$
$2$--complex with the Isolated Flats Property.
A \emph{preflat}~$\preflat$ is the closure of a
connected component of the interior of $\phi^{-1}(E)$ for some flat
plane~$E$ in~$X$.
\end{defn}

In this section we see that both the presence and the absence of preflats in
disc diagrams provides a source of ruffles.  In
Proposition~\ref{prop:AlongGeodesic} we show that if a geodesic segment
occurs along the boundary of a reduced disc diagram $D \to X$, then either the
diagram is ruffled along the geodesic
or some preflat in~$X$ comes
close to the geodesic.  So in this case, the absence of preflats
provides ruffles. In contrast, we see in Proposition~\ref{prop:AroundFlat} that
in the presence of the Isolated Flats Property
preflats themselves are surrounded by ruffles.

\begin{prop}\label{prop:AlongGeodesic}
Let $X$ be a proper, cocompact piecewise Euclidean $\CAT(0)$ $2$--complex.
For each $M_0 >0$, there are positive constants $R = R(M_0,X)$
and $\theta_0=\theta_0(M_0,X)$ so that the following property holds
for any positive $\theta \le \theta_0$.

Let $\phi\co D\to \hat{X}$ be a reduced disc diagram
where $\hat{X}$ is a subdivision of\/~$X$.
Suppose the boundary cycle of\/~$D$ decomposes as a concatenation of
two paths $\gamma$~and~$\alpha$ such that $\phi\of\gamma$ is a geodesic
in~$\hat{X}$.  Then either
\begin{enumerate}
   \item \label{item:ruffled} the pair $(D,\gamma)$ is
         $(8R,\theta)$--ruffled, or
   \item \label{item:FlatDiscII} some preflat~$\preflat$
         with inscribed radius at least $M_0$
         intersects the $R$--neighborhood of\/ $\Image(\gamma)$,
         while $\phi(\preflat)$
         is not contained in the $(R+M_0)$--neighborhood
         of\/ $\Image({\phi\of\gamma})$.
\end{enumerate}
\end{prop}

The proof of the preceding proposition uses several lemmas which we will
state and prove before giving the proof of the proposition.
Before we get to these lemmas, let us reflect for a moment on why it
is useful to allow subdivisions of~$X$ in the statement above.
The complex~$X$ has a fixed cell structure such that the group of
combinatorial isometries of~$X$ acts cocompactly.
In practice, we need to consider disc diagrams arising from nullhomotopies
of piecewise geodesic loops in~$X$ which do not respect the given cell
structure of~$X$.  In general, a subdivision of~$X$ is required
to ensure that a particular piecewise geodesic loop lies inside
the $1$--skeleton of~$X$.  We emphasize that
the constants $R$ and $\theta_0$ obtained in the proposition depend only on
the cell structure of~$X$ and not on the particular subdivision~$\hat{X}$
in question.

Our first lemma does not use the Isolated Flats Property
or even nonpositive curvature.

\begin{lem}\label{lem:ExistsFlat}
Let $X$ be a proper, cocompact, piecewise Euclidean $2$--complex.
For every $m>0$ there is an $n=n(m,X)$ so that if\/ $B$~is a flat disc
in~$X$
of radius~$n$, then the central subdisc~${B}'$ of radius~$m$ lies in
some flat plane of\/~$X$.
\end{lem}

\begin{proof}
Our argument is a combinatorial version of the diagonal argument in
the proof of
Lemma~\ref{lem:embedball}.
Suppose, by way of contradiction, that there is a number $m>0$ so that
for each integer $n\ge m$ there is a flat disc~$B_n$ of radius~$n$
whose central subdisc~${B}'_n$ of radius~$m$ does not lie
entirely within any flat
plane of~$X$.  Consider the sequence of embeddings $B_n\inclusion X$.
By composing each embedding with a suitable isometry of~$X$,
we may assume that the center~$p_n$ of every disc~$B_n$ lies in a common
finite subcomplex~$K$.

We proceed by induction to construct a sequence
\[
  \emptyset = C_0 \subseteq C_1 \subseteq C_2 \subseteq \dotsb
\]
of subcomplexes of~$X$ and a sequence
\[
  \N = \Set{S}_0 \supseteq \Set{S}_1 \supseteq \Set{S}_2 \supseteq \dotsb
\]
of infinite sets
such that $C_i$ lies inside each disc of the family
$\set{B_n}{n\in \Set{S}_i}$
and such that the union $\bigcup_i C_i$ is a flat plane.
Suppose $\Set{S}_{i-1}$ and $C_{i-1}$ have already been constructed.
Since $X$~has finitely many isometry types of $2$--cells, the cells
of~$X$ have a universally bounded radius.  It follows that
for each sufficiently large $n\in \Set{S}_{i-1}$ there is a $2$--cell inside
$B_n \setminus C_{i-1}$ which is a minimal distance~$R_{i-1}$
from the center~$p_n$.
The chosen $2$--cells all lie in a bounded (hence compact) closed
neighborhood of~$K$.  So there is an infinite set
$\Set{S}_i \subseteq \Set{S}_{i-1}$
such that for $n\in \Set{S}_i$ the chosen $2$--cells coincide in a single
cell~$e_i$.  Let $C_i = C_{i-1} \cup \bar e_i$.

Notice that for $n\in \Set{S}_i$ the disc of radius~$R_i$
in~$B_n$ centered
at~$p_n$ is contained in~$C_i$.  Furthermore, as $i$ tends to infinity,
the radii~$R_i$ become arbitrarily large.  It follows that the union
$\bigcup_i C_i$ is a flat plane~$E$.  If we choose $n_i\in \Set{S}_i$,
then $B_{n_i}$ is a disc of radius~$n_i$ whose central subdisc
of radius~$R_i$ lies in the flat plane~$E$, which is a contradiction.
\end{proof}

The following lemma gives a lower bound on the size of certain turning angles
occuring when a geodesic passes through a vertex.

\begin{lem}\label{lem:TurningAngle}
Let $X$ be a proper, cocompact piecewise Euclidean $\CAT(0)$
$2$--complex.
For each $R \ge 0$, there is a positive constant
$\theta = \theta(R)$ with the following
property.

Let $x$,~$y$, and~$z$ be distinct vertices of~$X$
such that $y$ lies on the geodesic $[x,z]$,
and suppose that each of the distances
$d(x,y)$ and $d(y,z)$ are less than~$R$.
Let $\pi \co X \setminus \{y\} \to \Lk(y,X)$
be the radial projection onto $\Lk(y,X)$.
Then any immersed path in $\Lk(y,X)$ connecting $\pi(x)$ and $\pi(z)$ which has
length strictly greater than $\pi$ actually has length at least $\pi + \theta$.
\end{lem}

\begin{proof}
Since $X$ is proper and cocompact
we can bring the points $x$,~$y$, and~$z$ into
a fixed compact set $K = K(R)$ by
applying a combinatorial isometry.
But $K$ contains only finitely many vertices.
So for each~$R$, only finitely many choices of $x$,~$y$,
and~$z$ concern us.  Consequently it suffices to prove the result
for a fixed triple of points $x$,~$y$, and~$z$.

Since $\Lk(y,X)$ is a finite metric graph, it contains only finitely many
immersed paths which connect $\pi(x)$ and $\pi(z)$ and which have length less
than $2\pi$.
Thus the paths with length greater than~$\pi$ have a minimum length
as desired.
\end{proof}

Recall that the Monodromy Theorem from complex analysis states that
an analytic continuation of an analytic function along a curve
depends only on the homotopy class of the curve in question
\cite{Rudin}.
The next lemma is an easy consequence of the Monodromy Theorem.

\begin{lem}\label{lem:Monodromy}
Let $S$ be a simply connected nonpositively curved surface
which is locally isometric to the Euclidean plane~$\E^2$.
Then $S$ admits a local isometry into~$\E^2$.
\end{lem}

\begin{proof}
By hypothesis, each point $x \in S$ is contained in an open disc~$D_x$
which admits an isometric embedding into~$\E^2$.
Fix a specific embedding of a single disc~$D_p$.
For each overlapping disc~$D_q$, we can compose the given embedding
$D_q \inclusion \E^2$ with an isometry of~$\E^2$ to make it agree
with the chosen embedding $D_p \inclusion \E^2$ on the overlap
$D_p \cap D_q$.
This procedure provides an analytic function
$D_p \cup D_q \inclusion \E^2$ extending the original map
$D_p \inclusion \E^2$.
Continuing in this fashion, a map $S \to \E^2$
can be defined by taking an analytic extension
of the original embedding $D_p \inclusion \E^2$ along various paths.
By the Monodromy Theorem,
these continuations fit together to give a globally defined
local isometry $S \to \E^2$.
\end{proof}

The following lemma roughly states that if a geodesic occurs as part of the
boundary of a disc diagram, then either
the diagram is ruffled along the geodesic
or some point of the geodesic has a neighborhood isometric to
a Euclidean half-disc.

\begin{lem}\label{lem:AlongGeodesic}
Let $X$~be a piecewise Euclidean $\CAT(0)$ $2$--complex, and let
$\phi\co D\to X$ be a reduced disc diagram.
Suppose the boundary path of\/~$D$ is a concatenation $\alpha\gamma$
where $\gamma$~is a geodesic in~$D$.
Let $p=\gamma(t)$ be an arbitrary point in the image of~$\gamma$.
If we let $B = \ball{p}{R}$, then either
\begin{enumerate}
   \item \label{item:NegCurvI} $B$~contains a vertex with negative
         curvature, or
   \item \label{item:HitsBoundaryI} $B$~intersects the image of~$\alpha$, or
   \item \label{item:FlatDiscI} $B$~is isometric to a Euclidean half-disc.
\end{enumerate}
\end{lem}

\begin{proof}
Suppose conditions (\ref{item:NegCurvI})~and~(\ref{item:HitsBoundaryI})
fail for~$B$.
Then the curvature at every vertex of~$B$ is zero.
So $B$~is locally flat.
In other words, each vertex of~$\interior{B}$
has a neighborhood isometric to a Euclidean disc, while
each vertex of~$\boundary B$ has a neighborhood
isometric to a Euclidean half-disc.

Since $B$~is a metric ball in a $\CAT(0)$ space, $B$~is a convex subspace
of~$D$.
So $B$ is simply connected, and hence
admits a local isometry into the
Euclidean plane by Lemma~\ref{lem:Monodromy}.
Since $\interior{B}$ is convex,
this local isometry is actually an isometry
from~$B$ to a convex set in the Euclidean plane.  Now $B$~is easily
seen to be isometric to a flat Euclidean half-disc of radius~$R$.
\end{proof}

We are now ready to prove Proposition~\ref{prop:AlongGeodesic}.

\begin{proof}[Proof of Proposition~\ref{prop:AlongGeodesic}]
By Lemma~\ref{lem:ExistsFlat}, we can choose $R=R(M_0,X)$
sufficiently large that
for any flat disc of radius~$R$ in~$X$ the central subdisc of
radius~$2M_0$ lies in a flat plane.
If some ball $B=\bigball{\gamma(t)}{2R}$ in~$D$ is isometric to
a Euclidean half-disc, then $B$~maps isometrically to~$X$ under~$\phi$.
Furthermore, $B$~contains a Euclidean disc of
radius~$R$ whose central subdisc of radius~$2M_0$ intersects the
$R$--neighborhood of~$\gamma$ but is not contained in the
$(R+M_0)$--neighborhood of~$\gamma$, as illustrated in
Figure~\ref{fig:FlatDisc}.
By our choice of~$R$, this subdisc lies in some preflat~$\preflat$.
Since $B$ maps isometrically to~$X$ under~$\phi$, it follows that
$\phi(\preflat)$ is not contained in the $(R+M_0)$--neighborhood
of $\Image(\phi\of\gamma)$.
\begin{figure}[ht!]
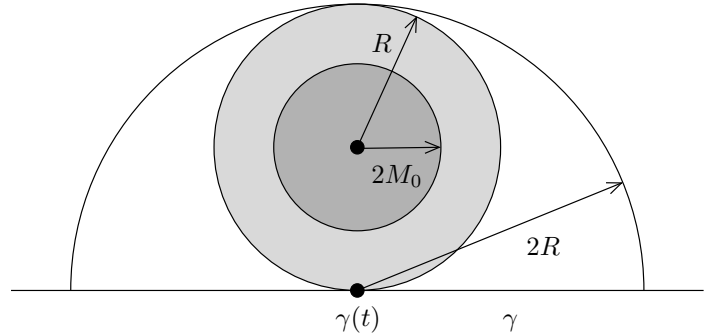

\drawflatdisc
\caption{The inner disc intersects the $R$--neighborhood
   of~$\gamma$ but is not contained in the $(R+M_0)$--neighborhood
   of $\gamma$.}
\label{fig:FlatDisc}
\end{figure}

We will now assume that
no such ball is isometric to a flat half-disc.
By Lemma~\ref{lem:AlongGeodesic}, each ball
$\bigball{\gamma(t)}{2R}$ either intersects the image
of $\alpha$, or contains
a vertex~$v$ with negative curvature.  It remains to find a positive constant
$\theta(M_0,X)$ so that (\ref{item:ruffled}) holds.

Let $v$ be a vertex in the interior of~$D$ with $\kappa(v)<0$.
Then $2\pi+\abs{\kappa(v)}$ is the length of
some locally geodesic loop~$\ell$ in
$\Lk \bigl( \phi(v), \hat{X} \bigr)$ since $\phi$ is
a reduced map.
If $\phi(v)$ is not a vertex of~$X$, then $\ell$ has length $n\pi$ for some
integer $n > 2$.
On the other hand, if $\phi(v)$ is a vertex of~$X$, then, since $X$ has
only finitely many isometry types of cells, there is a
positive constant~$\theta_1$, depending only on the space~$X$,
such that $\ell$ has length at least $2\pi + \theta_1$.
In either case, we see that
$\kappa(v) \le -\theta_1$.

Now suppose $v$ is a vertex on~$\gamma$ with $\kappa(v)<0$.
We may assume that $\phi(v)$ is a vertex of~$X$ since otherwise
$\kappa(v)$ would be a negative multiple of~$\pi$.
Let $x$~and~$y$ be points on~$\gamma$ at a distance $4R$ on either side
of~$v$, or the endpoints of~$\gamma$ if $v$ is too close to one of the ends.
If either $x$ or~$y$ is within $2R$ of~$\alpha$, then $v$ is
within $6R$ of~$\alpha$.
Otherwise, there are vertices $z$ and~$w$ with negative
curvature so that $d(x,z)$ and $d(y,w)$ are each less than $2R$.
If either $z$ or~$w$ is in the interior of~$D$, then its curvature is less
than~$-\theta_1$, and hence $v$ is within~$6R$ of a vertex
with curvature
less than~$-\theta_1$.
Now suppose that both $z$ and~$w$ lie on~$\gamma$.
If either point does not map under~$\phi$ to a vertex of~$X$, then the
curvature
at that point is at least~$\pi$ as above, so $v$ is within a distance~$6R$ of a
point with curvature at most~$-\pi$.
Finally assume that both $z$ and~$w$ map under~$\phi$ to vertices of~$X$.
Since $v$ also maps to a vertex of~$X$, and the distances
from $v$ to $z$ and~$w$ are each at most $6R$, then the curvature at~$v$ has
magnitude at least $\theta_2(R,X)$ by Lemma~\ref{lem:TurningAngle}.

We have now shown that every vertex with negative curvature either in
the interior of~$D$ or on~$\gamma$ is within $6R$ of either $\alpha$
or a vertex with negative curvature of magnitude at least
\[
  \theta_0=\min\{\theta_1,\theta_2\} >0.
\]
Since every point on~$\gamma$ is within $2R$ of either $\alpha$
or a vertex
with negative curvature (which does not lie on~$\alpha$),
the pair $(D,\gamma)$ must be $(8R,\theta_0)$--ruffled.
The result now holds for any positive $\theta \le \theta_0$.
\end{proof}

The second proposition of this section is the following result,
stating roughly that in the presence of the Isolated Flats Property
preflats are surrounded by ruffles.  More specifically, if the
interior of a preflat is removed from a disc diagram, then
the resulting diagram will be ruffled along the new
boundary cycle.

\begin{prop}\label{prop:AroundFlat}
Let $X$ be a proper, cocompact piecewise Euclidean $\CAT(0)$ $2$--complex
with the Isolated Flats Property.
There are positive constants $R(X)$ and $\theta(X)$
satisfying the following property.
Let $\phi\co D\to \hat{X}$ be a reduced disc diagram
where $\hat{X}$ is a subdivision of~$X$, and let
$D'=D \setminus \interior\preflat$ for some preflat~$\preflat$ in~$D$.
Then the pair $(D',\boundary\preflat)$
is $(R,\theta)$--ruffled.
\end{prop}

The proof of this proposition has many similarities with the proof
of Proposition~\ref{prop:AlongGeodesic}.
Before giving this proof, we will state and prove some related
lemmas describing important properties of preflats.
First we need to define the notion of a \emph{corner} of
a subdiagram.

\begin{defn}
Let $D$~be a disc diagram, and let $S$~be a subdiagram.  Let $v$ be a vertex
of~$D$.  A \emph{corner of\/~$S$ at~$v$} is a set of corners of $2$--cells
of~$S$ whose closure corresponds to an entire connected component
of $\Lk(v,S)$ for some $0$--cell~$v$.
An \emph{exterior corner of\/~$S$ in~$D$} is a set
of corners of $2$--cells of $D \setminus S$ whose closure corresponds to a
connected arc in $\Lk(v,D)$ which intersects $\Lk(v,S)$ only at the
endpoints
of the arc.  The \emph{angle} of a corner is the sum of the angles of its
elements.
\end{defn}

The first basic property of preflats that we need is their
convexity, which is established in the following lemma.

\begin{lem}\label{lem:BoundaryIsGeodesic}
Let $X$~be a piecewise Euclidean $\CAT(0)$ $2$--complex,
let $\phi\co D\to X$
be a reduced disc diagram, and let $E$~be any flat plane in~$X$.
Then each component of\/ $\phi^{-1}(E)$ is a convex subspace of\/~$D$.
Furthermore, each preflat~$\preflat$ is also convex in~$D$.
\end{lem}

\begin{proof}
It suffices to show that every exterior corner of $\phi^{-1}(E)$ in~$D$
has angle at least~$\pi$, since then $\phi^{-1}(E)$ is locally convex.
Let $v$~be a vertex in $\boundary \phi^{-1}(E)$,
let $I$~be a closed interval of~$\R$,
and let $I\to \Lk(v,D)$ be a local isometry such that the intersection
$\Image(I) \cap \Lk \bigl( v, \phi^{-1}(E) \bigr)$
consists of only the endpoints of $\Image(I)$.
Let $\Lk(v,D) \to \Lk \bigl( \phi(v), X \bigr)$
be the map induced by~$\phi$.  Then the composition
$I \to \Lk(v,D) \to \Lk \bigl( \phi(v), X \bigr)$ is a local isometry
whose image intersects $\Lk \bigl( \phi(v), E \bigr)$
only at its endpoints.  Notice that $\Lk \bigl( \phi(v), E \bigr)$ has
diameter~$\pi$. So $I$ has length at least~$\pi$, since otherwise there would
be
a locally geodesic loop $C\to \Lk \bigl( \phi(v), X \bigr)$ of length
less than~$2\pi$ contradicting the Link Condition for~$X$
(Theorem~\ref{thm:LinkCondition}).
It now follows easily that every exterior corner of $\phi^{-1}(E)$
in~$D$ has angle at least~$\pi$.
\end{proof}

The next lemma is quite similar to Lemma~\ref{lem:AlongGeodesic}.
We will show that each preflat is either surrounded by ruffles or
adjacent to a large Euclidean half-disc.

\begin{lem}\label{lem:AroundFlat}
Let $\phi\co D\to X$ be a reduced disc diagram where
$X$~is a piecewise Euclidean $\CAT(0)$ $2$--complex, and let $\preflat$~be
a preflat.
Let $D' = D \setminus \interior\preflat$, and let $p$ be any point
on~$\boundary \preflat$.
If\/ $B'$~is the open ball of radius~$R$ in~$D'$ centered at~$p$,
then either
\begin{enumerate}
   \item \label{item:NegCurvIII} $B'$~contains a vertex with negative
       curvature in~$D'$, or
   \item \label{item:HitsBoundaryIII} $B'$~intersects $\boundary D$, or
   \item \label{item:FlatDiscIII} $B'$~is isometric to a flat half-disc
       of radius~$R$.
\end{enumerate}
\end{lem}

\begin{proof}
Suppose conditions (\ref{item:NegCurvIII})~and~(\ref{item:HitsBoundaryIII})
fail for~$B'$.
Then each vertex on the interior of~$B'$ has a neighborhood isometric
to a Euclidean disc.

Now choose an arbitrary vertex $v$ in $B'\cap \boundary \preflat$.
By Lemma~\ref{lem:BoundaryIsGeodesic}, every exterior corner of~$\preflat$
in~$D$ has angle at least~$\pi$.  Since the curvature
$\kappa_{D'}(v)$ is zero, $v$~must have a neighborhood in~$B'$
isometric to a Euclidean half-disc.  Furthermore,
$\boundary \preflat$~intersects~$B'$
in a collection of disjoint segments which are geodesic in~$B'$
as well as in~$D$.
In particular, note that $\preflat$~does not lie entirely inside
the open ball~$B$ of radius~$R$ in~$D$ centered at~$p$.

By Lemma~\ref{lem:Monodromy} if $B'$~is simply connected, it admits
a local isometry into the Euclidean plane.
Recall that any metric ball in a $\CAT(0)$ space is convex.
We will show that $B'$ is contractible inside~$D'$, and hence lifts
isometrically
to the universal cover $\tilde D'$, where it is a metric ball in a
$\CAT(0)$ space, and must therefore be convex and contractible.

Suppose $B'$~is not contractible in~$D'$.  Notice that $\pi_1(D')$ is cyclic
generated by the boundary cycle of~$D$.
So $B'$~must contain a simple loop~$\ell$ which encloses~$\preflat$.
Since $B'\subset B$, it follows that $\ell$ also lies inside~$B$, which is
contractible.  Therefore $\preflat$~lies entirely within~$B$, which is
a contradiction.  It now follows that
$B'$~lifts isometrically to a convex subspace of~$\tilde D'$,
and therefore $B'$~is itself $\CAT(0)$.  Being simply connected
and locally flat, $B'$~admits
a local isometry to the plane, which must actually be an isometric
embedding since $B'$~is $\CAT(0)$.

Since $B'$~and~$\preflat$ are both convex in~$D$, their intersection
is connected.
Therefore $B'$~is isometric to a Euclidean half disc
of radius~$R$.
\end{proof}

To complete the proof of Proposition~\ref{prop:AroundFlat}
we need to show that in the presence of the Isolated Flats Property,
the third case of the preceding lemma cannot occur.

\begin{proof}[Proof of Proposition~\ref{prop:AroundFlat}]
By Theorem~\ref{thm:FlatTriplane}, we can choose~$R$ sufficiently large
that $X$~does not contain an isometrically embedded copy of
the space~$T_R$ (in the notation of Definition~\ref{def:triplane}).
Choose $\theta>0$ so that for each vertex $v\in X^{(0)}$
\begin{enumerate}
   \item every combinatorial reduced path in $\Lk(v,X)$ with length
      greater than~$\pi$ actually has length at least $\pi+\theta$, and
   \item every locally geodesic loop in $\Lk(v,X)$ with length
      greater than~$2\pi$ actually has length at least $2\pi+\theta$.
\end{enumerate}
Such a $\theta$ exists since $X$~has finitely many isometry types
of $2$--cells.

Now fix a reduced disc diagram $\phi\co D\to \hat{X}$ where
$\hat{X}$ is
a subdivision of~$X$, and let
$D'= D \setminus \interior\preflat$ for some preflat~$\preflat$ in~$D$.
Notice that if $v$ is a vertex in $D' \setminus \boundary D$ with
$\kappa_{D'}(v) < 0$, then by our choice of~$\theta$
in fact $\kappa_{D'}(v) \le -\theta$.
Thus by Lemma~\ref{lem:AroundFlat}, either $(D',\boundary\preflat)$ is
$(R,\theta)$--ruffled or there is some point $p\in\boundary\preflat$
so that the ball $B=\ball{p}{R}$ in~$D'$ is isometric to a Euclidean
half disc.

Notice that $\phi$~maps the boundary edge of~$B$ isometrically to a geodesic
segment in a flat plane~$E$ of~$X$.  Furthermore, $\phi^{-1}(E)$
does not have a local cut point anywhere along $\boundary B$,
since otherwise $\phi^{-1}(E)$ would fail to be convex,
contradicting Lemma~\ref{lem:BoundaryIsGeodesic}.
It follows that $\phi(B)$ intersects~$E$ only along its boundary edge.
But then $X$~contains an isometrically embedded copy of the space~$T_R$,
which is a contradiction.
\end{proof}

%%%%%%%%%%%%%%%%%%%%%%%%%%%%%%%%%%%%%%%%%%%%%%%%%%%%%%%%%%%%%%%%%%%%%%%%%%
\section{$2$--complexes with isolated flats have the Relatively Thin Triangle
Property}
\label{sec:RelThinTriangles}
%%%%%%%%%%%%%%%%%%%%%%%%%%%%%%%%%%%%%%%%%%%%%%%%%%%%%%%%%%%%%%%%%%%%%%%%%%

In this section we show that the Isolated Flats Property implies
the Relatively Thin Triangle Property, establishing
(1)~$\implies$~(2) of Theorem~\ref{thm:2dEquivalent}.
The proof combines the two main
propositions from the previous section about preflats and ruffles.

\begin{thm}\label{thm:IFP=>ThinTriangles}
Let $X$ be a proper, cocompact piecewise Euclidean $\CAT(0)$ $2$--complex.
If\/ $X$ satisfies the Isolated Flats Property, then
$X$ also satisfies the Relatively Thin Triangle Property.
\end{thm}

\begin{lem}\label{lem:RelThinTriangles}
Let $X$ be a proper, cocompact piecewise Euclidean $\CAT(0)$ $2$--complex with
the Isolated Flats Property.
There are positive constants $\mu$ and~$\delta$ such that the following
property
holds.
Let $\Delta(x,y,z)$ be a geodesic triangle in~$X$.
Form a subdivision~$\hat{X}$ of\/~$X$ such that
the sides of~$\Delta$ lie in the
$1$--skeleton of\/~$\hat{X}$.
Let $\phi \co D \to \hat{X}$ be a reduced disc diagram for $\Delta$.
Suppose $D$ contains a preflat~$P$ which maps into a flat~$F$ under~$\phi$.
If\/ $P$ has inscribed radius at least~$\mu$, then the triangle~$\Delta$ is
$\delta$--thin relative to the flat~$F$.
\end{lem}

\begin{proof}
To show that $\Delta$ is $\delta$--thin relative to~$F$ for some~$\delta$,
it suffices to show that each side of~$D$ lies in a $\delta$--neighborhood of
the
union of the other two sides and~$P$.
If we let $D'$ denote $D \setminus \interior{P}$, then
Proposition~\ref{prop:AroundFlat} provides constants $R$ and~$\theta$
so that $(D',\boundary P)$ is $(R,\theta)$--ruffled.

Suppose the boundary of~$P$ has length~$L$.  For each~$i$ with
\[
   1 \le i \le m = \lfloor L / 2R \rfloor,
\]
let $x_i$ be a point on $\boundary P$ so that the balls
$B_i = \ball{x_i}{R}$ in~$D'$ are pairwise disjoint.
Call $B_i$ \emph{defiant} if it contains a vertex~$v$ with
$\kappa_{D'} (v) \le -\theta$.
By the Combinatorial Gauss--Bonnet Theorem, the negative curvature
inside $D'$
has total magnitude at most~$\pi$.  So at most $\pi / \theta$ of the~$B_i$ are
defiant, and each nondefiant~$B_i$ intersects $\boundary D$.

We will now show that if~$P$ has a sufficiently large inscribed radius, then
each of the three sides of the triangular diagram~$D$ intersects at least one 
of
the~$B_i$.  For each side of~$D$ which intersects some ball~$B_i$, the set of
all such balls intersecting that side is a contiguous string of nondefiant
balls.  Thus the set of all balls~$B_i$ decomposes into alternate strings of
defiant and nondefiant balls, each string of nondefiant balls containing at
most $\pi / \theta$ balls.

If no side of~$D$ intersects any ball~$B_i$ then every $B_i$ is defiant
as in Figure~\ref{fig:defiance}(a).
In this case, we see that $\boundary P$ has length
less than $\epsilon = 2R ( \pi/\theta + 1)$, which bounds the inscribed radius
of~$P$.

\begin{figure}[ht!]
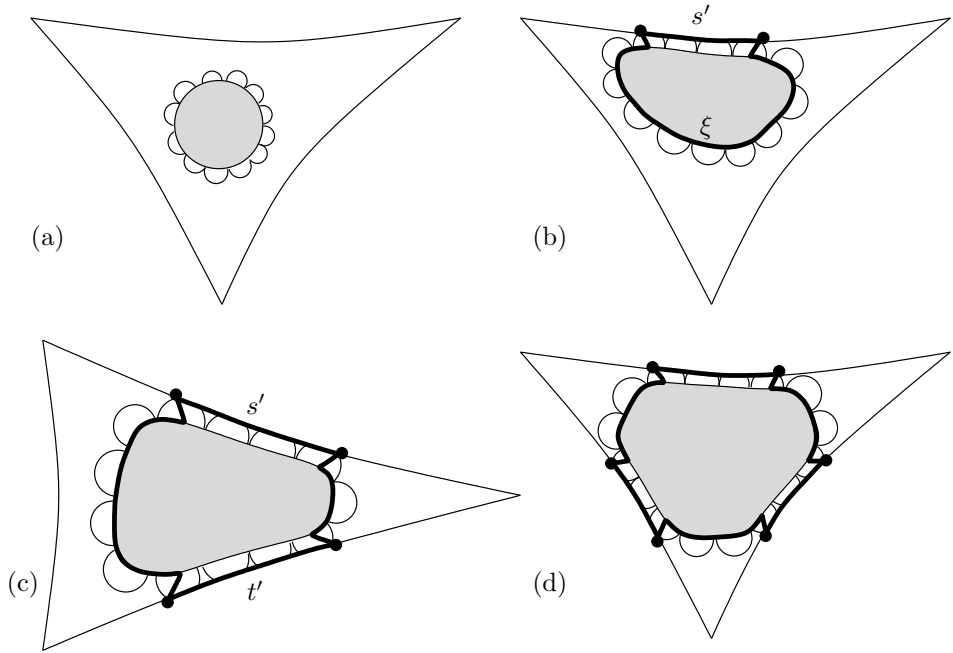

\drawdefiance
\caption[Four cases in the proof of Lemma~\ref{lem:RelThinTriangles}.]%
{(a) No sides of $D$ intersect the balls $B_i$\qua
(b) One side of~$D$ intersects the $B_i$\qua
(c) Two sides of~$D$ intersect the $B_i$\qua
(d) Three sides of~$D$ intersect the $B_i$.
Only in this last case can the preflat have an arbitrarily
large inscribed radius.}
\label{fig:defiance}
\end{figure}

If exactly one side~$s$ of~$D$ intersects the balls~$B_i$, let $s'$ be the
shortest subsegment of~$s$ containing $\bigcup_i (s \cap B_i)$.
The two ends of~$s'$ are connected by a path~$\xi$  which travels around
the ``nondefiant'' side of $\boundary P$ as
illustrated in Figure~\ref{fig:defiance}(b).  Since $\xi$ has length
at most $\epsilon + 4R$, the geodesic~$s'$ also has length
at most $\epsilon + 4R$.
So the inscribed radius of~$P$ is bounded by some constant
which depends only on $\epsilon$ and~$R$.

If exactly two sides $s$ and~$t$ of~$D$ intersect the balls~$B_i$, then as
before, let $s'$ (resp.~$t'$) denote the shortest subsegment of $s$
(resp.~$t$) containing
\[
  \bigcup_i s \cap B_i \quad \Big( \text{resp. }\bigcup_i t \cap B_i \Big).
\]
As in the previous case, the endpoints of $s'$ and~$t'$ are separated by a
distance of at most $\epsilon + 4R$, as shown in
Figure~\ref{fig:defiance}(c).  So again we have a
bound on the inscribed radius of~$P$.

Thus we see that for arbitrarily large inscribed radii, each side of~$D$
intersects some ball~$B_i$.  For each side $s$, let $s'$ be the subsegment
defined as above.  It is clear that the endpoints of these three segments come
in pairs each of which is a distance at most $\epsilon + 4R$ apart as in
Figure~\ref{fig:defiance}(d).
It now follows easily that the triangular diagram~$D$ is
$(\epsilon + 4R)$--thin whenever its inscribed radius is larger than
some constant~$\mu$ depending only on the space~$X$.
\end{proof}

\begin{proof}[Proof of Theorem~\ref{thm:IFP=>ThinTriangles}]
Choose a geodesic triangle $\Delta$ in~$X$, and let $D \to \hat{X}$ be a
reduced disc diagram with $\hat{X}$ a subdivision of~$X$.
By the previous lemma, we have constants $\mu$ and~$\delta$ so that if
$D$ contains a preflat $P$ with inscribed radius at least~$\mu$, then
$\Delta$ is $\delta$--thin relative to some flat~$F$.

Suppose $D$ does not contain such a preflat.  Then by
Proposition~\ref{prop:AlongGeodesic}, there are constants $R$ and~$\theta$,
depending on~$\mu$,
so that for each side~$s$ of~$D$ the pair $(D,s)$ is $(R,\theta)$--ruffled.
But then Proposition~\ref{prop:ThinTriangles} gives that
$D$ is $\delta'$--thin for some constant~$\delta'$ depending only on $R$
and~$\theta$.
\end{proof}

%%%%%%%%%%%%%%%%%%%%%%%%%%%%%%%%%%%%%%%%%%%%%%%%%%%%%%%%%%%%%%%%%%%%%%%%%%
\section{The fellow travelling of quasigeodesics and ruffled geodesics}
\label{sec:Divergence}
%%%%%%%%%%%%%%%%%%%%%%%%%%%%%%%%%%%%%%%%%%%%%%%%%%%%%%%%%%%%%%%%%%%%%%%%%%

Our objective for the rest of this article is to prove
Theorem~\ref{thm:RelativeFTP},
which states that $\CAT(0)$ $2$--complexes with the Isolated Flats Property also
have the Relative Fellow Traveller Property.
In this section, we prove Proposition~\ref{prop:RuffledFTP}, which is
a special case of
Theorem~\ref{thm:RelativeFTP}.
Let $D$ be a $\CAT(0)$ disc diagram whose boundary is a concatenation of a
geodesic and a quasigeodesic.  Proposition~\ref{prop:RuffledFTP} states that
the geodesic and quasigeodesic track close together provided that $D$ is
ruffled along the geodesic.

Note that the conclusion of Proposition~\ref{prop:RuffledFTP}
is stronger than the conclusion of
Theorem~\ref{thm:RelativeFTP}, where we only get that two paths
fellow travel relative to flats.
Recall that a $\delta$--hyperbolic disc diagram is ruffled throughout.  In
particular if such a diagram contains a geodesic on its boundary, then the
diagram will be ruffled along that geodesic.
So from Proposition~\ref{prop:RuffledFTP} we recover a $2$--dimensional version
of Theorem~\ref{thm:GromovFTP} which states that quasigeodesics and geodesics
in a $\delta$--hyperbolic space asynchronously fellow travel.

In fact the structure of our argument in this section is inspired by
Cooper, Lustig, and Mihalik's proof in \cite{ABC91} of this
fellow traveler theorem for $\delta$--hyperbolic spaces.
Their strategy is to first prove an ``exponential
divergence'' theorem \cite[Theorem~2.19]{ABC91}
which roughly states that geodesics diverge at an
exponential rate in
a $\delta$--hyperbolic space.
They then use the fact that this divergence
is \emph{superlinear}
to prove the fellow traveller property in \cite[Proposition~3.3]{ABC91}.

The main tool used in our proof of Proposition~\ref{prop:RuffledFTP}
is Proposition~\ref{prop:QuadraticDivergence}, which is essentially a
``quadratic divergence'' theorem analogous to the exponential
divergence result alluded to above.

At this point, it may be useful to compare our divergence theorem
to some related results in the literature.  A \emph{divergence function}
for a geodesic space
in the sense of \cite{ABC91} is, roughly speaking,
a function that provides a lower bound
on the rate of divergence of all pairs of geodesic rays in the space.
The details of the definition are such that Euclidean space does not admit
an unbounded divergence function.
Papasoglu shows in \cite{Papasoglu95} that a space admitting an unbounded
divergence function admits an exponential divergence function,
and is hence $\delta$-hyperbolic by \cite{ABC91}.
Since a space with isolated flats (that truly contains flats)
fails to be $\delta$-hyperbolic,
it cannot admit a quadratic divergence function in this sense.
The main distinction between Papasoglu's result and Proposition~\ref{prop:QuadraticDivergence}
is that we restrict the type of geodesic rays considered.  Proposition~\ref{prop:QuadraticDivergence}
roughly says that a pair of ``ruffled'' geodesics must diverge
at least quadratically.  However, the proposition makes no conclusion
about arbitrary pairs of rays.

The proof of Proposition~\ref{prop:QuadraticDivergence}
uses the notion of a \emph{broom}, which is a type of disc diagram
that occurs when two geodesics have a common initial segment and then
separate from each other.
After establishing
Proposition~\ref{prop:QuadraticDivergence}, the proof of
\cite[Proposition~3.3]{ABC91}
can be applied almost verbatim to prove Proposition~\ref{prop:RuffledFTP}.

\begin{figure}[ht!]
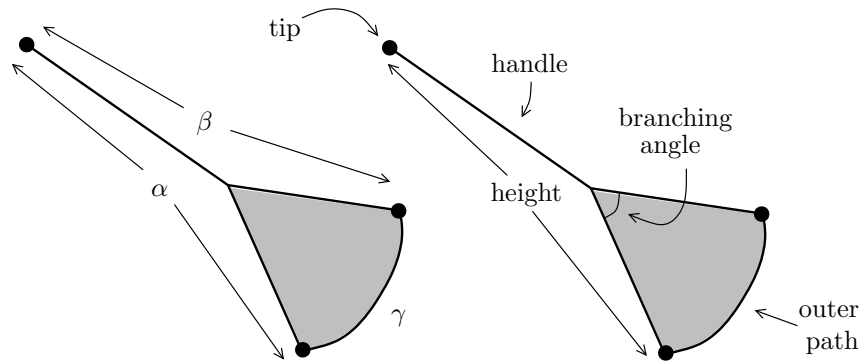

\drawbroom
\caption[Two pictures of a broom with geodesic sides $\alpha$ and~$\beta$.]%
{Two pictures of a broom with geodesic sides $\alpha$ and~$\beta$.
On the right, the different parts of the broom are labeled.}
\label{fig:broom}
\end{figure}

\begin{defn}\label{def:Broom}
A \emph{broom}~$B$ is a piecewise Euclidean $\CAT(0)$ disc diagram whose
boundary is a composition $\alpha\beta\gamma$ where $\alpha$~and~$\beta$
are geodesics, as illustrated in Figure~\ref{fig:broom}.  The \emph{handle}
of~$B$ is the intersection $\alpha \cap \beta$.  The \emph{tip} of~$B$ is the
common endpoint of $\alpha$~and~$\beta$.  The \emph{height} of~$B$ is the
minimum of the lengths of $\alpha$~and~$\beta$.
The \emph{branching angle} of~$B$
is the angle between $\alpha$~and~$\beta$ at the point where they first
separate
(in the degenerate case where
either $\alpha$~or~$\beta$ is equal to the handle,
then the branching angle is defined to be zero).
The \emph{outer path} of~$B$ is the boundary path~$\gamma$.
\end{defn}

\begin{lem}[Linear divergence of brooms]\label{lem:BroomDivergence}
Let $B$~be a broom of height at least~$M$,
branching angle at least~$\theta$, a handle of length at most~$N$,
and an outer path~$\gamma$.
Then
\[
   \ell(\gamma) \ge (M - N) \theta/2.
\]
\end{lem}

\proof
Let $\alpha$~and~$\beta$ be the geodesic sides of~$B$.
Let $\alpha(0)=\beta(0)=p$ be the tip of~$B$, let $q$ be the other end of the
handle,
let $x$ be the common endpoint of $\alpha$ and~$\gamma$,
and let $y$ be the common endpoint of $\beta$ and~$\gamma$.
Let $\Delta$ be the geodesic triangle with vertices $q$,~$x$, and~$y$,
and let $\bar{\Delta}$ be a comparison triangle in the Euclidean plane
with vertices $\bar{q}$,~$\bar{x}$, and~$\bar{y}$.  Let $\bar{\theta}$
denote the angle at the vertex~$\bar{q}$ of~$\bar{\Delta}$.

Since $\alpha$ and~$\beta$ each have length at least~$M$,
we can set $x'=\alpha(M)$ and $y'=\beta(M)$, and let $\bar{x'}$
and~$\bar{y'}$ be the corresponding comparison points on~$\bar{\Delta}$.
By the Law of Cosines,
\[
   d ( \bar{x'},\bar{y'} ) \le d(\bar{x},\bar{y}).
\]
Using the fact that $\sin(\phi) \ge \phi/2$ whenever $0 \le \phi \le \pi/2$,
we conclude that
\begin{align*}
   \ell(\gamma) & \ge d(x,y) \\
                & =   d(\bar{x},\bar{y}) \\
                & \ge d(\bar{x'},\bar{y'}) \\
                & \ge 2(M - N) \sin ( \bar{\theta} /2 ) \\
                & \ge (M - N) \bar{\theta}/2 \\
                & \ge (M - N) \theta/2.   \tag*{\qed}
\end{align*}

The following proposition is a quadratic divergence theorem
analogous to the exponential divergence theorem for $\delta$--hyperbolic spaces
\cite[Theorem~2.19]{ABC91}.
Bridson and Haefliger reformulated this exponential divergence theorem
as an exponential lower bound on the length of any path that a geodesic stays
far away from \cite[Proposition~III.H.1.6]{BH99}.
The precise statement of our quadratic divergence result
is closer to the form given by Bridson and Haefliger.

\begin{prop}[Quadratic divergence]\label{prop:QuadraticDivergence}
For each choice of positive constants $R$ and~$\theta$, there is a quadratic
function $Q\co \R^{+} \to \R$ so that the following property holds.
Let $D$ be any piecewise Euclidean $\CAT(0)$ disc diagram whose boundary
is a concatenation $\alpha\gamma$, where $\gamma$ is a geodesic, and
$(D,\gamma)$ is $(R,\theta)$--ruffled.
Then for any point $p\in \Image(\gamma)$, we have
\[
  \ell(\alpha) \ge Q \bigl( d(p, \Image \alpha) \bigr).
\]
\end{prop}

\begin{proof}
Fix a point $p\in \Image(\gamma)$, and let $r= d(p,\Image \alpha)$.
For convenience, replace $\theta$ with $\min \{\theta,\pi\}$.
We will show the existence of a quadratic function $Q(r)$
independent of our choice of diagram~$D$.
For any point~$q \in D$, let $\Sh(q)$ denote the \emph{shadow} of~$q$
on~$\alpha$, ie, the set
of points~$\alpha(t)$ such that the geodesic from $p$ to~$\alpha(t)$
passes through~$q$.  Notice that each shadow is connected.

Set $Q(r)=0$ for $r < 2R$.  Henceforth we assume
that $r\ge 2R$.
For $1 \le i \le k = \lfloor r / 2R \rfloor$, let $x_i$ be a point
of $\Image(\gamma)$ with $d(p,x_i)=(2i-1)R$.  Notice that the open
balls $B_i = \ball{x_i}{R}$ are pairwise disjoint and do not
intersect the image of~$\alpha$,
as illustrated in Figure~\ref{fig:quadraticdiv}.
Since $(D,\gamma)$ is $(R,\theta)$--ruffled, $B_i$ contains a point~$y_i$
with $\kappa(y_i) \le -\theta$.
\begin{figure}[ht!]
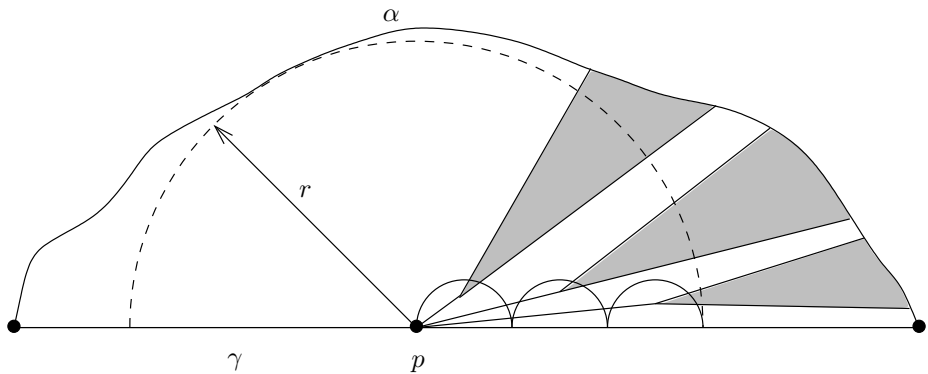

\drawquadraticdiv
\caption{The length of~$\alpha$ is bounded below by a quadratic function
of~$r$.}
\label{fig:quadraticdiv}
\end{figure}

Since $\kappa(y_i) \le -\theta$, we can find
a pair of geodesics $\beta_i$ and~$\beta'_i$
starting at~$p$ and passing
through~$y_i$, where they separate with an angle of~$\theta$.
Since geodesics are extendible to the boundary in any $\CAT(0)$ disc
diagram,
we may assume that the given geodesics continue until they eventually hit
points $z_i$ and~$z'_i$ on the image of~$\alpha$.

Notice that the subdiagram~$B_i$ of~$D$ bounded by $\beta_i$,~$\beta'_i$,
and~$\Sh(y_i)$ is a broom with height at least $2R k$,
branching angle at least~$\theta$, a handle of length at most $2R i$,
and outer path $\Sh(y_i)$.
Applying Lemma~\ref{lem:BroomDivergence} to~$B_i$ gives that
\[
  \ell \bigl( \Sh(y_i) \bigr) \ge R(k-i)\theta.
\]
If we knew that the shadows of the $y_i$ had pairwise disjoint interiors,
then we would be done, since the sum of their lengths satisfies
\[
   \sum_{i=1}^k \ell \bigl( \Sh(y_i) \bigr)  \ge R \theta \sum_{i=1}^k (k-i)
                   %   = R \theta \sum_{j=0}^{k-1} j
                      = R \theta k(k-1) \big/ 2 \text{,}
\]
which is a quadratic function of $k=\lfloor r/2R \rfloor$, as desired.
However, shadows are not in general disjoint.

If two distinct shadows intersect, then one of them is a subset
of the other, say $\Sh(y_i) \subset \Sh(y_j)$.  Subdividing $B_j$
along the geodesics $\beta_i$ and~$\beta'_i$ gives a decomposition
of $B_j$ into three brooms: $B_i$ and two others $C_j$ and~$C'_j$
which branch at~$y_j$
and have branching angles adding up to at least~$\theta$,
as shown in Figure~\ref{fig:shadows}.
Applying Lemma~\ref{lem:BroomDivergence} to these three brooms separately,
we see that
\[
  \ell \bigl( \Sh(y_j) \bigr) \ge R(k-i)\theta + \ell \bigl( \Sh(y_i) \bigr).
\]
Repeatedly subdividing brooms in this manner until they all have disjoint
interiors shows that, in fact,
\[
  \ell(\gamma) \ge R \theta k(k-1) \big/ 2,
\]
as desired.
\end{proof}
\begin{figure}[ht!]
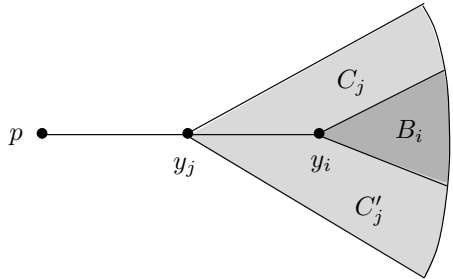

\drawshadows
\caption{Two shadows that intersect}
\label{fig:shadows}
\end{figure}

\begin{prop}[Ruffled Fellow Traveller Property]\label{prop:RuffledFTP}
Given positive constants $R$,~$\theta$, $\lambda$, and~$\epsilon$,
there is a constant~$L$ such that the following property holds.
Let $D$ be any piecewise Euclidean $\CAT(0)$ disc diagram
whose boundary is a concatenation $\alpha\gamma$,
where $\gamma$ is a geodesic and $\alpha$ is a
$(\lambda,\epsilon)$--quasigeodesic parametrized by arclength.
Suppose also that $(D,\gamma)$ is $(R,\theta)$--ruffled.
Then the Hausdorff distance between $\Image(\alpha)$ and
$\Image(\gamma)$ is less than~$L$.
\end{prop}

\begin{proof}
We follow the proof of \cite[Proposition~3.3]{ABC91} almost verbatim.
Let $Q$ be the quadratic function guaranteed by
Proposition~\ref{prop:QuadraticDivergence} corresponding to the ruffling
constants $R$~and~$\theta$.  Fix a disc diagram~$D$ as in the
statement of the proposition.
Let
\[
  r = \sup_{x\in \Image \gamma} \bigl\{ d(x,\Image \alpha) \bigr\},
\]
and choose a point $x \in \Image(\gamma)$ where this supremum is achieved.
Then
\[
  \Image(\alpha) \cap \ball{x}{r} = \emptyset.
\]
Let $y$~and~$z$ be points on $\Image(\gamma)$ at a distance~$r$
from $x$, and let $y'$ and~$z'$ be points on $\Image(\gamma)$ at
a distance~$2r$ from $x$ (or the endpoints of~$\gamma$ if these are
closer to~$x$ than~$2r$).  Choose points $u$~and~$v$ on $\Image(\alpha)$
such that $d(y',u) \le r$ and $d(z',v) \le r$, as in
Figure~\ref{fig:ruffledftp}.
Notice that
\[
  \bigl( [y',u] \cup [z',v] \bigr) \cap \ball{x}{r} = \emptyset.
\]
Following a path by way of $y'$,~$x$, and~$z'$, we see that
$d(u,v) \le 6r$.  However, since $\alpha$ is a
$(\lambda,\epsilon)$--quasigeodesic parametrized by arclength,
we have that the length of~$\alpha$ from $u$ to~$v$ is at most
$6\lambda r + \epsilon$.
Hence there is a path~$\beta$ of length at most
$4r + 6\lambda r + \epsilon$ from $y$ to~$z$ which stays outside
$\ball{x}{r}$.  Furthermore, $\beta$ together with $[y,z]$
bounds a subdiagram~$D'$ of~$D$ such that $\bigl( D',[y,z] \bigr)$
is $(R,\theta)$--ruffled by Lemma~\ref{lem:SubdiagramRuffled}.
But Proposition~\ref{prop:QuadraticDivergence} says that $\beta$~has
length at least $Q(r)$.  Therefore $r$ is bounded above by
some constant $L_0(R,\theta,\lambda,\epsilon)$,
and we see that $\Image(\gamma) \subseteq \bignbd{\Image(\alpha)}{L_0}$.
\begin{figure}[ht!]
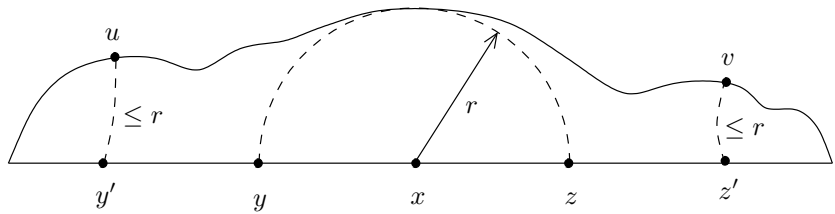

\drawruffledftp
\caption{The Ruffled Fellow Traveller Property}
\label{fig:ruffledftp}
\end{figure}

Now suppose that $\Image(\alpha) \nsubseteq \bignbd{\Image(\gamma)}{L_0}$.
Then each component of $\Image(\alpha) \setminus \bignbd{\Image(\gamma)}{L_0}$
is a path $\alpha'$ with endpoints $u$ and~$v$ at a distance $L_0$
from points $y$~and~$z$ on $\Image(\gamma)$.
Suppose $\alpha$ is a map $[0,a] \to D$, and
$\alpha' = \alpha \big| [t_0,t_1]$.
Then each point of $\Image(\gamma)$ is within~$L_0$ of some point of
$\alpha \bigl( [0,t_0] \bigr) \cup \alpha \bigl( [t_1,a] \bigr)$ by the first
part of the proof. So there must be some point~$x$ on $\Image(\gamma)$ which is
within $L_0$ of some point $u_0$ on $\alpha \bigl( [0,t_0] \bigr)$ and also
within $L_0$ of some point $u_1$ on $\alpha \bigl( [t_1,a] \bigr)$.
Thus $d(u_0,u_1) \le 2L_0$, and we see that the length of~$\alpha$
from $u_0$ to~$u_1$ is at most $2L_0 \lambda + \epsilon$,
which bounds the length of~$\alpha'$.
It now follows that every point on~$\alpha'$ is at most a distance~$L$
from $\Image(\gamma)$, where $L = L_0 + L_0 \lambda + \epsilon /2$.
\end{proof}

%%%%%%%%%%%%%%%%%%%%%%%%%%%%%%%%%%%%%%%%%%%%%%%%%%%%%%%%%%%%%%%%%%%%%%%%%%
\section{Convex hulls involving preflats}
\label{sec:ConvexHull}
%%%%%%%%%%%%%%%%%%%%%%%%%%%%%%%%%%%%%%%%%%%%%%%%%%%%%%%%%%%%%%%%%%%%%%%%%%

Recall that we need to prove Theorem~\ref{thm:RelativeFTP}
in order to complete the proof of Theorem~\ref{thm:2dEquivalent}.
Our strategy for proving Theorem~\ref{thm:RelativeFTP}
is to first consider the special case of a geodesic and
quasigeodesic with common endpoints.

This special case consists of showing that a quasigeodesic
stays close to the union of a geodesic and the flats that come near the
geodesic.  In order to prove this special case, it will be useful
to understand in detail the structure of the convex hull of the union of
two flats and also the convex hull of the union of a point and a flat.
In the present diagrammatic setting, it suffices to examine
the convex hulls of the corresponding objects inside a
reduced disc diagram.

In Proposition~\ref{prop:PointPlaneHull}
we give a detailed examination of the convex hull
of the union of a preflat~$P$ and a point~$x$ outside the preflat.
We determine that this convex hull lies in a $\delta$--neighborhood
of the union of~$P$ and the shortest geodesic connecting it
to~$x$.
We also conclude that the portion of the convex hull close to the preflat
is surrounded by ruffles.

We then examine the convex hull of the union of two preflats.
In Proposition~\ref{prop:DisjointHull} we consider the case of two disjoint
preflats, while in Proposition~\ref{prop:IntersectingHull}
we consider the case of two intersecting preflats.
In each case, our conclusion is similar to the conclusion of
Proposition~\ref{prop:PointPlaneHull}.
We determine that the convex hull lies in a $\delta$--neighborhood
of the union of the two preflats and a shortest geodesic connecting them.
We also conclude that the portion of the convex hull close to
either preflat is surrounded by ruffles.

\begin{figure}[ht!]
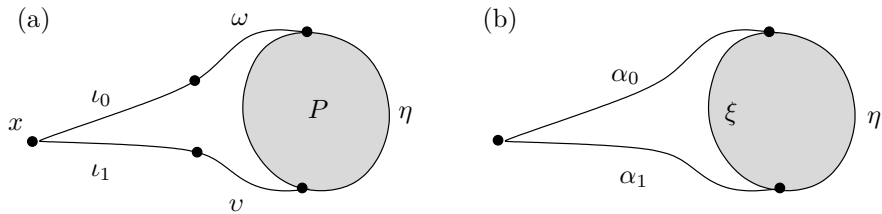

\drawPointPlaneHull
\caption[The convex hull of a point and a preflat.]%
{(a) The convex hull of a point and a preflat\qua
(b) The paths $\alpha_0$,~$\alpha_1$ and~$\xi$ bound a triangular diagram
which is $\delta$--thin.}
\label{fig:PointPlaneHull}
\end{figure}

\begin{prop}\label{prop:PointPlaneHull}
Given a proper, cocompact piecewise Euclidean $\CAT(0)$ $2$--complex $X$
with the Isolated Flats Property,
there are positive constants $\delta_1(X)$, and $\theta_1(X)$
such that the following property holds
for any positive $\delta\ge\delta_1$ and $\theta\le\theta_1$.

Let $\phi\co D \to \hat{X}$ be a reduced disc diagram,
where $\hat{D}$ is a subdivision of\/~$X$.
Let $P$ be a preflat in~$D$, let $x$ be a point in
$D \setminus P$, and let $H$ be the convex hull of $P \cup \{x\}$.

Then $x \in \boundary H$, and the boundary path~$\gamma$
of~$H$ beginning and ending at the point~$x$ is a local geodesic
in $D' = D \setminus \interior{P}$ (when considered as a
map with domain an interval, rather than as a map $S^1 \to \boundary H$).
This local goedesic
is a concatenation $\iota_0 \omega \eta \upsilon \iota_1$
such that, as illustrated in Figure~\ref{fig:PointPlaneHull}(a),
\begin{enumerate}
  \item $\eta$ is a subpath of the boundary cycle of\/~$P$,
  \item $\Image(\upsilon)$ and $\Image(\omega)$
        each lie in $\nbd{P}{\delta}$,
  \item $(D - \interior H, \omega\eta\upsilon)$
        is ($\delta, \theta)$--ruffled, and
  \item the Hausdorff distance between $\Image(\iota_0)$ and
        $\Image(\iota_1)$ is at most~$\delta$.
\end{enumerate}
\end{prop}

\begin{proof}
We first check that $x \in \boundary H$.  Notice that, away from~$x$,
the boundary cycle~$C$ of~$H$ is a local geodesic in~$D'$.
Recall that preflats are convex by Lemma~\ref{lem:BoundaryIsGeodesic}.
If $x \notin \boundary H$ then $C$ and $\boundary P$ are homotopic
local geodesics in~$D'$, which are therefore identical.
But this is absurd since $C$ encloses~$x$ and $\boundary P$ does not.
So $x \in \boundary H$ as desired.

Let $\gamma$ denote the boundary cycle of~$H$ considered as a
path which starts and ends at~$x$ (as opposed to~$C$ which is a
map $S^1 \to H$).  We have seen that $\gamma$ is a local geodesic.
We now check that $\Image(\gamma)$ intersects $\boundary P$
in a single connected arc.
If $\Image(\gamma) \cap \boundary P = \emptyset$, then $\gamma$ is a
local geodesic in~$D$.  But $D$ is simply connected, so $\gamma$ must
be a geodesic, which contradicts the fact that $\gamma$
is a closed nonconstant loop.
So $\Image(\gamma) \cap \boundary P$ consists of a finite (positive) number
of arcs.  Suppose there were more than one such arc.
Then $\gamma$ would contain a subpath which does not involve~$x$
and which intersects $\boundary P$ only at the endpoints of the subpath.
Such a subpath would be a geodesic in~$D$, contradicting the
convexity of~$P$.

We have now determined that $\gamma$ is a concatenation
$\alpha_0 \eta \alpha_1$ such that $\eta$ is the maximal subpath of~$\gamma$
lying inside $\boundary P$, and such that $\alpha_i$ is a
geodesic in~$D$ connecting $x$ and~$P$.  Let $\xi$ be the path
so that the boundary cycle of~$P$ is a concatenation $\eta \xi$,
as shown in Figure~\ref{fig:PointPlaneHull}(b).
By Proposition~\ref{prop:AroundFlat}, there are global constants
$R(X)$ and $\theta(X)$ so that $(D', \xi)$ is $(R, \theta)$--ruffled.
The three $D'$--geodesics $\alpha_0$, $\xi$, and~$\alpha_1$
bound a triangular diagram~$\Delta$,
which is a subdiagram of a subdivision of~$D'$.
By Lemma~\ref{lem:SubdiagramRuffled}, we have that
$(\Delta, \xi)$ is also $(R, \theta)$--ruffled.
Furthermore, we see by Proposition~\ref{prop:ThinTriangles} that
$\Delta$ is $\delta$--thin for some constant $\delta = \delta(R,\theta)$.
So there are points $y_i \in \Image(\alpha_i)$
with $d(y_0,y_1) < \delta$.
Therefore, $\alpha_0$ is a concatenation $\iota_0 \omega$
and $\alpha_1$ is a concatenation $\upsilon \iota_1$
such that (2) and (4) are satisfied.

It remains to verify (3).  We already know that
$(D \setminus \interior{H}, \eta)$ is $(R, \theta)$--ruffled.
We now show that $(D \setminus \interior{H})$
is $(R', \theta')$--ruffled for some global constants $R'$ and~$\theta'$.
Let $a$ be the common endpoint of $\upsilon$ and~$\xi$,
let $b$ be the common endpoint of $\upsilon$ and~$\iota_1$,
and let $c$ be a point on $\Image(\xi)$ within a distance~$\delta$
of~$b$.  The points $a$,~$b$, and~$c$ are the vertices of a triangular
subdiagram of~$\Delta$
which has $\upsilon$ and a portion of~$\xi$ as two sides
and with a third side $[b,c]$ of length less than~$\delta$.
Lemma~\ref{lem:RuffledTriangle} now shows that
$(D \setminus \interior{H}, \upsilon)$ is $(R',\theta/2)$--ruffled
for some constant $R'$ depending only on $R$,~$\theta$, and~$\delta$.
By a nearly identical argument, we see that
$(D \setminus \interior{H},\omega)$ is also $(R', \theta/2)$--ruffled.
Since $D \setminus \interior{H}$ is now known to be ruffled along
each of $\omega$,~$\eta$, and~$\upsilon$,
Lemma~\ref{lem:ConcatenatedRuffles} shows that
$(D \setminus \interior{H}, \omega\eta\upsilon)$
is $(3R', \theta/2)$--ruffled.
Setting $\theta_1 = \theta/2$ and
$\delta_1 = \max \{\delta, 3R'\}$ gives (4) for any positive
$\delta \ge \delta_1$ and $\theta \le \theta_1$, as desired.
\end{proof}

\begin{figure}[ht!]
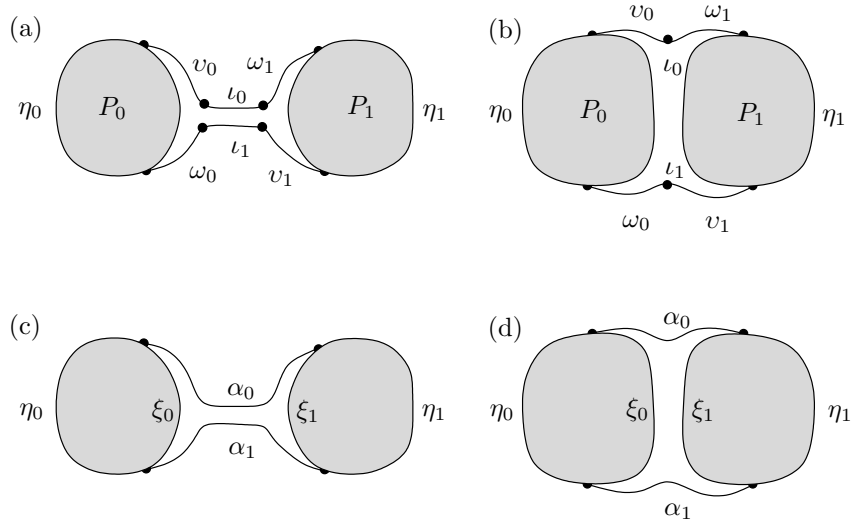

\drawdisjoint
\caption[The convex hull of the union of two disjoint preflats.]%
{The convex hull of two disjoint preflats $P_0$ and~$P_1$.
In (a), the Hausdorff distance between $\iota_0$ and $\iota_1$ is
at most~$\delta$.
In (b), each $\iota_i$ has image a single point.
In (c) and (d) the quadrilateral bounded by $\xi_0$,~$\alpha_0$,
$\xi_1$, and~$\alpha_1$ is $\delta$--thin.
The two pictures indicate the two
possible shapes of a thin quadrilateral.}
\label{fig:disjoint}
\end{figure}

\begin{prop}\label{prop:DisjointHull}
Given a proper, cocompact piecewise Euclidean $\CAT(0)$ $2$--complex $X$
with the Isolated Flats Property,
there are positive constants $\delta_2(X)$, and $\theta_2(X)$
such that the following property holds
for any positive $\delta\ge\delta_2$ and $\theta\le\theta_2$.

Let $\phi\co D\to \hat{X}$ be a reduced disc diagram,
where $\hat{X}$ is a subdivision of\/~$X$.  Let
$\preflat_0$ and~$\preflat_1$ be disjoint preflats in~$D$.
Let $H$ be the convex hull of\/ $\preflat_0 \cup \preflat_1$.

Then the boundary cycle $C\to D$ of\/~$H$ is a local geodesic
in
$ D' = D \setminus \bigl( \interior\preflat_0 \cup \interior\preflat_1 \bigr)$
which is a concatenation
\[
   \omega_0 \eta_0 \upsilon_0 \iota_0 \omega_1 \eta_1 \upsilon_1 \iota_1
\]
such that, as illustrated in Figures~\ref{fig:disjoint}(a) and~(b),
\begin{enumerate}
  \item\label{item:DHa} $\eta_i$ is a subpath of the boundary cycle
     of\/~$\preflat_i$,
  \item\label{item:DHb} $\Image(\upsilon_i)$ and $\Image(\omega_i)$
        each lie in $\nbd{\preflat_i}{\delta}$,
  \item\label{item:DHd} $(D\setminus\interior{H}, \omega_i \eta_i \upsilon_i)$
        is $(\delta,\theta)$--ruffled, and
  \item\label{item:DHc} either
     \begin{enumerate}
       \item[\rm(a)]\label{item:DH:wide}the Hausdorff distance between
          $\Image(\iota_0)$ and $\Image(\iota_1)$ is at most~$\delta$, or
       \item[\rm(b)]\label{item:DH:tall}
          $P_0$ and~$P_1$ correspond to distinct flats
          in~$X$, and each~$\iota_i$ has image a single point~$x_i$.
     \end{enumerate}
\end{enumerate}
\end{prop}

\begin{proof}
Since the present proof is extremely similar to the proof of
Proposition~\ref{prop:PointPlaneHull}, we omit the redundant details.
The only new idea occurs during the verification of~(4).
By an argument similar to the one in the previous proof, we obtain
a $\delta$--thin quadrilateral~$Q$
with sides $\xi_0 \alpha_0 \xi_1 \alpha_1$ as shown in
Figures~\ref{fig:disjoint}(c) and~(d).
Recall that $\delta$--thin quadrilaterals have two possible shapes:
either $\alpha_0$ and $\alpha_1$ come close together, or
$\xi_0$ and $\xi_1$ come close together.
This dichotomy accounts for the two distinct shapes in
Figure~\ref{fig:disjoint} and described in the statement of
condition~(4).

We need to verify the assertion in~(4)
that if $P_0$ and $P_1$ correspond to the same flat in~$X$
then the Hausdorff distance between $\Image(\iota_0)$
and $\Image(\iota_1)$ is at most~$\delta$.
It suffices to consider the following situation.
Suppose $P_0$ and~$P_1$ correspond to the same flat in~$X$,
and there exist points $x_i \in \Image(\alpha_i)$
so that for each $i,j$ the distances $d \bigl( x_i, \Image(\xi_j) \bigr)$
are less than~$\delta$.
In order to complete the proof of~(4), we need to bound the
distance $d(x_0,x_1)$.  We can then set the constant $\delta_2$
to be at least as large as this distance.

\begin{figure}[ht!]
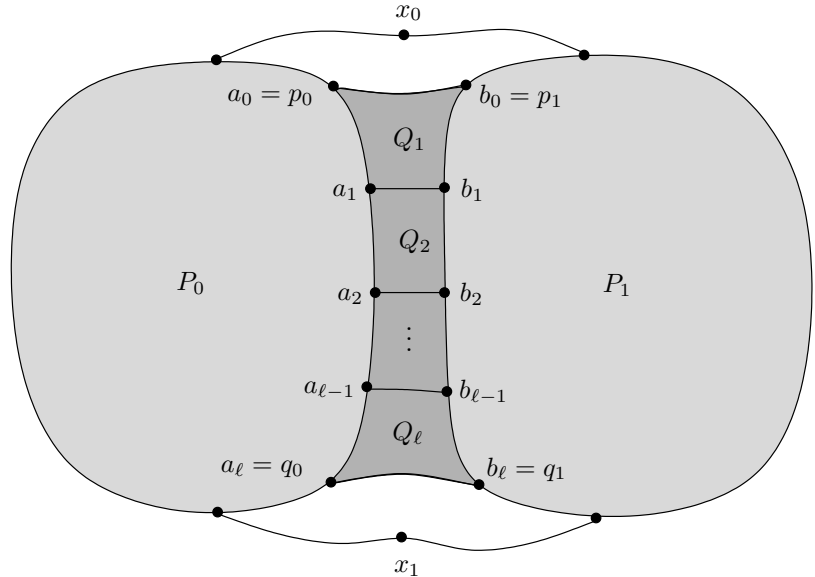

\drawBetweenPreflats
\caption[A geodesic quadrilateral between two preflats.]%
{The preflats $P_0$ and~$P_1$ map to the same flat in~$X$.
The darkened diagram is a geodesic quadrilateral~$Q$ with corners
$p_0$,~$p_1$, $q_0$, and~$q_1$.
This quadrilateral has been subdivided into $\ell$ smaller
quadrilaterals~$Q_k$.}
\label{fig:BetweenPreflats}
\end{figure}

For the sake of notation, suppose $x_0$ is within a distance~$\delta$
of points $p_j \in \Image(\xi_j)$ and $x_1$ is within a distance~$\delta$
of points $q_j \in \Image(\xi_j)$, as shown in Figure~\ref{fig:BetweenPreflats}.
In order to bound $d(x_0, x_1)$, it suffices to bound the quantity
\[
   \ell = \left\lfloor \frac{ \min \bigl\{ d(p_0, q_0), d(p_1,q_1) \bigr\} }
                            {4\delta} \right\rfloor.
\]
If one of the distances $d(p_j,q_j)$ is less than $4\delta$, then $\ell=0$
and we are done.  Otherwise,
let $a_0, \dots, a_\ell$ be equally spaced points on $\xi_0$
so that $a_0 = p_0$ and $a_\ell = q_0$.
Similarly, let $b_0, \dots, b_\ell$ be equally spaced points on~$\xi_1$
so that $b_0 = p_1$ and $b_\ell = q_1$.
By our choice of~$\ell$, for $0 \le k \le \ell - 1$ we have
\begin{equation}
   d(a_k, a_{k+1}) \ge 4 \delta \quad \text{and} \quad
   d(b_k, b_{k+1}) \ge 4 \delta. \tag{$*$}
\end{equation}
Furthermore $a_k$,~$a_{k+1}$, $b_k$, and $b_{k+1}$ are the vertices of a
quadrilateral diagram~$Q_k$ as shown in Figure~\ref{fig:BetweenPreflats}.

We first rule out the degenerate case that $Q_k$ is not a single line segment.
Note that by Theorem~\ref{thm:CAT0convexity}, the distances $d(a_k,b_k)$ and
$d(a_{k+1},b_{k+1})$ are each less than $2\delta$.
If the four corners of~$Q_k$ were colinear, then by~($*$)
the
geodesics $[a_k,a_{k+1}]$ and $[b_k,b_{k+1}]$ would intersect,
contradicting the fact that $\xi_0$ and $\xi_1$ are disjoint paths.

Consequently, each~$Q_k$ contains some vertex with negative
curvature.  For if not, then by Corollary~\ref{cor:diameterpi}
the diagram~$Q_k$ would map isometrically to~$X$ under~$\phi$.
But that would mean that $Q_k$ maps into the same flat that
the preflats $P_0$ and~$P_1$ map to, which is absurd.

By Theorem~\ref{thm:CGB}, the quadrilateral~$Q$ with vertices
$p_0$,~$p_1$, $q_0$, and~$q_1$ has total curvature $2\pi$.
The only vertices in~$Q$ that could possibly have positive curvature
are these four corners.  At each of these, the curvature is at most~$\pi$.
So the positive curvature in~$Q$ totals at most $4\pi$.
As a result, the negative curvature in~$Q$ has total magnitude at most
$2\pi$.  So the curvature at each negatively curved vertex of~$Q$
has magnitude at most $2\pi$.

We now show that by our choice of~$\ell$,
each vertex with negative curvature in~$Q$ intersects at most two of the
quadrilaterals~$Q_k$.
Suppose $Q_{k-1}$, $Q_k$, and $Q_{k+1}$ share a common point~$z$.
Then $z$ is a cut point of~$Q_k$ at which the sides $[a_k,b_k]$ and
$[a_{k+1},b_{k+1}]$ meet.  But by the triangle inequality, these sides cannot
meet, since their lengths are each less than $2\delta$ and we also have~($*$).

Except along the paths $[a_0,b_0]$ and $[a_\ell,b_\ell]$, every point of
negative curvature in~$Q$ has curvature with magnitude at least $\psi$, for
some constant~$\psi$ depending only on~$X$.
So there are at most $2\pi/\psi$ points of negative curvature
within the quadrilaterals $Q_2, \dots, Q_{\ell-1}$, which shows that
$\ell \le (4\pi / \psi) + 2$.
\end{proof}

%The following proposition is proved in essentially the same way as
%Proposition~\ref{prop:PointPlaneHull}, without the additional complications
%that arose in the proof of Proposition~\ref{prop:DisjointHull}.

\begin{figure}[ht!]
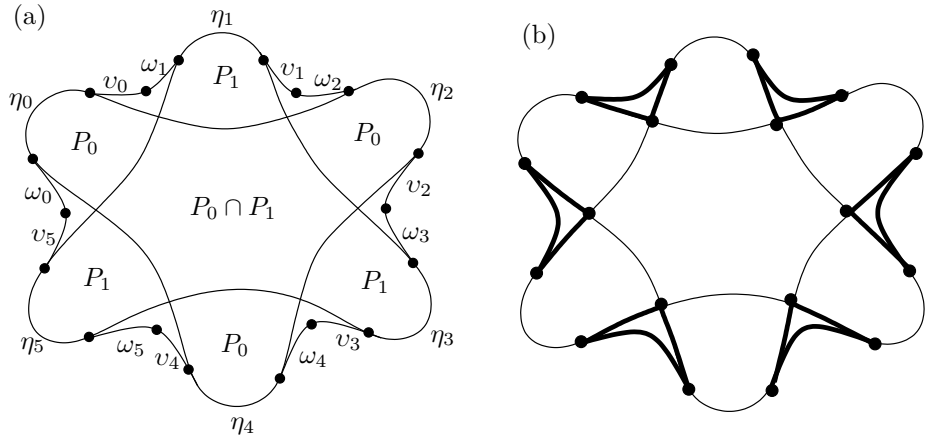

\drawintersectingone
\caption[The convex hull of the union of two intersecting preflats.]%
{(a)~The convex hull of the union of two intersecting preflats\qua
(b)~The darkened triangles are $\delta$--thin.}
\label{fig:intersectingone}
\end{figure}

\begin{prop}\label{prop:IntersectingHull}
Given a proper, cocompact piecewise Euclidean $\CAT(0)$ $2$--complex $X$
with the Isolated Flats Property, there are positive constants
$\delta_3(X)$ and $\theta_3(X)$
such that the following property holds
for any positive $\delta\ge\delta_3$ and $\theta\le \theta_3$.

Let $\phi\co D \to \hat{X}$ be a reduced disc diagram,
where $\hat{X}$ is a subdivision of\/~$X$. Let $\preflat_0$
and~$\preflat_1$ be preflats in~$D$ with a nonempty intersection.
Let $H$ be the convex hull of\/
$\preflat_0 \cup \preflat_1$.

Then the boundary cycle $C\to D$ of\/~$H$ is a local geodesic in
$D' = D \setminus \bigl( \interior\preflat_0 \cup \interior\preflat_1 \bigr)$,
which can be expressed as a concatenation
\[
  \omega_0 \eta_0 \upsilon_0 \omega_1 \eta_1 \upsilon_1 \dotsm
  \omega_{2m-1} \eta_{2m-1} \upsilon_{2m-1}
\]
such that, as illustrated in Figure~\ref{fig:intersectingone}(a),
\begin{enumerate}
   \item\label{item:IHa} $\eta_i$ is a subpath of the boundary cycle
      of\/~$\preflat_{(i \bmod 2)}$,
   \item\label{item:IHb} $\Image(\upsilon_i)$ and $\Image(\omega_i)$
      each lie in $\nbd{\preflat_{(i \bmod 2)}}{\delta}$, and
   \item\label{item:IHc}
      $(D \setminus \interior H, \omega_i \eta_i \upsilon_i)$ is
      $(\delta,\theta)$--ruffled.
\end{enumerate}
\end{prop}

\begin{proof}
The key idea is to observe that the triangles shown in
Figure~\ref{fig:intersectingone}(b) are $\delta$--thin.
Using this observation,
the proof is essentially the same as the proof of
Proposition~\ref{prop:PointPlaneHull} without the additional complications
that arose in the proof of Proposition~\ref{prop:IntersectingHull}.
Since the details are nearly identical, we omit the verification.
\end{proof}

It is important to note that condition~(\ref{item:IHc}) does not in general
imply that $(D \setminus \interior H, \boundary H)$ is $(R,\psi)$--ruffled
for constants $R$ and $\psi$ depending only on the $2$--complex~$X$.
Although the boundary cycle of~$H$ is a concatenation of the paths
$\omega_i\eta_i\upsilon_i$ for $0 \le i \le 2m-1$,
Lemma~\ref{lem:ConcatenatedRuffles} can be applied only if one has either a
lower bound on the lengths of the concatenated paths or an upper bound on the
number of these paths.  In general we do not have control over either of these
quantities.

%%%%%%%%%%%%%%%%%%%%%%%%%%%%%%%%%%%%%%%%%%%%%%%%%%%%%%%%%%%%%%%%%%%%%%%%%%%%
\section{The Flat Closure Lemma}
\label{sec:FlatClosure}
%%%%%%%%%%%%%%%%%%%%%%%%%%%%%%%%%%%%%%%%%%%%%%%%%%%%%%%%%%%%%%%%%%%%%%%%%%%%

Recall that we will prove Theorem~\ref{thm:RelativeFTP}
by first considering the special case of a geodesic and quasigeodesic
with common endpoints, and then deriving the general case of
two quasigeodesics from this special case.

This section is devoted to a proof of the following proposition,
which will be used in the next section in the proof of
Theorem~\ref{thm:RelativeFTP}.
The reader should imagine that the diagram~$D$ in the statement of this
proposition is a reduced disc diagram such that $\gamma$ is a geodesic and
$\alpha$ a quasigeodesic,
although for the purposes of this section
we do not need the path $\alpha$ to have any special properties.

The name Flat Closure refers to the fact that we construct a subdiagram
of~$D$ containing all preflats that come within a certain distance
of the geodesic~$\alpha$.

\begin{prop}[Flat Closure]\label{prop:FlatClosure}
Let $X$ be a proper, cocompact
$\CAT(0)$ $2$--complex with the Isolated Flats Property.
There are positive constants $L$, $K$ and~$\theta$ such that the following
property holds.

Let $\phi\co D \to \hat X$ be any reduced disc diagram
where $\hat X$ is a subdivision of~$X$ such that the
boundary path of~$D$ is a concatenation $\alpha\gamma$ where $\phi\of\gamma$
is a geodesic in~$X$.
Then there is a path~$\beta$ in~$D$ with the same endpoints
as $\alpha$ and~$\gamma$ satisfying the following three conditions:
\begin{enumerate}
\renewcommand{\theenumi}{\Roman{enumi}}
   \item \label{item:FC:main3}
         $\beta$ is a geodesic in~$D_\alpha$, where $D_\alpha$
         is the subdiagram of~$D$ bounded by $\beta$ and~$\alpha$,
   \item \label{item:FC:main1}
         $\phi\of\gamma$ and $\phi\of\beta$ \ $L$--fellow
         travel relative to flats, and
   \item \label{item:FC:main2}
         $(D_\alpha, \beta)$ is $(K,\theta)$--ruffled.
\end{enumerate}
\end{prop}

\begin{proof}
We begin by carefully choosing several constants
which will be necessary
for the construction of~$\beta$.
First let $\delta=\delta(X)$ be the largest of the three constants
$\delta_1(X)$,~$\delta_2(X)$, and~$\delta_3(X)$
guaranteed by Propositions \ref{prop:PointPlaneHull},
\ref{prop:DisjointHull}, and~\ref{prop:IntersectingHull}
respectively.
By the Isolated Flats Property,
there is a constant
$M(\delta,X) \ge 6\delta$
such that for any two flat planes $E_1,E_2$ in~$X$, the intersection
$\nbd{E_1}{\delta} \cap \nbd{E_2}{\delta}$
has diameter less than~$M$.
Using $M_0=2M$, choose constants $R(M_0,X)$ and $\theta_0(M_0,X)$
satisfying the conclusion of Proposition~\ref{prop:AlongGeodesic}.
Let $\theta = \min\{\theta_0, \theta_1, \theta_2, \theta_3\}$, where
$\theta_1$,~$\theta_2$, and~$\theta_3$
are the constants given by Propositions \ref{prop:PointPlaneHull},
\ref{prop:DisjointHull} and~\ref{prop:IntersectingHull}.

The next step is the construction of~$\beta$, which proceeds as follows.
Let $\Set{P}$ be the set of all preflats in the diagram~$D$ which intersect
$\nbd{\Image{\gamma}}{R}$, where $R$ is some positive constant
to be determined in the course of the proof.
Form a subdiagram $D'$ from~$D$ by removing the interior of
$\bigcup_{P \in \Set{P}} P$.  In that subdiagram,
let $\beta$ be the local geodesic homotopic to~$\alpha$ (rel.
endpoints).  Notice that $\beta$ cuts $D$ into two subdiagrams
$D_\alpha$ and $D_\gamma$ which intersect along~$\beta$
as shown in Figure~\ref{fig:FlatClosureOne}.
Notice that $D_\gamma$ is the convex hull of
$\Image(\gamma) \cup \bigl(\bigcup_{P \in \Set{P}} P\bigr)$.
Condition~(\ref{item:FC:main3}) follows easily from the fact that
$D_\gamma$ is convex.
\begin{figure}[ht!]
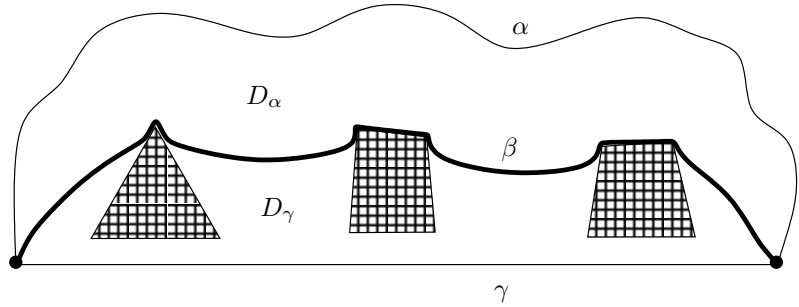

\drawFlatClosureOne
\caption{The path~$\beta$ cuts $D$ into two subdiagrams $D_\alpha$
and~$D_\gamma$ that meet along~$\beta$.}
\label{fig:FlatClosureOne}
\end{figure}

\begin{claim}\label{claim:FC1}
$\beta$ is a concatenation
\[
   \tau_0 \pi_1 \rho_1 \sigma_1 \tau_1 \pi_2 \rho_2 \sigma_2 \tau_2
   \dotsm \pi_n \rho_n \sigma_n \tau_n,
\]
as shown in Figure~\ref{fig:flatclosure},
such that for some sequence of flats $(E_1, \dots, E_k)$
and some constant $\delta'(X)$,
the following properties hold for all $i$:
\begin{description}
   \item[\rm A($i$):]\label{item:FCa} $\Image (\phi\of\rho_i) \subseteq
      \nbd{E_i}{\delta}$, and
      $(D_\alpha, \rho_i)$ is $(\delta,\theta)$--ruffled.
   \item[\rm B($i$):]\label{item:FCb} $\Image (\phi\of\sigma_i) \subseteq
      \nbd{E_i}{\delta}$, and
      $(D_\alpha, \sigma_i)$ is $(\delta,\theta)$--ruffled.
   \item[\rm C($i$):]\label{item:FCc} $\Image (\phi\of\tau_i) \subseteq
      \bignbd{\Image(\phi\of\gamma)}{R+M}$, and
      $(D_\alpha, \tau_i)$ is $(\delta',\theta / 4)$--ruffled.
   \item[\rm D($i$):]\label{item:FCd} $\Image (\phi\of\pi_{i+1}) \subseteq
      \nbd{E_{i+1}}{\delta}$, and
      $(D_\alpha, \pi_{i+1})$ is $(\delta,\theta)$--ruffled.
\end{description}
\end{claim}

\begin{figure}[ht!]
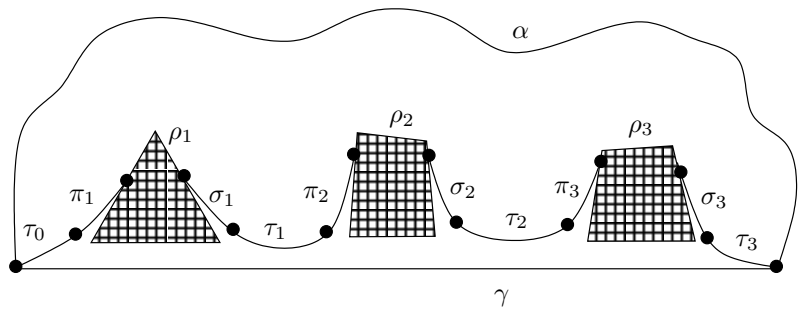

\drawflatclosure
\caption{The path~$\beta$ stays close to the union of the preflats
in~$\Set{P}$ and $\Image(\gamma)$.}
\label{fig:flatclosure}
\end{figure}

\begin{proof}[Proof of Claim~\ref{claim:FC1}]
For each preflat~$P \in \Set{P}$, the preimage $\beta^{-1}(P)$ consists of a
finite number of closed intervals.  Consider the collection
of all such intervals for all choices of~$P \in \Set{P}$.
In general, these intervals intersect each other in a complicated fashion.
To simplify matters, let $\bigset{[a_i,b_i]}{1\le i\le n}$ be a minimal
subcollection of these intervals with the same union as the original collection.
Then for each~$i$, the path $\beta \bigl( [a_i,b_i] \bigr)$ lies in
some~$\preflat_i$.
By the minimality of the collection, no interval is contained within
another,
and there are no nonempty triple intersections among the intervals.
We may assume that the intervals are ordered so that
\[
  a_1 < a_2 < \dots < a_n
\]
and
\[
  b_1 < b_2 < \dots < b_n.
\]
Our strategy is to verify properties
A($i$),~B($i$),
C($i$), and~D($i$)
by examining in detail the path that $\beta$ follows between two consecutive
intervals $[a_i,b_i]$ and $[a_{i+1}, b_{i+1}]$.
We will have two cases depending on whether these intervals intersect or are
disjoint.
The second case breaks into two cases depending on whether the associated
preflats $P_i$ and~$P_{i+1}$ intersect or are disjoint.
The disjoint preflats case further subdivides into two cases depending on
whether these preflats satisfy Condition (4a)
or~(4b) of Proposition~\ref{prop:DisjointHull}.

We will also examine in detail what happens to~$\beta$ between the
initial point $\beta(0)$ and the first interval $[a_1,b_1]$, as well as an
identical case at the terminal end of $\beta$.
We also consider the degenerate case when the set of intervals $[a_i,b_i]$ is
empty.  This degenerate case arises precisely when $\Set{P} = \emptyset$.

To summarize, we have the following situations to examine:
\begin{description}
   \item[Case~0] $\bigl\{ [a_i,b_i] \bigr\} = \emptyset$.
   \item[Case~1] Between $\beta(0)$ and $[a_1,b_1]$ (and similarly
                  between $[a_n,b_n]$ and the terminal point of~$\beta$).
   \item[Case~2] Between two consecutive intervals such that
                 $P_i \cap P_{i+1} = \emptyset$ and these preflats satisfy
                 Condition~(4a) of
                 Proposition~\ref{prop:DisjointHull}.
   \item[Case~3] Between two consecutive intervals such that
                 $P_i \cap P_{i+1} = \emptyset$ and these preflats satisfy
                 Condition~(4b) of
                 Proposition~\ref{prop:DisjointHull}.
   \item[Case~4] Between two disjoint consecutive intervals such that
                 $P_i \cap P_{i+1} \ne \emptyset$.
   \item[Case~5] Between two intersecting intervals such that
                 $P_i \cap P_{i+1} \ne \emptyset$.
\end{description}

Our examination of these cases will establish properties
A($i$), B($i$),
C($i$), and D($i$)
for all values of~$i$.

\textbf{Case 0}\qua 
We first consider the degenerate case in which the collection
of intervals $\bigl\{ [a_i,b_i] \bigr\}$ is empty. In other words, the image
of~$\beta$ does not intersect any~$\preflat_k$.
Then $\beta$ is a local geodesic in the simply connected diagram~$D$.
It follows that $\beta$ is a global geodesic which, therefore,
coincides with~$\gamma$.
Since $\beta$~and~$\alpha$ are homotopic in~$D'$,
it follows that $D$~and~$D'$ are equal.  In other words,
there are no preflats intersecting $\bignbd{\Image(\gamma)}{R}$.
Setting $\tau_0=\beta=\gamma$, we see that
C($0$) follows
from Proposition~\ref{prop:AlongGeodesic}, using $\delta'=8R$.
Since $n=0$, all other properties in the statement of Claim~\ref{claim:FC1}
are vacuous.

\textbf{Case 1}\qua
Let $H$ denote the convex hull of $\beta(0) \cup \preflat_1$.
By Proposition~\ref{prop:PointPlaneHull},
we see that the boundary cycle of~$H$ is a local geodesic in
$D\setminus\interior\preflat_1$ except at the point $\beta(0)$.
Furthermore, by our choice of $\delta$ and~$\theta$
this boundary cycle is a concatenation
$\iota_0 \omega \eta \upsilon \iota_1$ as shown in
Figure~\ref{fig:FCstagezero}(a), satisfying
\begin{enumerate}
   \item \label{item:FC:0a}
         $\eta$ is a subpath of the boundary cycle of~$\preflat_1$,
   \item \label{item:FC:0b}
         $\Image(\upsilon)$ and $\Image(\omega)$ each lie in
         $\nbd{\preflat_1}{\delta}$,
   \item \label{item:FC:0d}
         $(D\setminus\interior H, \omega\eta\upsilon)$ is
         $(\delta,\theta)$--ruffled, and
   \item \label{item:FC:0c}
         the Hausdorff distance between $\Image(\iota_0)$ and $\Image(\iota_1)$
         is at most~$\delta$.
\end{enumerate}
If we set $\tau_0 = \iota_0$ and $\pi_1 = \omega$, then
D($0$) follows easily
from (\ref{item:FC:0b}) and~(\ref{item:FC:0d}),
since $D_\alpha$ is a subdiagram of $D\setminus\interior H$.
\begin{figure}[t!]
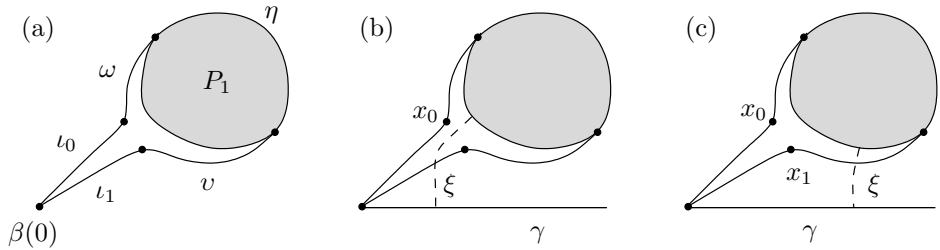

\drawFCstagezero
\caption[The convex hull of the union of the preflat~$P_1$ and
    the point $\beta(0)$.]{(a)~The convex hull of the preflat $P_1$ and
    the point $\beta(0)$\qua
    (b)~The path~$\xi$ intersects $\iota_1$.\qua
    (c)~The path~$\xi$ intersects $\eta \upsilon$.}
\label{fig:FCstagezero}
\end{figure}

It remains to verify C($0$).  Recall that $\preflat_1$
intersects the $R$--neighborhood of $\Image(\gamma)$.
Let $\xi$ be a geodesic from $\preflat_1$
to $\Image(\gamma)$ with length at most $R$.
Since $\preflat_1$ and $\Image(\gamma)$ each lie inside the convex
subdiagram $D_\gamma$, it follows that $\Image(\xi)$ also lies
inside $D_\gamma$.  Therefore $\Image(\xi)$ intersects the image of one of
the paths $\iota_1$, $\upsilon$, and~$\eta$.
If $\Image(\xi)$ intersects $\Image(\iota_1)$ as in
Figure~\ref{fig:FCstagezero}(b),
then $\Image(\xi)$ passes within
a distance~$\delta$ of the common endpoint~$x_0$ of $\iota_0$ and~$\omega$.
So $d(x_0, \Image \gamma) \le R+\delta$.
On the other hand, if $\Image(\xi)$ intersects the image of $\eta\upsilon$
as in Figure~\ref{fig:FCstagezero}(c),
then by the convexity of the metric, the common endpoint~$x_1$
of $\iota_1$ and~$\upsilon$ lies within a distance~$R$ of $\Image(\gamma)$.
So again, $d(x_0,\Image\gamma)  \le R + \delta$.
Thus, in either case, $\Image(\tau_0) \subseteq \bignbd{\Image(\gamma)}{R+M}$,
since $M \ge \delta$.

We now only have to find a constant $\delta'$ so that
$(D_\alpha, \tau_0)$ is $(\delta',\theta / 4)$--ruffled.
Let $x_2$ be the point on $\Image(\gamma)$ which is closest to~$x_0$,
as illustrated in Figure~\ref{fig:FCstagezerob},
and let $\Delta$ be the triangular subdiagram of~$D$ with corners $\beta(0)$,
$x_0$, and~$x_2$. Let $\hat{D} = D_\alpha \cup \Delta$,
and let $\hat{\gamma}$ be the subsegment
of~$\gamma$ connecting $\beta(0)$ and~$x_2$.
Recall that by construction any preflat in $D$ intersecting
$\bignbd{\Image(\gamma)}{R}$ lies inside the subdiagram $D_\gamma$.
The composition $\hat{D} \inclusion D \to X$ gives a natural notion
of a preflat in~$\hat{D}$ with the property that every preflat of~$\hat{D}$
lies inside a preflat of~$D$.  So any preflat in~$\hat{D}$ intersecting
$\bignbd{\Image(\hat{\gamma})}{R}$ lies inside the subdiagram~$\Delta$.
But $\Delta$ lies inside $\bignbd{\Image(\hat{\gamma})}{R+2M}$,
so we see that any preflat~$\preflat$ in~$\hat{D}$ intersecting
$\bignbd{\Image(\hat{\gamma})}{R}$ lies inside
$\bignbd{\Image(\hat{\gamma})}{R+2M}$.
Therefore, $\phi(\preflat)$ lies inside
$\bignbd{\Image (\phi \of \hat{\gamma})}{R+2M}$.
By our choice of $R$, $\theta$, and $M_0=2M$,
Proposition~\ref{prop:AlongGeodesic}
implies that $(\hat{D},\hat{\gamma})$ is $(8R,\theta)$--ruffled.
Since $[x_0,x_2]$ has length at most $R+\delta$,
it follows from Lemma~\ref{lem:RuffledTriangle} that $(D_\alpha,\tau_0)$
is $(\delta',\theta/2)$--ruffled for some constant $\delta'(R,\theta,\delta)$,
which completes the proof of C($0$).
\begin{figure}[ht!]
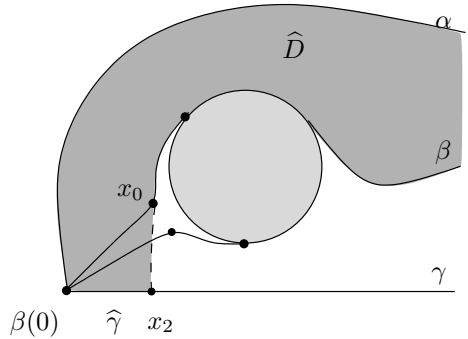

\drawFCstagezerob
\caption[The subdiagram~$\hat{D}$ consists of the union of
$D_\alpha$ and~$\Delta$.]%
{The darkened subdiagram~$\hat{D}$ consists of the union of $D_\alpha$
and the triangular diagram~$\Delta$ with vertices
$\beta(0)$, $x_0$, and~$x_2$.}
\label{fig:FCstagezerob}
\end{figure}

\textbf{Case~2}\qua Suppose $\preflat_{i} \cap \preflat_{i+1}= \emptyset$,
and the preflats in question satisfy Condition~(4a) of
Proposition~\ref{prop:DisjointHull}.
Let $H$ denote the convex hull of $\preflat_{i} \cup \preflat_{i+1}$.
Then the boundary cycle of~$H$
is a local geodesic in
$
  D \setminus \bigl( \interior \preflat_{i} \cup \interior
      \preflat_{i+1} \bigr)
$
freely homotopic to $\boundary D$.
By our choice of $\delta$ and~$\theta$, the cycle $\boundary H$ is a
concatenation $\omega_0 \eta_0 \upsilon_0 \iota_0 \omega_1 \eta_1 \upsilon_1
\iota_1$ satisfying the conclusion of Proposition~\ref{prop:DisjointHull}
including Condition~(4a),
as illustrated in Figure~\ref{fig:disjoint}(a).
Since $D_\alpha$ is a subdiagram of $D \setminus \interior H$,
Conditions B($i$) and~D($i$) follow easily if we set
$\sigma_i = \upsilon_0$, \  $\tau_i = \iota_0$
and $\pi_{i+1}=\omega_1$.

It remains to verify C($i$).
Recall that $\preflat_i$ and $\preflat_{i+1}$ each intersect
the $R$--neighborhood of $\Image(\gamma)$.
For each $j\in \{0,1\}$ let $\xi_j$ be a $D$--geodesic from $\preflat_{i+j}$
to $\Image(\gamma)$ with length at most~$R$, and let $C_j$ denote the
path $\omega_j \eta_j \upsilon_j$.  As in Case~1, the image of $\xi_j$
is inside the convex subdiagram $D_\gamma$, so it must intersect the image of
at least one of the paths $C_0$, $C_1$, or $\iota_1$.

Suppose $\Image(\xi_j)$ intersects $\Image(C_{j+1})$ for some $j\in \Z_2$,
as in Figure~\ref{fig:FCstageia}(a).
Since $H$ is convex, it follows that
$\Image(\iota_0) \subseteq \bignbd{\Image(\xi_j)}{\delta}$.
But $\xi_j$ has length at most~$R$,
so $\iota_0$ has length at most $R+2\delta$
and $\Image(\iota_0) \subseteq\ball{z_j}{R+\delta}$,
where $z_j$ is the point where
the images of $\gamma$ and~$\xi_j$ intersect.
Since $\tau_i = \iota_0$, it now follows immediately
that $(D_\gamma,\iota_0)$
is $(\delta + R/2,\theta / 4)$--ruffled.  Since $M \ge \delta$,
we now have C($i$),
using $\delta'=\delta + R/2$.
\begin{figure}[ht!]
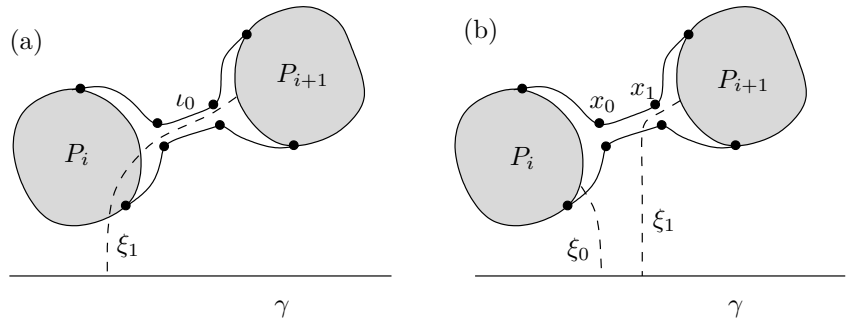

\drawFCstageia
\caption{(a) The path $\xi_1$ intersects $C_0$.\qua  (b) For each $j \in \Z_2$,
the path $\xi_j$ intersects either $C_j$ or $\iota_1$.}
\label{fig:FCstageia}
\end{figure}

On the other hand, suppose for each $j \in \Z_2$ that $\Image(\xi_j)$
intersects either $\iota_1$ or~$C_j$, as in Figure~\ref{fig:FCstageia}(b).
Then as in Case~1, the convexity
of the metric gives that the endpoints $x_0$ and~$x_1$ of $\iota_0$
lie within a distance $R+\delta$ of points $y_0$ and~$y_1$ respectively
on $\Image(\gamma)$ as shown in Figure~\ref{fig:FCstageib}.
The proof of C($i$) now follows from an argument which is
essentially the same as that used in the proof of C($0$) in
Case~1 above.
\begin{figure}[ht!]
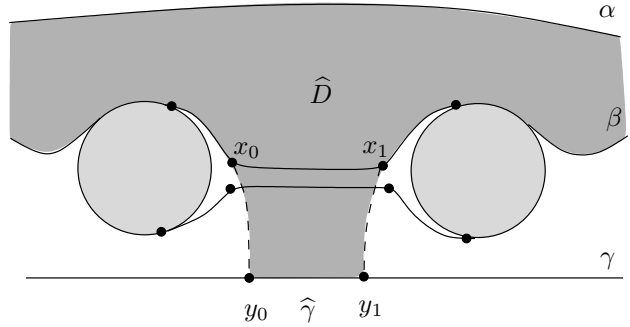

\drawFCstageib
\caption[The subdiagram~$\hat{D}$ is the union of
   $D_\alpha$ and~$\square$.]%
   {The shaded subdiagram~$\hat{D}$ consists of the union of $D_\alpha$
   with the quadrilateral diagram~$\square$ with corners $x_0$,~$x_1$,
   $y_1$, and~$y_0$.}
\label{fig:FCstageib}
\end{figure}

\textbf{Case~3}\qua  Suppose $P_i \cap P_{i+1} = \emptyset$, and the preflats in
question satisfy Condition~(4b) of
Proposition~\ref{prop:DisjointHull}.
As in Case~2, let $H$ be the convex hull of $P_i \cup P_{i+1}$.
Then by our choice of $\delta$ and~$\theta$, the cycle $\boundary H$
is a concatenation
$\omega_0 \eta_0 \upsilon_0 \iota_0 \omega_1 \eta_1 \upsilon_1 \iota_1$
satisfying the conclusion of Proposition~\ref{prop:DisjointHull}
including Condition~(4b), as illustrated in
Figure~\ref{fig:disjoint}(b).
In particular, note that the map $\phi\co D \to X$
sends the preflats $P_i$ and
$P_{i+1}$ to distinct flats $E_i$ and $E_{i+1}$ in~$X$, and each path $\iota_j$
has image a single point~$x_j$.  Conditions B($i$)
and~D($i$) follow as in Case~2 if we set $\sigma_i = \upsilon_0$, \
$\tau_i = \iota_0$, and $\pi_{i+1} = \omega_1$.

Since $\tau_i = \iota_0$ has length zero, the requirement
that $(D_\alpha, \tau_i)$ be ruffled is satisfied vacuously.
Therefore,
in order to prove C($i$), it suffices to show that the distance
in~$X$ between $\phi(x_0)$ and $\Image(\phi\of\gamma)$ is at most $R+M$.
As in Case~2,
for each $j\in \Z_2$, let $\xi_j$ denote a geodesic of length
at most~$R$ connecting $\preflat_{i+j}$ with $\Image(\gamma)$, and
let $C_j$ denote the path $\omega_j \eta_j \upsilon_j$.

\begin{figure}[ht!]
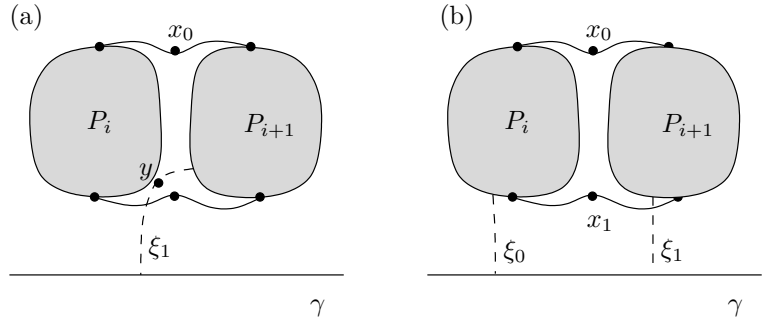

\drawFCstageic
\caption{(a) The path $\xi_1$ intersects $C_0$.\qua
(b) For each $j\in \Z_2$, the path $\xi_j$ intersects $C_j$.}
\label{fig:FCstageic}
\end{figure}
Suppose for some $j\in \Z_2$ that $\Image(\xi_j)$ intersects
$\Image(C_{j+1})$ as shown in Figure~\ref{fig:FCstageic}(a).
Then there is a point $y \in \Image(\xi_j)$ with
\[
  y \in \nbd{\preflat_{i}}{\delta} \cap
        \nbd{\preflat_{i+1}}{\delta},
\]
and $d(y,\Image\gamma) \le R$.
But
\[
  x_0 \in \nbd{\preflat_i}{\delta} \cap \nbd{\preflat_{i+1}}{\delta}.
\]
By our choice of~$M$, the intersection
$\nbd{E_i}{\delta} \cap \nbd{E_{i+1}}{\delta}$
has diameter less than~$M$.  So
$d \bigl( \phi(x_0),\phi(y) \bigr) \le M$, and hence
\[
  d \bigl( \phi(x_0), \Image(\phi\of\gamma) \bigr) \le R + M.
\]
On the other hand, if for each $j \in \Z_2$, the image of~$\xi_j$
intersects $\Image(C_j)$ in some point~$y_j$ as shown in
Figure~\ref{fig:FCstageic}(b), then by convexity of the
metric, it follows that $x_1$ is within a distance~$R$
of $\Image(\gamma)$.
Since $x_0$ and~$x_1$ each lie in
$\nbd{\preflat_{i}}{\delta} \cap \nbd{\preflat_{i+1}}{\delta}$,
it follows that $d \bigl( \phi(x_0), \phi(x_1) \bigr) \le M$,
and hence that $d \bigl( \phi(x_0), \Image(\phi\of\gamma) \bigr) \le R + M$.

\textbf{Case~4}\qua Suppose the intervals $[a_i,b_i]$ and $[a_{i+1},b_{i+1}]$
are disjoint and the corresponding preflats $P_i$ and $P_{i+1}$ intersect.
As in the previous case, the path $\iota_i$ will have image a single point, and
the map $\phi\co D \to X$ sends the preflats to distinct flats in~$X$.
Conditions B($i$), C($i$), and~D($i$)
follow almost exactly as in Case~3.

\textbf{Case~5}\qua The only difference between this case and the previous case is
that the intervals $[a_i,b_i]$ and $[a_{i+1},b_{i+1}]$ intersect.
In this case, an argument similar to the one given in Case~3 shows that
$\phi$ maps their intersection $[a_{i+1},b_i]$ into
$\bignbd{\Image(\phi\of\gamma)}{R+M}$.  The paths $\sigma_i$, \ $\tau_i$, and
$\pi_i$ can be chosen to all have image the same single point, which can be any
point in the intersection $[a_{i+1},b_i]$.
Conditions B($i$), C($i$), and~D($i$)
now follow as before.

This completes the proof of Claim~\ref{claim:FC1}
\end{proof}

Notice that Claim~\ref{claim:FC1}
establishes that $\beta$ satisfies~(\ref{item:FC:main1}).
It also gives a subdivision of~$\beta$ into subpaths
such that $D_\alpha$ is ruffled along
each subpath.  In order to conclude that $D_\alpha$ is ruffled
along~$\beta$, we need to apply Lemma~\ref{lem:ConcatenatedRuffles}
which deals with concatenations of ruffled boundary paths.
In that lemma, one can bound the constants associated to the ruffles
if one has either an upper bound on the number of segments
being concatenated or a lower bound on the lengths of the
concatenated segments.  In Claim~\ref{claim:FC1}, we have control over
neither of these two quantities.

Our strategy for establishing~(\ref{item:FC:main2}) is to
replace $\beta$ with a smoother path $\check{\beta}$ which tracks
close to~$\beta$, using the following claim.

\begin{claim}\label{claim:FC2}
There are universal constants $L(X)$ and $K(X)$ and a path~$\check\beta$ in~$D$
with the same endpoints as~$\beta$ so that
\begin{enumerate}
\renewcommand{\theenumi}{\ensuremath{ \check{\textup{\Roman{enumi}}} }}
\item \label{item:FC:Claim2:z}
   $\check\beta$ is a geodesic in $\check D_\alpha$,
   where $\check D_\alpha$ is the subdiagram of~$D$ bounded by $\check\beta$
   and~$\alpha$,
\item \label{item:FC:Claim2:a}
   the paths $\phi\of\check\beta$ and $\phi\of\gamma$ \  $L$--fellow travel
   relative to some sequence of flats, and
\item \label{item:FC:Claim2:b}
   $(D'_\alpha, \check\beta)$ is $(K,\theta/4)$--ruffled.
\end{enumerate}
\end{claim}

The path~$\check\beta$ will be a concatenation of
\emph{long} paths along each of which $\check D_\alpha$ is ruffled.
The idea is to mimic the construction of~$\beta$, this time using
only those preflats that come close to~$\Image(\gamma)$
and also extend far away from $\Image(\gamma)$.
When some of the preflats involved intersect each other, this modification
is not quite enough to prove
Claim~\ref{claim:FC2}, but it is close in spirit to the modification
we actually use in the proof of that lemma.

\begin{proof}[Proof of Claim~\ref{claim:FC2}]
We use the subdivision of~$\beta $ provided by Claim~\ref{claim:FC1}
as a foundation for our construction of~$\check\beta$.
For each $i$, let us call the path $\pi_i \rho_i \sigma_i$
\emph{nonwandering} if the image of $\phi \of (\pi_i \rho_i \sigma_i)$
in~$X$ lies inside the
$(R+2M)$--neighborhood of $\Image(\phi\of\gamma)$
and call it \emph{wandering} otherwise.
The idea is that wandering paths correspond to ``tall'' preflats
and nonwandering paths correspond to ``short'' preflats.
An interval $[r,s] \subset \Z$ is a \emph{maximal nonwandering interval}
if it is a maximal interval such that for each $i$ in $[r,s]$,
the path $\pi_i \rho_i \sigma_i$ is nonwandering.
The path $\check\beta$ is formed from $\beta$ by cutting out each part
of~$\beta$ corresponding to a maximal nonwandering interval and replacing it
with a path~$\xi$ which is described below.  The path $\xi$ ``smooths'' out
the nonwandering part of~$\beta$.
More precisely, for each maximal nonwandering interval $[r,s]$
replace the subpath
$\sigma_{r-1} \tau_{r-1} \pi_r \dotsm \sigma_{s} \tau_s \pi_{s+1}$
of~$\beta$ with a path $\xi$ constructed as follows.

{\bf Case 1}\qua Suppose the preflats $\preflat_{r-1}$ and $\preflat_{s+1}$
are disjoint.  Let $H$ be the convex hull of
$\preflat_{r-1} \cup \preflat_{s+1}$.
Then $\xi$ is defined to be the subpath of $\boundary H$ connecting
$\rho_{r-1}$ with $\rho_{s+1}$, as shown in Figure~\ref{fig:FCclaimtwoa}.
As in the proof of Cases 2~and~3 from Claim~\ref{claim:FC1},
the path~$\xi$ is a concatenation $\sigma \tau \pi$
so that $\Image(\sigma) \subseteq \nbd{\preflat_{r-1}}{\delta}$,
\ $\Image (\phi\of\tau) \subseteq \bignbd{\Image (\phi\of\gamma)}{R+M}$,
and $\Image(\pi) \subseteq \nbd{\preflat_{s+1}}{\delta}$.
Furthermore, $(D \setminus \interior H, \sigma)$ and
$(D \setminus \interior H, \pi)$ are
$(\delta,\theta)$--ruffled and $(D \setminus \interior H, \tau)$ is
$(\delta',\theta /4)$--ruffled.
\begin{figure}[ht!]
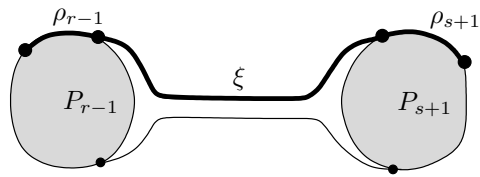

\drawFCclaimtwoa
\caption{The new path $\check\beta$ follows $\xi$ from $\rho_{r-1}$
         to $\rho_{s+1}$.}
\label{fig:FCclaimtwoa}
\end{figure}

{\bf Case 2}\qua Suppose the preflats $\preflat_{r-1}$ and $\preflat_{s+1}$
intersect.  Let $H$ be the convex hull of
$\preflat_{r-1} \cup \preflat_{s+1}$.
Then the boundary path of~$H$ is a concatenation
\[
   \omega_0 \eta_0 \upsilon_0 \dotsm \omega_{2m-1}
   \eta_{2m-1} \upsilon_{2m-1}
\]
satisfying the conclusion of Proposition~\ref{prop:IntersectingHull}.
We may assume that $\rho_{r-1}$ is a subpath of $\eta_0$
and that $\rho_{s+1}$ is a subpath of $\eta_{2i-1}$ for some~$i$.
As shown in Figure~\ref{fig:FCclaimtwob},
let $\sigma$ denote the path which starts at the end of $\rho_{r-1}$
and follows the rest of $\eta_0$ and then follows $\upsilon_0$.
Let $\pi$ denote the path which follows $\omega_1$ and then follows
$\boundary \preflat_{s+1}$ until it reaches the first endpoint of
$\eta_{2i-1}$.  Let $\tau$ be the constant path with image the common
endpoint of $\sigma$ and~$\pi$.
Then $\Image(\sigma) \subseteq \nbd{\preflat_{r-1}}{\delta}$,
\ $\Image (\phi\of\tau) \subseteq \bignbd{\Image (\phi\of\gamma)}{R+M}$,
and $\Image(\pi) \subseteq \nbd{\preflat_{s+1}}{\delta}$.
Furthermore, $(D'_\alpha, \sigma)$ and $(D'_\alpha,\pi)$
are each $(\delta,\theta)$--ruffled, and
$(D'_\alpha, \tau)$ is $(\delta',\theta/4)$--ruffled.
\begin{figure}[ht!]
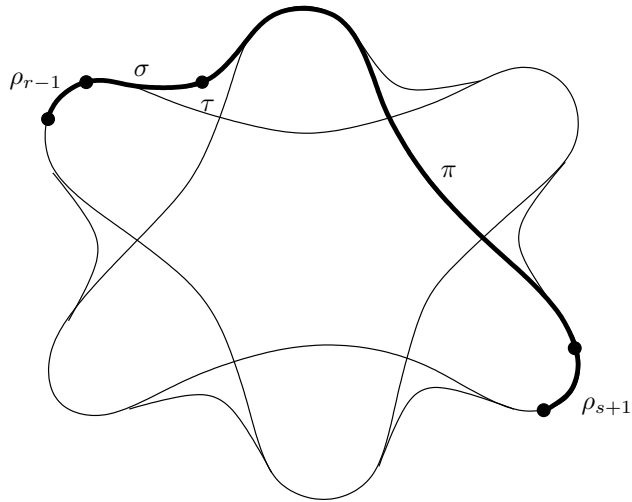

\drawFCclaimtwob
\caption{The new path $\check\beta$ follows $\xi = \sigma \tau \pi$ from
         $\rho_{r-1}$ to $\rho_{s+1}$.}
\label{fig:FCclaimtwob}
\end{figure}

It now follows that in either case, we may relabel the subpaths of~$\check\beta$
and the preflats $\{ \preflat_i \}$ so that $\check\beta$ is a concatenation
\[
   \tau_0 \pi_1 \rho_1 \sigma_1 \tau_1 \dotsm
   \pi_k \rho_k \sigma_k \tau_k
\]
satisfying the following properties
\begin{enumerate}
\item\label{item:FC:Claim2:1}
   $(\hat D_\alpha, \pi_i)$ is $(\delta,\theta)$--ruffled
   and $\Image(\phi\of\pi_1) \subset \nbd{E_i}{\delta}$,
\item\label{item:FC:Claim2:2}
   $(\hat D_\alpha, \rho_i)$ is $(\delta,\theta)$--ruffled
   and $\Image(\phi\of\rho_i) \subset \nbd{E_i}{\delta}$,
\item\label{item:FC:Claim2:3}
   $(\hat D_\alpha, \sigma_i)$ is $(\delta,\theta)$--ruffled
   and $\Image(\phi\of\sigma_i) \subset \nbd{E_i}{\delta}$,
\item\label{item:FC:Claim2:4}
   $(\hat D_\alpha, \tau_i)$ is $(\delta', \theta / 4)$--ruffled
   and $\Image(\phi\of\tau_i) \subseteq \bignbd{\Image(\phi\of\gamma)}{R+M}$,
   and
\item\label{item:FC:Claim2:5}
   $(\pi_i \rho_i \sigma_i)$ has endpoints in
   $\bignbd{\Image(\phi\of\gamma)}{R+M}$
   but its image does not lie entirely inside
   $\bignbd{\Image(\phi\of\gamma)}{R+2M}$.
\end{enumerate}

Clearly the paths $\phi\of\check\beta$ and $\phi\of\gamma$ \  $L$--fellow travel
relative to the sequence $(E_1, \dots, E_k)$ of flats, where $L=R+M$.
Since the pairs $(D'_\alpha,\pi_i)$, $(D'_\alpha,\rho_i)$,
and $(D'_\alpha,\sigma_i)$ are each $(\delta,\theta)$--ruffled,
it follows from Lemma~\ref{lem:ConcatenatedRuffles}
that $(D'_\alpha, \pi_i \rho_i \sigma_i)$ is $(6\delta,\theta)$--ruffled.
Furthermore, by~(\ref{item:FC:Claim2:5}) it is clear that
$(\pi_i \rho_i \sigma_i)$ has length at least $2M$,
which is at least $12\delta$ by our choice of~$M$.
A computation similar to the proof of Lemma~\ref{lem:ConcatenatedRuffles}
now easily shows that $(D'_\alpha,\check\beta)$ is $(K, \theta/4)$--ruffled,
where $K=2M+2\delta$, completing the proof of Claim~\ref{claim:FC2}.
\end{proof}

Replacing $\beta$ with $\check\beta$ and $\theta$ with $\theta/4$
completes the proof of Proposition~\ref{prop:FlatClosure}.
\end{proof}

%%%%%%%%%%%%%%%%%%%%%%%%%%%%%%%%%%%%%%%%%%%%%%%%%%%%%%%%%%%%%%%%%%%%%%%%%%
\section{$2$--complexes with isolated flats have the Relative Fellow Traveller
Property}
\label{sec:RelativeFTP}
%%%%%%%%%%%%%%%%%%%%%%%%%%%%%%%%%%%%%%%%%%%%%%%%%%%%%%%%%%%%%%%%%%%%%%%%%%

In this section, we are finally ready to prove the following theorem.

\begin{thm}\label{thm:RelativeFTP}
A proper, cocompact piecewise Euclidean $\CAT(0)$ $2$--complex
with the Isolated Flats Property
satisfies the Relative Fellow Traveller Property.
\end{thm}

Recall that the Relative Fellow Traveller Property deals with the fellow
travelling of a pair of quasigeodesics with common endpoints.  The proof uses
the following special case in which we consider a geodesic and quasigeodesic
with common endpoints.  This special case is an easy consequence of
Propositions \ref{prop:RuffledFTP} and~\ref{prop:FlatClosure}.

\begin{prop}\label{prop:G-QG}
Let $X$ be a proper, cocompact $\CAT(0)$ space
with the Isolated Flats Property.
Then for each fixed $\lambda$ and $\epsilon$ there is a
constant $R(\lambda,\epsilon,X)$
so that any geodesic and $(\lambda,\epsilon)$--quasigeodesic in~$X$ with
common endpoints $R$--fellow travel relative to flats.
\end{prop}

The proof of Proposition~\ref{prop:G-QG} uses the following standard technical
result, which allows one to ignore the local pathologies of a quasigeodesic
by approximating it with a piecewise geodesic path.
For a proof of this lemma, see \cite[Lemma~III.H.1.11]{BH99}.

\begin{lem}[Taming quasigeodesics]\label{lem:Taming}
Let $X$ be a geodesic space.  Given any $(\lambda,\epsilon)$--quasigeodesic
$\alpha$
in~$X$, one can find a continuous path~$\alpha'$
satisfying the following properties:
\begin{enumerate}
\item $\alpha$ and $\alpha'$ have the same endpoints,
\item $\alpha'$ is piecewise geodesic,
\item $\alpha'$ is a $(\lambda',\epsilon')$--quasigeodesic
      when parametrized by arclength, where $\lambda'$ and~$\epsilon'$
      depend only on $\lambda$ and~$\epsilon$, and
\item the Hausdorff distance between the images of $\alpha$ and~$\alpha'$
      is less than $\lambda + \epsilon$.
\end{enumerate}
\end{lem}

\begin{proof}[Proof of Proposition~\ref{prop:G-QG}]
Let $\alpha$ be a $(\lambda,\epsilon)$--quasigeodesic in~$X$, and
let $\gamma$
be the geodesic connecting its endpoints.
Let $\alpha'$ be a tame $(\lambda',\epsilon')$--quasigeodesic as in
Lemma~\ref{lem:Taming}.
Let $\gamma$ be the geodesic connecting the endpoints of~$\alpha'$.
Notice that $\alpha'$ and~$\gamma'$ both lie in the $1$--skeleton
of some subdivision~$\hat{X}$ of~$X$.
Since the concatenation $\alpha' \bar{\gamma}$ is a nullhomotopic loop
in~$X$,
Theorem~\ref{thm:ReducedDiagram} gives a reduced disc diagram
$\phi\co D \to \hat{X}$ for this loop.
Let $\tilde{\alpha}'$ and $\tilde{\gamma}$ be paths in~$D$
so that $\phi\of \tilde{\alpha}' = \alpha'$ and
$\phi\of \tilde{\gamma} = \gamma$.  Since $\phi$ is distance nonincreasing,
it follows that $\tilde{\gamma}$ is a geodesic and
$\tilde{\alpha}'$ is a $(\lambda',\epsilon')$--quasigeodesic
parametrized by arclength.

Let $\beta$ be a path in~$D$ with the same endpoints as $\tilde{\alpha}'$
and $\tilde{\gamma}$ satisfying the conclusion of
Proposition~\ref{prop:FlatClosure}.
Then applying Proposition~\ref{prop:RuffledFTP} to the subdiagram $D_\alpha$
bounded by $\beta$ and $\tilde{\alpha}'$, we see that
the Hausdorff distance between the images of $\beta$ and $\tilde{\alpha}'$
is at most $N=N(\lambda,\epsilon, X)$.
Since $\phi\of\beta$ and~$\gamma$ \  $L$--fellow travel relative to flats
for some $L$ depending only on~$X$,
it follows that $\gamma$ and~$\alpha$ \  $R$--fellow travel relative
to flats,
where $R$ depends only on $\lambda$, $\epsilon$, and~$X$ as desired.
\end{proof}

We are now ready to prove Theorem~\ref{thm:RelativeFTP} using the special
case proved in Proposition~\ref{prop:G-QG} and the Isolated
Flats Property.

\begin{proof}[Proof of Theorem~\ref{thm:RelativeFTP}]
We need to show that given constants $\lambda$ and~$\epsilon$
we can find a constant~$L$ so that any pair of
$(\lambda,\epsilon)$--quasigeodesics $\alpha$ and~$\alpha'$ with common
endpoints $L$--fellow
travel relative to flats.
So fix a pair $\alpha$ and $\alpha'$ of such quasigeodesics,
and let $\gamma$ denote the geodesic connecting the endpoints
of $\alpha$.

By Proposition~\ref{prop:G-QG}, we know that the paths $\alpha$ and~$\gamma$
and the paths $\alpha'$ and~$\gamma$ each $R$--fellow travel
relative to some sequence of flats, where $R=R(\lambda,\epsilon,X)$.
The main difficulty with pasting together
these two facts is that the sequences of flats involved
may not be the same.  In fact, the sequences of flats for each pair
are not even well-defined.  For instance, if a pair of paths
travels for a sufficiently short distance
in some flat, then that flat can be inserted or deleted from a
sequence of flats without affecting whether the paths fellow travel relative
to the sequence.

To circumvent the difficulties alluded to above, we construct a canonical
sequence of flats using only properties of~$\gamma$
so that each of the given pairs fellow travels relative to the constructed
sequence.  Furthermore, the partition of~$\gamma$ in
Definition~\ref{def:RelativeFellowTravelling} will be canonical,
so that the two pairs can be pasted together coherently.  It will then
follow
that $\alpha$ and~$\alpha'$ fellow travel relative to the canonical sequence
of flats.

By the Isolated Flats Property,
there is a constant $K=K(R,X)$ such that given any pair of flat planes
$E_1$ and~$E_2$ in~$X$ the intersection $\nbd{E_1}{R} \cap \nbd{E_2}{R}$
has diameter less than~$K$.
Let us call a flat plane~$E$ in~$X$ an \emph{essential $\gamma$--flat}
if $\Image(\gamma)$ intersects the $R$--neighborhood of~$E$
in a segment~$\xi$ of length at least~$K$.
In this case, the segment~$\xi$ will be called a
\emph{maximal flat segment} of~$\gamma$.
Removing the subsegment of length $K/2$ from each end of~$\xi$
gives a shorter segment~$\eta$ called
a \emph{shortened flat segment} of~$\gamma$.
Notice that any two distinct shortened flat segments $\eta_1$ and~$\eta_2$
of~$\gamma$ are disjoint, since otherwise there would be a pair of distinct
flats
$E_1$ and $E_2$ whose $R$--neighborhoods intersect in a set of diameter
at least~$K$.

Our canonical sequence of flats will be the essential
$\gamma$--flats, and the canonical partition of~$\gamma$ will
consist of the shortened flat segments of~$\gamma$
alternating with the segments of~$\gamma$ that connect two consecutive
shortened flat segments.
We need to verify that $\gamma$ and $\alpha$ actually fellow travel
relative to this canonical data using the given
(non-canonical)
sequence of flats and our given (non-canonical) partition of~$\gamma$

Our argument consists of two directions.  First we show that each flat
of our given sequence appears in the canonical sequence, unless it is
very small.  Then we show that every flat of the canonical sequence
which does not correspond to a flat of our given sequence
may be added to that sequence without creating problems.

For the first direction,
suppose the quasigeodesic~$\alpha$ contains a subpath~$\pi$ whose
endpoints $x$ and~$y$ lie within a distance~$R$ of points $w$ and~$z$
in the image of~$\gamma$, and suppose further that the image of~$\pi$
lies in an $R$--neighborhood of some flat~$E$.
Finally suppose the segment $[w,z]$ lies in an $R$--neighborhood of~$E$.

If the distance from $w$ to~$z$ is less than~$K$,
then the distance from $x$ to~$y$ is less than $K+2R$.
An easy computation using the definition of quasigeodesic
then shows that $\Image(\pi)$ has diameter less than
$\lambda^2 (K+2R) + 2\epsilon\lambda$.  So in this case, $\Image(\pi)$
lies inside a $R_1$--neighborhood of $[w,z]$, where
$R_1=R+\lambda^2 (K+2R) +2\epsilon\lambda$.  Similarly, since the diameter
of $[w,z]$ is less than~$K$, it follows that $[w,z]$ lies in a
$R_2$--neighborhood of $\Image(\pi)$, where $R_2=R+K$.  So the Hausdorff
distance between $\Image(\pi)$ and $[w,z]$ is less than
$R'=\max \{R_1,R_2\}$.
On the other hand, if $d(w,z) \ge K$, then $[w,z]$ lies inside
a unique maximal flat segment of~$\gamma$.

Now for the second direction,
suppose $\alpha$ contains a subpath~$\beta$
with endpoints $x'$ and~$y'$ such that $\beta$ is within
a Hausdorff distance~$R$ from a subpath $[w',z']$ of~$\gamma$.
Suppose further that $[w',z']$ lies in an $R$--neighborhood of some flat
plane~$E$.  Then $\Image(\beta)$ lies in a $2R$--neighborhood of~$E$.

It is now clear that we can choose partitions
\[
   0 = t_0 \le s_0 \le t_1 \le s_1 \le \dots \le t_n \le s_n = a
\]
and
\[
   0 =t'_0 \le s'_0\le t'_1\le s'_1\le \dots \le t'_n\le s'_n = a'
\]
of the domains of $\gamma$ and~$\alpha$ respectively,
so that the subpaths $\gamma \big| [s_{i-1},t_{i}]$
are precisely the collection of all shortened flat segments of~$\gamma$.
Furthermore, the Hausdorff distance between
$\gamma \bigl( [t_i,s_i] \bigr)$
and $\alpha \bigl( [t'_i,s'_i] \bigr)$ is at most $R''(\lambda,\epsilon,X)$,
while the sets $\gamma \bigl( [t_{i-1},s_i] \bigr)$ and
$\alpha \bigl( [t'_{i-1},s'_i] \bigr)$ lie in a $2R$--neighborhood of a unique
essential $\gamma$--flat $E_i$.

Since the partition of~$\gamma$ and the sequence of essential $\gamma$--flats
described above is independent of the choice of~$\alpha$,
it follows that $\alpha$ and $\alpha'$ \  $L$--fellow travel
relative to the sequence of flats $(E_1, \dots, E_n)$
for some constant $L$ depending only on $\lambda$, $\epsilon$, and~$X$.
\end{proof}

%%%%%%%%%%%%%%%%%%%%%%%%%%%%%%%%%%%%%%%%%%%%%%%%%%%%%%%%%%%%%%%%%%%%%%%%
%%                  BIBLIOGRAPHY
%%%%%%%%%%%%%%%%%%%%%%%%%%%%%%%%%%%%%%%%%%%%%%%%%%%%%%%%%%%%%%%%%%%%%%%%


\begin{thebibliography}

\bibitem{Aitchison}
\textbf{I Aitchison}, \emph{Canonical flat structures on $3$-manifolds},
  preprint

\bibitem{ABC91}
\textbf{J\,M Alonso}, \textbf{T Brady}, \textbf{D Cooper}, \textbf{V Ferlini},
  \textbf{M Lustig}, \textbf{M Mihalik}, \textbf{H Short}, \emph{Notes on word
  hyperbolic groups}, (H Short, editor),
  from: ``Group theory from a geometrical viewpoint
  \textup{(}Trieste, 1990\textup{)}'', ({\'E} Ghys, A Haefliger, A Verjovsky,
  editors), World Sci.\ Publishing, River Edge, NJ (1991)  3--63

\bibitem{Ballmann90}
\textbf{W Ballmann}, \emph{Singular spaces of nonpositive curvature}, from:
  ``Sur les groupes hyperboliques d'apr\`es Mikhael Gromov \textup{(}Bern,
  1988\textup{)}'', ({\'E} Ghys, P\,de~la Harpe, editors), Birkh\"auser Boston,
  Boston, MA (1990)  189--201

\bibitem{BallmannBrin94}
\textbf{W Ballmann}, \textbf{M Brin}, \emph{Polygonal complexes and
  combinatorial group theory}, Geom. Dedicata 50 (1994) 165--191

\bibitem{BangertSchroeder91}
\textbf{V Bangert}, \textbf{V Schroeder}, \emph{Existence of flat tori in
  analytic manifolds of nonpositive curvature}, Ann. Sci. \'Ecole Norm. Sup.
  (4) 24 (1991) 605--634

\bibitem{Benakli94}
\textbf{N Benakli}, \emph{Polygonal complexes. {I}. {C}ombinatorial and
  geometric properties}, J. Pure Appl. Algebra 97 (1994) 247--263

\bibitem{BowditchRelHyp}
\textbf{B\,H Bowditch}, \emph{Relatively hyperbolic groups} (1999), preprint,
  University of Southampton

\bibitem{Bridson91}
\textbf{M\,R Bridson}, \emph{Geodesics and curvature in metric simplicial
  complexes}, from: ``Group theory from a geometrical viewpoint
  \textup{(}{T}rieste, 1990\textup{)}'', ({\'E} Ghys, A Haefliger, A Verjovsky,
  editors), World Sci.\ Publishing, River Edge, NJ (1991)  373--463

\bibitem{Bridson95}
\textbf{M\,R Bridson}, \emph{On the existence of flat planes in spaces of
  nonpositive curvature}, Proc.\ Amer.\ Math.\ Soc. 123 (1995) 223--235

\bibitem{BH99}
\textbf{M\,R Bridson}, \textbf{A Haefliger}, \emph{Metric spaces of
  non-positive curvature}, Springer-Verlag, Berlin (1999)

\bibitem{CrokeKleiner00}
\textbf{C\,B Croke}, \textbf{B Kleiner}, \emph{Spaces with nonpositive
  curvature and their ideal boundaries}, Topology 39 (2000) 549--556

\bibitem{CrokeKleiner02}
\textbf{C\,B Croke}, \textbf{B Kleiner}, \emph{The geodesic flow of a
  nonpositively curved graph manifold}, Geom. Funct. Anal. 12 (2002) 479--545

\bibitem{Dehn87}
\textbf{M Dehn}, \emph{Papers on group theory and topology}, Springer-Verlag,
  New York (1987), translated from the German by J.~Stillwell

\bibitem{Eberlein73}
\textbf{P Eberlein}, \emph{Geodesic flows on negatively curved manifolds.
  {I}{I}}, Trans.\ Amer.\ Math.\ Soc. 178 (1973) 57--82

\bibitem{EfromovichTihomirova63}
\textbf{V\,A Efromovich}, \textbf{E\,S Tihomirova}, \emph{Continuation of an
  equimorphism to infinity}, Soviet Math. Dokl. 4 (1963) 1494--1496

\bibitem{ECHLPT92}
\textbf{D\,B\,A Epstein}, \textbf{J\,W Cannon}, \textbf{D\,F Holt},
  \textbf{S\,V\,F Levy}, \textbf{M\,S Paterson}, \textbf{W\,P Thurston},
  \emph{Word processing in groups}, Jones and Bartlett Publishers, Boston, MA
  (1992)

\bibitem{Farb98}
\textbf{B Farb}, \emph{Relatively hyperbolic groups}, Geom.\ Funct.\ Anal. 8
  (1998) 810--840

\bibitem{Gersten87}
\textbf{S\,M Gersten}, \emph{Reducible diagrams and equations over groups},
  from: ``Essays in group theory'', (S\,M Gersten, editor), Springer, New York
  (1987)  15--73

\bibitem{GerstenShort90Automatic}
\textbf{S\,M Gersten}, \textbf{H\,B Short}, \emph{Small cancellation theory and
  automatic groups}, Invent. Math. 102 (1990) 305--334

\bibitem{GerstenShort91Automatic}
\textbf{S\,M Gersten}, \textbf{H\,B Short}, \emph{Small cancellation theory and
  automatic groups. {II}}, Invent. Math. 105 (1991) 641--662

\bibitem{Gromov87}
\textbf{M Gromov}, \emph{Hyperbolic groups}, from: ``Essays in group theory'',
  (S\,M Gersten, editor), Springer, New York (1987)  75--263

\bibitem{Haglund91}
\textbf{F Haglund}, \emph{Les poly\`edres de {G}romov}, C. R. Acad. Sci. Paris
  S\'er. I Math. 313 (1991) 603--606

\bibitem{Heber87}
\textbf{J Heber}, \emph{Hyperbolische geodatische {R}aume}, Diplomarbeit,
  Univ.\ Bonn (1987)

\bibitem{HruskaGeometric}
\textbf{G\,C Hruska}, \emph{Geometric invariants of spaces with isolated
  flats}, preprint available at
  \texttt{http://www.math.uchicago.edu/{\~{}}chruska/papers}

\bibitem{HruskaRelHyp}
\textbf{G\,C Hruska}, \emph{On the relative hyperbolicity of nonpositively
  curved groups with isolated flats}, in preparation

\bibitem{KapovichLeeb95}
\textbf{M Kapovich}, \textbf{B Leeb}, \emph{On asymptotic cones and
  quasi-isometry classes of fundamental groups of $3$-manifolds}, Geom.\
  Funct.\ Anal. 5 (1995) 582--603

\bibitem{KariPapasoglu99}
\textbf{J Kari}, \textbf{P Papasoglu}, \emph{Deterministic aperiodic tile
  sets}, Geom. Funct. Anal. 9 (1999) 353--369

\bibitem{Lang96}
\textbf{U Lang}, \emph{Quasigeodesics outside horoballs}, Geom.\ Dedicata 63
  (1996) 205--215

\bibitem{Lyndon66}
\textbf{R\,C Lyndon}, \emph{On {D}ehn's algorithm}, Math.\ Ann. 166 (1966)
  208--228

\bibitem{LS77}
\textbf{R\,C Lyndon}, \textbf{P\,E Schupp}, \emph{Combinatorial group theory},
  Ergebnisse der Mathematik und ihrer Grenzgebiete, Band 89, Springer-Verlag,
  Berlin (1977)

\bibitem{MasurMinsky99}
\textbf{H\,A Masur}, \textbf{Y\,N Minsky}, \emph{Geometry of the complex of
  curves. {I}. {H}yperbolicity}, Invent. Math. 138 (1999) 103--149

\bibitem{McCammondWise02}
\textbf{J\,P McCammond}, \textbf{D\,T Wise}, \emph{Fans and ladders in small
  cancellation theory}, Proc.\ London Math.\ Soc. \textup{(}3\textup{)} 84
  (2002) 599--644

\bibitem{Morse24}
\textbf{H\,M Morse}, \emph{A fundamental class of geodesics on any closed
  surface of genus greater than one}, Trans.\ Amer.\ Math.\ Soc. 26 (1924)
  25--60

\bibitem{Moussong88}
\textbf{G Moussong}, \emph{Hyperbolic {C}oxeter groups}, PhD thesis, Ohio
  State Univ. (1988)

\bibitem{NibloReeves98}
\textbf{G\,A Niblo}, \textbf{L\,D Reeves}, \emph{The geometry of cube complexes
  and the complexity of their fundamental groups}, Topology 37 (1998) 621--633

\bibitem{Papasoglu95}
\textbf{P Papasoglu}, \emph{Strongly geodesically automatic groups are
  hyperbolic}, Invent. Math. 121 (1995) 323--334

\bibitem{Pride88}
\textbf{S\,J Pride}, \emph{Star-complexes, and the dependence problems for
  hyperbolic complexes}, Glasgow Math.\ J. 30 (1988) 155--170

\bibitem{Rebbechi01}
\textbf{D\,Y Rebbechi}, \emph{Algorithmic properties of relatively hyperbolic
  groups}, PhD thesis, Rutgers Univ. (2001)

\bibitem{Rudin}
\textbf{W Rudin}, \emph{Real and complex analysis}, third edition, McGraw-Hill
  Book Co., New York (1987)

\bibitem{SelaProblems}
\textbf{Z Sela}, \emph{Diophantine geometry over groups: {A} list of research
  problems}, available at {\texttt{http://www.ma.huji.ac.il/{\~{}}zlil/}}

\bibitem{Short91}
\textbf{H Short}, \emph{Quasiconvexity and a theorem of {H}owson's}, from:
  ``Group theory from a geometrical viewpoint \textup{(}{T}rieste,
  1990\textup{)}'', ({\'E} Ghys, A Haefliger, A Verjovsky, editors), World Sci.
  Publishing, River Edge, NJ (1991)  168--176

\bibitem{Tukia94}
\textbf{P Tukia}, \emph{Convergence groups and {G}romov's metric hyperbolic
  spaces}, New Zealand J. Math. 23 (1994) 157--187

\bibitem{vanKampen33}
\textbf{E\,R van Kampen}, \emph{On some lemmas in the theory of groups}, Amer.\
  J. Math. 55 (1933) 268--273

\bibitem{Weinbaum71}
\textbf{C\,M Weinbaum}, \emph{The word and conjugacy problems for the knot
  group of any tame, prime, alternating knot}, Proc. Amer. Math. Soc. 30 (1971)
  22--26

\bibitem{WilsonBoundary}
\textbf{J Wilson}, \emph{A \textup{CAT(0)} group with uncountably many distinct
  boundaries}, preprint

\bibitem{Wise96}
\textbf{D\,T Wise}, \emph{Non-positively curved squared complexes, aperiodic
  til\-ings, and non-residually finite groups}, PhD thesis, Princeton Univ.
  (1996)

\bibitem{WiseFigure8}
\textbf{D\,T Wise}, \emph{Subgroup separability of the figure 8 knot group}
  (1998), pre\-print available at \texttt{http://www.math.mcgill.ca/%
  \linebreak[0]wise/\linebreak[0]papers.\linebreak[0]html}

\end{thebibliography}
\end{document}